\documentclass[12pt,reqno]{amsart}

\usepackage{amssymb}

\input texdraw
\catcode`\@=11

\thinlines
\newskip\Einheit \Einheit=0.6cm
\newcount\xcoord \newcount\ycoord
\newdimen\xdim \newdimen\ydim \newdimen\PfadD@cke \newdimen\Pfadd@cke
\def\PfadDicke#1{\PfadD@cke#1 \divide\PfadD@cke by2 \Pfadd@cke\PfadD@cke \multiply\PfadD@cke by2}
\long\def\LOOP#1\REPEAT{\def\BODY{#1}\ITERATE}
\def\ITERATE{\BODY \let\next\ITERATE \else\let\next\relax\fi \next}
\let\REPEAT=\fi
\def\Punkt{\hbox{\raise-2pt\hbox to0pt{\hss\scriptsize$\bullet$\hss}}}
\def\DuennPunkt(#1,#2){\unskip
  \raise#2 \Einheit\hbox to0pt{\hskip#1 \Einheit
          \raise-2.5pt\hbox to0pt{\hss\normalsize$\bullet$\hss}\hss}}
\def\NormalPunkt(#1,#2){\unskip
  \raise#2 \Einheit\hbox to0pt{\hskip#1 \Einheit
          \raise-3pt\hbox to0pt{\hss\large$\bullet$\hss}\hss}}
\def\DickPunkt(#1,#2){\unskip
  \raise#2 \Einheit\hbox to0pt{\hskip#1 \Einheit
          \raise-4pt\hbox to0pt{\hss\Large$\bullet$\hss}\hss}}
\def\Kreis(#1,#2){\unskip
  \raise#2 \Einheit\hbox to0pt{\hskip#1 \Einheit
          \raise-4pt\hbox to0pt{\hss\Large$\circ$\hss}\hss}}
\def\Diagonale(#1,#2)#3{\unskip\leavevmode
  \xcoord#1\relax \ycoord#2\relax
      \raise\ycoord \Einheit\hbox to0pt{\hskip\xcoord \Einheit
         \unitlength\Einheit
         \line(1,1){#3}\hss}}
\def\AntiDiagonale(#1,#2)#3{\unskip\leavevmode
  \xcoord#1\relax \ycoord#2\relax 
      \raise\ycoord \Einheit\hbox to0pt{\hskip\xcoord \Einheit
         \unitlength\Einheit
         \line(1,-1){#3}\hss}}
\def\Pfad(#1,#2),#3\endPfad{\unskip\leavevmode
  \xcoord#1 \ycoord#2 \thicklines\ZeichnePfad#3\endPfad\thinlines}
\def\ZeichnePfad#1{\ifx#1\endPfad\let\next\relax
  \else\let\next\ZeichnePfad
    \ifnum#1=1
      \raise\ycoord \Einheit\hbox to0pt{\hskip\xcoord \Einheit
         \vrule height\Pfadd@cke width1 \Einheit depth\Pfadd@cke\hss}%
      \advance\xcoord by 1
    \else\ifnum#1=2
      \raise\ycoord \Einheit\hbox to0pt{\hskip\xcoord \Einheit
        \hbox{\hskip-\PfadD@cke\vrule height1 \Einheit width\PfadD@cke depth0pt}\hss}%
      \advance\ycoord by 1
    \else\ifnum#1=3
      \raise\ycoord \Einheit\hbox to0pt{\hskip\xcoord \Einheit
         \unitlength\Einheit
         \line(1,1){1}\hss}
      \advance\xcoord by 1
      \advance\ycoord by 1
    \else\ifnum#1=4
      \raise\ycoord \Einheit\hbox to0pt{\hskip\xcoord \Einheit
         \unitlength\Einheit
         \line(1,-1){1}\hss}
      \advance\xcoord by 1
      \advance\ycoord by -1
    \else\ifnum#1=5
      \advance\xcoord by -1
      \raise\ycoord \Einheit\hbox to0pt{\hskip\xcoord \Einheit
         \vrule height\Pfadd@cke width1 \Einheit depth\Pfadd@cke\hss}%
    \else\ifnum#1=6
      \advance\ycoord by -1
      \raise\ycoord \Einheit\hbox to0pt{\hskip\xcoord \Einheit
        \hbox{\hskip-\PfadD@cke\vrule height1 \Einheit width\PfadD@cke depth0pt}\hss}%
    \else\ifnum#1=7
      \advance\xcoord by -1
      \advance\ycoord by -1
      \raise\ycoord \Einheit\hbox to0pt{\hskip\xcoord \Einheit
         \unitlength\Einheit
         \line(1,1){1}\hss}
    \else\ifnum#1=8
      \advance\xcoord by -1
      \advance\ycoord by +1
      \raise\ycoord \Einheit\hbox to0pt{\hskip\xcoord \Einheit
         \unitlength\Einheit
         \line(1,-1){1}\hss}
    \fi\fi\fi\fi
    \fi\fi\fi\fi
  \fi\next}
\def\hSSchritt{\leavevmode\raise-.4pt\hbox to0pt{\hss.\hss}\hskip.2\Einheit
  \raise-.4pt\hbox to0pt{\hss.\hss}\hskip.2\Einheit
  \raise-.4pt\hbox to0pt{\hss.\hss}\hskip.2\Einheit
  \raise-.4pt\hbox to0pt{\hss.\hss}\hskip.2\Einheit
  \raise-.4pt\hbox to0pt{\hss.\hss}\hskip.2\Einheit}
\def\vSSchritt{\vbox{\baselineskip.2\Einheit\lineskiplimit0pt
\hbox{.}\hbox{.}\hbox{.}\hbox{.}\hbox{.}}}
\def\DSSchritt{\leavevmode\raise-.4pt\hbox to0pt{%
  \hbox to0pt{\hss.\hss}\hskip.2\Einheit
  \raise.2\Einheit\hbox to0pt{\hss.\hss}\hskip.2\Einheit
  \raise.4\Einheit\hbox to0pt{\hss.\hss}\hskip.2\Einheit
  \raise.6\Einheit\hbox to0pt{\hss.\hss}\hskip.2\Einheit
  \raise.8\Einheit\hbox to0pt{\hss.\hss}\hss}}
\def\dSSchritt{\leavevmode\raise-.4pt\hbox to0pt{%
  \hbox to0pt{\hss.\hss}\hskip.2\Einheit
  \raise-.2\Einheit\hbox to0pt{\hss.\hss}\hskip.2\Einheit
  \raise-.4\Einheit\hbox to0pt{\hss.\hss}\hskip.2\Einheit
  \raise-.6\Einheit\hbox to0pt{\hss.\hss}\hskip.2\Einheit
  \raise-.8\Einheit\hbox to0pt{\hss.\hss}\hss}}
\def\SPfad(#1,#2),#3\endSPfad{\unskip\leavevmode
  \xcoord#1 \ycoord#2 \ZeichneSPfad#3\endSPfad}
\def\ZeichneSPfad#1{\ifx#1\endSPfad\let\next\relax
  \else\let\next\ZeichneSPfad
    \ifnum#1=1
      \raise\ycoord \Einheit\hbox to0pt{\hskip\xcoord \Einheit
         \hSSchritt\hss}%
      \advance\xcoord by 1
    \else\ifnum#1=2
      \raise\ycoord \Einheit\hbox to0pt{\hskip\xcoord \Einheit
        \hbox{\hskip-2pt \vSSchritt}\hss}%
      \advance\ycoord by 1
    \else\ifnum#1=3
      \raise\ycoord \Einheit\hbox to0pt{\hskip\xcoord \Einheit
         \DSSchritt\hss}
      \advance\xcoord by 1
      \advance\ycoord by 1
    \else\ifnum#1=4
      \raise\ycoord \Einheit\hbox to0pt{\hskip\xcoord \Einheit
         \dSSchritt\hss}
      \advance\xcoord by 1
      \advance\ycoord by -1
    \else\ifnum#1=5
      \advance\xcoord by -1
      \raise\ycoord \Einheit\hbox to0pt{\hskip\xcoord \Einheit
         \hSSchritt\hss}%
    \else\ifnum#1=6
      \advance\ycoord by -1
      \raise\ycoord \Einheit\hbox to0pt{\hskip\xcoord \Einheit
        \hbox{\hskip-2pt \vSSchritt}\hss}%
    \else\ifnum#1=7
      \advance\xcoord by -1
      \advance\ycoord by -1
      \raise\ycoord \Einheit\hbox to0pt{\hskip\xcoord \Einheit
         \DSSchritt\hss}
    \else\ifnum#1=8
      \advance\xcoord by -1
      \advance\ycoord by 1
      \raise\ycoord \Einheit\hbox to0pt{\hskip\xcoord \Einheit
         \dSSchritt\hss}
    \fi\fi\fi\fi
    \fi\fi\fi\fi
  \fi\next}
\def\Koordinatenachsen(#1,#2){\unskip
 \hbox to0pt{\hskip-.5pt\vrule height#2 \Einheit width.5pt depth1 \Einheit}%
 \hbox to0pt{\hskip-1 \Einheit \xcoord#1 \advance\xcoord by1
    \vrule height0.25pt width\xcoord \Einheit depth0.25pt\hss}}
\def\Koordinatenachsen(#1,#2)(#3,#4){\unskip
 \hbox to0pt{\hskip-.5pt \ycoord-#4 \advance\ycoord by1
    \vrule height#2 \Einheit width.5pt depth\ycoord \Einheit}%
 \hbox to0pt{\hskip-1 \Einheit \hskip#3\Einheit 
    \xcoord#1 \advance\xcoord by1 \advance\xcoord by-#3 
    \vrule height0.25pt width\xcoord \Einheit depth0.25pt\hss}}
\def\Gitter(#1,#2){\unskip \xcoord0 \ycoord0 \leavevmode
  \LOOP\ifnum\ycoord<#2
    \loop\ifnum\xcoord<#1
      \raise\ycoord \Einheit\hbox to0pt{\hskip\xcoord \Einheit\Punkt\hss}%
      \advance\xcoord by1
    \repeat
    \xcoord0
    \advance\ycoord by1
  \REPEAT}
\def\Gitter(#1,#2)(#3,#4){\unskip \xcoord#3 \ycoord#4 \leavevmode
  \LOOP\ifnum\ycoord<#2
    \loop\ifnum\xcoord<#1
      \raise\ycoord \Einheit\hbox to0pt{\hskip\xcoord \Einheit\Punkt\hss}%
      \advance\xcoord by1
    \repeat
    \xcoord#3
    \advance\ycoord by1
  \REPEAT}
\def\Label#1#2(#3,#4){\unskip \xdim#3 \Einheit \ydim#4 \Einheit
  \def\lo{\advance\xdim by-.5 \Einheit \advance\ydim by.5 \Einheit}%
  \def\llo{\advance\xdim by-.25cm \advance\ydim by.5 \Einheit}%
  \def\loo{\advance\xdim by-.5 \Einheit \advance\ydim by.25cm}%
  \def\o{\advance\ydim by.25cm}%
  \def\ro{\advance\xdim by.5 \Einheit \advance\ydim by.5 \Einheit}%
  \def\rro{\advance\xdim by.25cm \advance\ydim by.5 \Einheit}%
  \def\roo{\advance\xdim by.5 \Einheit \advance\ydim by.25cm}%
  \def\l{\advance\xdim by-.30cm}%
  \def\r{\advance\xdim by.30cm}%
  \def\lu{\advance\xdim by-.5 \Einheit \advance\ydim by-.6 \Einheit}%
  \def\llu{\advance\xdim by-.25cm \advance\ydim by-.6 \Einheit}%
  \def\luu{\advance\xdim by-.5 \Einheit \advance\ydim by-.30cm}%
  \def\u{\advance\ydim by-.30cm}%
  \def\ru{\advance\xdim by.5 \Einheit \advance\ydim by-.6 \Einheit}%
  \def\rru{\advance\xdim by.25cm \advance\ydim by-.6 \Einheit}%
  \def\ruu{\advance\xdim by.5 \Einheit \advance\ydim by-.30cm}%
  #1\raise\ydim\hbox to0pt{\hskip\xdim
     \vbox to0pt{\vss\hbox to0pt{\hss$#2$\hss}\vss}\hss}%
}
\catcode`\@=12

\edef\pKratBN {1}
\edef\pKratBN {\pKratBN , 1}
\edef\pKratBN {\pKratBN , 1}
\edef\pKratBN {\pKratBN , 1}
\edef\pKratBN {\pKratBN , 2}
\edef\pKratBN {\pKratBN , 2}
\edef\pStanBI {2}
\edef\pKratBN {\pKratBN , 2}
\edef\pKratBN {\pKratBN , 2}
\edef\pBresAO {3}
\edef\pKratBN {\pKratBN , 3}
\edef\pZeilAP {3}
\edef\pPeWZAA {3}
\edef\pZeilAP {\pZeilAP , 3}
\edef\pKupeAH {3}
\edef\pCoBWAA {3}
\edef\pKratBN {\pKratBN , 3}
\edef\pKratBN {\pKratBN , 3}
\edef\pRubeAD {3}
\edef\pPropAH {4}
\edef\pMacMAA {5}
\edef\pDT {5}
\edef\pCiEKAA {6}
\edef\pCiEKAA {6}
\edef\pKratAY {6}
\edef\pKratBK {6}
\edef\pCiEKAA {\pCiEKAA , 7}
\edef\pCiKrAA {7}
\edef\pCiKrAC {7}
\edef\pCiKrAD {7}
\edef\pEisTAA {7}
\edef\pEisTAB {7}
\edef\pEisTAF {7}
\edef\pFiscAA {7}
\edef\pKratBY {7}
\edef\pOkKrAA {7}
\edef\pCiEKAA {\pCiEKAA , 7}
\edef\pKupeAG {7}
\edef\pKastAA {7}
\edef\pKupeAG {\pKupeAG , 7}
\edef\pKaMGAB {7}
\edef\pKaMGAC {7}
\edef\pLindAA {7}
\edef\pFishAA {7}
\edef\pJoSaAB {7}
\edef\pGrJSAA {7}
\edef\pGeViAA {7}
\edef\pGeViAB {7}
\edef\pSlatZY {8}
\edef\pSlatZZ {8}
\edef\pGeViAA {9}
\edef\pGeViAB {9}
\edef\pStemAE {9}
\edef\pCiucAH {9}
\edef\pCiEKAA {\pCiEKAA , 9}
\edef\pCiKrAA {\pCiKrAA , 9}

\edef\pCiKrAD {\pCiKrAD , 9}
\edef\pEisTAA {\pEisTAA , 9}
\edef\pEisTAB {\pEisTAB , 9}
\edef\pEisTAF {\pEisTAF , 9}
\edef\pFiscAA {\pFiscAA , 9}
\edef\pFuKrAC {9}
\edef\pKratBN {\pKratBN , 9}
\edef\pKratBN {\pKratBN , 9}
\edef\pCiEKAA {\pCiEKAA , 9}
\edef\pKratBD {9}
\edef\pCiKrAB {9}
\edef\pKratBN {\pKratBN , 9}
\edef\pGosper {9}
\edef\pBellard {9}
\edef\pPlouffe {9}
\edef\pBellard {\pBellard , 10}
\edef\pBaBoPl {10}
\edef\pAlKPAA {10}
\edef\pLeLLAA {10}
\edef\pCohHAA {10}
\edef\pAAR {12}
\edef\pAlKPAA {\pAlKPAA , 13}
\edef\pAlKPAA {\pAlKPAA , 15}
\edef\pKratBN {\pKratBN , 16}
\edef\pKratBN {\pKratBN , 17}
\edef\pKratBN {\pKratBN , 17}
\edef\pBresAO {\pBresAO , 17}
\edef\pKnutAF {17}
\edef\pAmdeAD {17}
\edef\pKratBN {\pKratBN , 18}
\edef\pKratBN {\pKratBN , 18}
\edef\pKratBD {\pKratBD , 18}
\edef\pKratBI {19}
\edef\pKratBD {\pKratBD , 21}
\edef\pSinVAA {21}
\edef\pPeWZAA {21}
\edef\pZeilAP {21}
\edef\pZeilAM {21}
\edef\pZeilAV {21}
\edef\pKratBD {\pKratBD , 22}
\edef\pChSaAA {22}
\edef\pWiZeAC {22}
\edef\pAlKPAA {\pAlKPAA , 28}
\edef\pKratBN {\pKratBN , 28}
\edef\pKratBN {\pKratBN , 28}
\edef\pKratBN {\pKratBN , 28}
\edef\pKratBN {\pKratBN , 28}
\edef\pKratBN {\pKratBN , 28}
\edef\pKratBN {\pKratBN , 29}
\edef\pKratBN {\pKratBN , 29}
\edef\pAmZeAB {29}
\edef\pKratBN {\pKratBN , 29}
\edef\pAmZeAB {\pAmZeAB , 29}
\edef\pAmZeAB {\pAmZeAB , 29}
\edef\pAmZeAB {\pAmZeAB , 29}
\edef\pKupeAH {30}
\edef\pElKLAB {30}
\edef\pKupeAD {30}
\edef\pKupeAH {\pKupeAH , 30}
\edef\pIzerAA {30}
\edef\pKupeAH {\pKupeAH , 30}
\edef\pMuirAB {30}
\edef\pSchuAA {30}
\edef\pStemAE {\pStemAE , 30}
\edef\pStemAE {\pStemAE , 30}
\edef\pKratBN {\pKratBN , 30}
\edef\pKupeAH {\pKupeAH , 30}
\edef\pLaLTAA {30}
\edef\pStemAE {\pStemAE , 30}
\edef\pKupeAH {\pKupeAH , 30}
\edef\pOkadAJ {32}
\edef\pIsOTAA {32}
\edef\pLascAT {32}
\edef\pRoScAC {32}
\edef\pKratBN {\pKratBN , 32}
\edef\pKratBN {\pKratBN , 33}
\edef\pKratBN {\pKratBN , 33}
\edef\pRoScAC {\pRoScAC , 33}
\edef\pRoScAC {\pRoScAC , 33}
\edef\pSchlAB {33}
\edef\pSchlAF {33}
\edef\pRoScAC {\pRoScAC , 33}
\edef\pRoScAC {\pRoScAC , 34}
\edef\pRoScAC {\pRoScAC , 34}
\edef\pSchlAB {\pSchlAB , 34}
\edef\pSchlAG {34}
\edef\pKratBN {\pKratBN , 35}
\edef\pHaKrAA {35}
\edef\pKratBN {\pKratBN , 35}
\edef\pScotAB {35}
\edef\pMuirAB {\pMuirAB , 35}
\edef\pFulmAF {35}
\edef\pForrAC {35}
\edef\pFulmAF {\pFulmAF , 35}
\edef\pKratBN {\pKratBN , 36}
\edef\pGaRaAA {36}
\edef\pFoHaAL {36}
\edef\pKoorAG {36}
\edef\pScheAA {37}
\edef\pKratBN {\pKratBN , 37}
\edef\pKratBN {\pKratBN , 37}
\edef\pJohWAF {37}
\edef\pKratBN {\pKratBN , 37}
\edef\pJohWAF {\pJohWAF , 37}
\edef\pCiucAL {37}
\edef\pCiucAL {\pCiucAL , 37}
\edef\pMuirAD {39}
\edef\pHaKrAA {\pHaKrAA , 39}
\edef\pKratBN {\pKratBN , 39}
\edef\pMinaAA {40}
\edef\pKedlAA {40}
\edef\pStWiAA {40}
\edef\pKratBN {\pKratBN , 40}
\edef\pChuWBG {40}
\edef\pChuWBG {\pChuWBG , 40}
\edef\pChuWBG {\pChuWBG , 40}
\edef\pChuWBG {\pChuWBG , 40}
\edef\pJohWAE {40}
\edef\pArboAA {40}
\edef\pCraiAA {40}
\edef\pArboAA {\pArboAA , 40}
\edef\pLacrAA {40}
\edef\pJohWAE {\pJohWAE , 40}
\edef\pCraiAA {\pCraiAA , 40}
\edef\pChuWBG {\pChuWBG , 41}
\edef\pChuWBG {\pChuWBG , 41}
\edef\pKratBN {\pKratBN , 41}
\edef\pKratBN {\pKratBN , 41}
\edef\pKratBN {\pKratBN , 41}
\edef\pWallCF {41}
\edef\pVienAE {41}
\edef\pVienAE {\pVienAE , 42}
\edef\pVienAE {\pVienAE , 42}
\edef\pKratBN {\pKratBN , 42}
\edef\pIsStAB {43}
\edef\pIsStAC {43}
\edef\pIsStAB {\pIsStAB , 43}
\edef\pIsStAC {\pIsStAC , 43}
\edef\pKratBN {\pKratBN , 43}
\edef\pKoSwAA {43}
\edef\pIsStAB {\pIsStAB , 43}
\edef\pIsStAC {\pIsStAC , 43}
\edef\pKratBN {\pKratBN , 43}
\edef\pJoThAA {43}
\edef\pVienAE {\pVienAE , 43}
\edef\pTammAA {44}
\edef\pEgRiAA {44}
\edef\pGeXiAB {44}
\edef\pEgRiAA {\pEgRiAA , 44}
\edef\pMiRRAD {46}
\edef\pBresAO {\pBresAO , 46}
\edef\pKupeAH {\pKupeAH , 46}
\edef\pBresAO {\pBresAO , 46}
\edef\pEgRiAA {\pEgRiAA , 46}
\edef\pKratBN {\pKratBN , 46}
\edef\pKratBN {\pKratBN , 46}
\edef\pGhKrAA {46}
\edef\pKratBN {\pKratBN , 46}
\edef\pWimpAB {46}
\edef\pKratZZ {46}
\edef\pCvRIAA {46}
\edef\pAnWiAA {46}
\edef\pAignAA {46}
\edef\pCiglAM {46}
\edef\pCiglAO {46}
\edef\pCiglAV {46}
\edef\pEhreAB {46}
\edef\pLuThAB {47}
\edef\pLascAZ {47}
\edef\pLascAZ {\pLascAZ , 47}
\edef\pLascAZ {\pLascAZ , 47}
\edef\pHoLMAA {47}
\edef\pHoLMAB {47}
\edef\pKratBN {\pKratBN , 47}
\edef\pKratBN {\pKratBN , 47}
\edef\pAndrAO {47}
\edef\pAndrAN {47}
\edef\pKratBN {\pKratBN , 47}
\edef\pLascAZ {47}
\edef\pAndrAO {\pAndrAO , 48}
\edef\pMiRRAD {\pMiRRAD , 48}
\edef\pKratBN {\pKratBN , 48}
\edef\pCiEKAA {\pCiEKAA , 48}
\edef\pKratBN {\pKratBN , 48}
\edef\pAndrAN {\pAndrAN , 48}
\edef\pAndrAN {\pAndrAN , 49}
\edef\pAndrAO {\pAndrAO , 49}
\edef\pCiEKAA {\pCiEKAA , 49}
\edef\pEisTAF {\pEisTAF , 50}
\edef\pEisTAG {50}
\edef\pCiEKAA {\pCiEKAA , 51}
\edef\pCiEKAA {\pCiEKAA , 51}
\edef\pCiEKAA {\pCiEKAA , 51}
\edef\pCiKrAC {\pCiKrAC , 51}
\edef\pFuKrAC {\pFuKrAC , 52}
\edef\pKratBN {\pKratBN , 52}
\edef\pBacRAA {52}
\edef\pAmZeAB {\pAmZeAB , 53}
\edef\pKratBN {\pKratBN , 53}
\edef\pPeWiAA {53}
\edef\pZaPeAA {53}
\edef\pBacRAA {\pBacRAA , 53}
\edef\pKratBU {53}
\edef\pKratBU {\pKratBU , 53}
\edef\pKratBU {\pKratBU , 53}
\edef\pNeuwAE {53}
\edef\pPeWZAA {\pPeWZAA , 53}
\edef\pWegsAA {53}
\edef\pWiZeAC {53}
\edef\pZeilAM {\pZeilAM , 53}
\edef\pZeilAV {\pZeilAV , 53}
\edef\pNeuwAE {\pNeuwAE , 54}
\edef\pNeuwAE {\pNeuwAE , 54}
\edef\pHumpAC {54}
\edef\pKratBN {\pKratBN , 54}
\edef\pKratBN {\pKratBN , 54}
\edef\pHumpAC {\pHumpAC , 55}
\edef\pKratBN {\pKratBN , 55}
\edef\pVarcAC {55}
\edef\pReivAC {55}
\edef\pAdRoAC {56}
\edef\pAdRoAC {\pAdRoAC , 56}
\edef\pDGGGAA {56}
\edef\pSaWaAA {56}
\edef\pAdRoAC {\pAdRoAC , 57}
\edef\pAdRoAC {\pAdRoAC , 57}
\edef\pAdRoAC {\pAdRoAC , 57}
\edef\pAdRoAC {\pAdRoAC , 57}
\edef\pAdBRAA {57}
\edef\pKratBN {\pKratBN , 58}
\edef\pSoloAA {58}
\edef\pChoCAA {58}
\edef\pLamTAA {58}
\edef\pAdBRAA {\pAdBRAA , 58}
\edef\pAdBRAB {58}
\edef\pBagnAA {58}
\edef\pBernAA {58}
\edef\pBernAB {58}
\edef\pBiagAA {58}
\edef\pBiCaAA {58}
\edef\pFoHaAM {58}
\edef\pFoHaAN {58}
\edef\pFoHaAO {58}
\edef\pHaLRAA {58}
\edef\pReRoAA {58}
\edef\pBiagAA {\pBiagAA , 58}
\edef\pBiCaAA {\pBiCaAA , 58}
\edef\pReivAC {\pReivAC , 58}
\edef\pReRoAB {59}
\edef\pSchFAA {59}
\edef\pDahaAA {59}
\edef\pTuttAC {59}
\edef\pKratBN {\pKratBN , 59}
\edef\pStanAP {59}
\edef\pSimiAD {59}
\edef\pReivAD {59}
\edef\pSchFAA {\pSchFAA , 60}
\edef\pSchFAA {\pSchFAA , 60}
\edef\pSchFAA {\pSchFAA , 60}
\edef\pSchFAA {\pSchFAA , 60}
\edef\pSlPlAA {60}
\edef\pSloaAA {60}
\edef\pSchFAA {\pSchFAA , 61}
\edef\pLickAA {61}
\edef\pKoSmAA {61}
\edef\pDiFrAA {61}
\edef\pKratBN {\pKratBN , 61}
\edef\pCoSSAA {61}
\edef\pSchFAA {\pSchFAA , 61}
\edef\pAtReAA {61}
\edef\pReivAD {\pReivAD , 61}
\edef\pLindAB {61}
\edef\pKratBN {\pKratBN , 61}
\edef\pStanAP {\pStanAP , 61}
\edef\pSchFAA {\pSchFAA , 62}
\edef\pSmitAA {62}
\edef\pAlSTAA {62}
\edef\pStanAP {\pStanAP , 62}
\edef\pAlSTAA {\pAlSTAA , 62}
\edef\pBrMoAA {63}
\edef\pBrKMAA {63}
\edef\pBrMoAB {63}
\edef\pBrKMAA {\pBrKMAA , 63}
\edef\pBeRRAA {63}
\edef\pStanBI {\pStanBI , 64}
\edef\pBeRRAA {\pBeRRAA , 64}
\edef\pKrRiAA {64}
\edef\pKrRiAAA {64}
\edef\pLascAZ {\pLascAZ , 64}
\edef\pMacdAC {64}
\edef\pStanBI {\pStanBI , 64}
\edef\pKrRiAA {\pKrRiAA , 64}
\edef\pKratBN {\pKratBN , 65}
\edef\pKrRiAA {\pKrRiAA , 65}
\edef\pMacdAC {\pMacdAC , 66}
\edef\pLascAZ {\pLascAZ , 66}
\edef\pMacdAE {66}
\edef\pKoVeAA {67}
\edef\pLascAZ {\pLascAZ , 68}
\edef\pJameAA {69}
\edef\pKratBN {\pKratBN , 69}
\edef\pJameAA {\pJameAA , 69}
\edef\pJaKeAA {69}
\edef\pSagaAQ {69}
\edef\pMacdAC {\pMacdAC , 69}
\edef\pBeOSAA {70}
\edef\pOlssAB {70}
\edef\pBeOlAE {70}
\edef\pWW {70}
\edef\pWW {\pWW , 70}
\edef\pGaRaAA {70}
\edef\pKratBN {\pKratBN , 70}
\edef\pMilnAO {70}
\edef\pMilnAP {70}
\edef\pKratBN {\pKratBN , 70}
\edef\pKN {71}
\edef\pRainAA {71}
\edef\pRoseAA {71}
\edef\pRo {71}
\edef\pRoScAB {71}
\edef\pSp {71}
\edef\pWarnAG {71}
\edef\pH {71}
\edef\pRu {71}
\edef\pRoseAB {71}
\edef\pFr {71}
\edef\pRu {\pRu , 71}
\edef\pKN {\pKN , 71}
\edef\pRainAA {\pRainAA , 71}
\edef\pRoseAB {\pRoseAB , 71}
\edef\pF {71}
\edef\pAmdeAC {71}
\edef\pOkadAK {71}
\edef\pRoScAC {\pRoScAC , 71}
\edef\pKratBN {\pKratBN , 71}
\edef\pMacdAA {71}
\edef\pRoScAC {\pRoScAC , 71}
\edef\pHumpAC {\pHumpAC , 72}
\edef\pRoScAC {\pRoScAC , 72}
\edef\pRoScAC {\pRoScAC , 72}
\edef\pRoScAC {\pRoScAC , 73}
\edef\pMacdAA {\pMacdAA , 73}
\edef\pWarnAG {\pWarnAG , 74}
\edef\pWarnAG {\pWarnAG , 74}
\edef\pKratBN {\pKratBN , 74}
\edef\pKratBN {\pKratBN , 74}
\edef\pRoseAA {\pRoseAA , 74}
\edef\pRo {\pRo , 74}
\edef\pRoScAB {\pRoScAB , 74}
\edef\pSp {\pSp , 74}
\edef\pWarnAG {\pWarnAG , 74}
\edef\pKratBN {\pKratBN , 74}
\edef\pKratBN {\pKratBN , 75}
\edef\pKratBN {\pKratBN , 75}
\edef\pRainAA {\pRainAA , 75}
\edef\pRoScAC {\pRoScAC , 75}
\edef\pRoScAC {\pRoScAC , 75}
\edef\pSchlAB {\pSchlAB , 75}
\edef\pSchlAF {\pSchlAF , 75}
\edef\pRoScAC {\pRoScAC , 76}
\edef\pKratBN {\pKratBN , 76}
\edef\pKratBN {\pKratBN , 76}
\edef\pTV {76}
\edef\pRoScAC {\pRoScAC , 76}
\edef\pRoScAC {\pRoScAC , 76}
\edef\pRoScAC {\pRoScAC , 76}
\edef\pAnStAA {77}
\edef\pKratBN {\pKratBN , 77}
\edef\pWarnAG {\pWarnAG , 77}
\edef\pMilnAO {\pMilnAO , 77}
\edef\pMilnAP {\pMilnAP , 77}
\edef\pWimpAB {\pWimpAB , 82}

\numberwithin{equation}{section}

\newtheorem{Theorem}{Theorem}

\newtheorem{Lemma}[Theorem]{Lemma}

\newtheorem{Conjecture}[Theorem]{Conjecture}
\newtheorem{Problem}[Theorem]{Problem}
\newtheorem*{DefinitionA}{``Definition"}

\theoremstyle{definition}

\theoremstyle{remark}

\catcode`\@=11
\newwrite\Seiten
\immediate\openout\Seiten detcomp.sei\relax

\def\machSeite#1{\@ifundefined{@p#1}%
{\write\Seiten{\string\edef\csname p#1\endcsname{\thepage}}%
\@namedef{@p#1}{}}%
{\write\Seiten{\string\edef\csname p#1\endcsname{\csname p#1\endcsname , 
\thepage}}}%
}

\def\@secnumfont{\bfseries}
\def\section{\@startsection{section}{1}%
  \z@{.7\linespacing\@plus\linespacing}{.5\linespacing}%
  {\normalfont\bfseries}}
\def\subsection{\@startsection{subsection}{2}%
  \z@{.5\linespacing\@plus.7\linespacing}{-.5em}%
  {\normalfont\scshape}}
\catcode`\@=12

\catcode`\@=11
\def\sideset#1#2#3{%
  \@mathmeasure\z@\displaystyle{#3}%
  \global\setbox\@ne\vbox to\ht\z@{}\dp\@ne\dp\z@
  \setbox\tw@\box\@ne
  \@mathmeasure4\displaystyle{\copy\tw@#1}%
  \@mathmeasure6\displaystyle{#3{#2}}%
  \dimen@-\wd6 \advance\dimen@\wd4 \advance\dimen@\wd\z@
  \hbox to\dimen@{}\mathop{\kern-\dimen@\box4\box6}%
}
\catcode`\@=12

{\obeylines
\gdef\MATH{\begingroup\parindent0pt\parskip0pt plus 0pt\obeylines%
        \def^^M{\vskip4pt}%
        \obeyspaces\tt\small}%
}
\def\goodbreakpoint{\par\penalty-5000%
         \vrule height10pt depth2pt width0pt\leavevmode}
\def\endMATH{\endgroup}
\def\MATHphi{\leavevmode
        \hbox to 0pt{\hbox to 5.24995pt{\hss$\phi$\hss}\hss}}
\def\MATHGamma{\leavevmode
        \hbox to 0pt{\hbox to 5.24995pt{\hss$\Gamma$\hss}\hss}}
\def\MATHpi{\leavevmode
        \hbox to 0pt{\hbox to 5.24995pt{\hss$\pi$\hss}\hss}}
\def\MATHinfty{\leavevmode
        \hbox to 0pt{\hbox to 5.24995pt{\hss$\infty$\hss}\hss}}
\def\MATHhStrich{\leavevmode
        \hbox to 0pt{\hbox to 5.24995pt{\vrule height4.5pt depth-3.5pt width5.24995pt}\hss}}
\def\MATHluEck{\leavevmode
        \hbox to 0pt{\hbox to 5.24995pt{\hskip2.12497pt
         \vrule height4.5pt depth1pt width1pt
         \vrule height4.5pt depth-3.5pt width2.12498pt}\hss}}
\def\MATHruEck{\leavevmode
        \hbox to 0pt{\hbox to 5.24995pt{%
         \vrule height4.5pt depth-3.5pt width2.12497pt
         \vrule height4.5pt depth1pt width1pt
         \hskip2.12498pt}\hss}}
\def\MATHloEck{\leavevmode
        \hbox to 0pt{\hbox to 5.24995pt{\hskip2.12497pt
         \vrule height9pt depth-3.5pt width1pt
         \vrule height4.5pt depth-3.5pt width2.12498pt}\hss}}
\def\MATHroEck{\leavevmode
        \hbox to 0pt{\hbox to 5.24995pt{%
         \vrule height4.5pt depth-3.5pt width2.12497pt
         \vrule height9pt depth-3.5pt width1pt
         \hskip2.12498pt}\hss}}
\def\MATHvStrich{\leavevmode
        \hbox to 0pt{\hbox to 5.24995pt{\hskip2.12497pt
         \vtop to 0pt{\hsize1pt\vss%
                \vrule height17pt depth6pt width1pt\vskip8pt\vss\par}%
         \hskip2.12498pt}\hss}}
\def\MATHtStueck{\leavevmode
        \hbox to 0pt{\hbox to 5.24995pt{%
         \vrule height4.5pt depth-3.5pt width2.12497pt
         \vrule height4.5pt depth2pt width1pt
         \vrule height4.5pt depth-3.5pt width2.12498pt}\hss}}
\def\MATHbackslash{\leavevmode
        \hbox to 0pt{\hbox to 5.24995pt{\hss$\backslash$\hss}\hss}}
\def\MATHlbrace{\leavevmode
        \hbox to 0pt{\hbox to 5.24995pt{\hss$\{$\hss}\hss}}
\def\MATHrbrace{\leavevmode
        \hbox to 0pt{\hbox to 5.24995pt{\hss$\}$\hss}\hss}}
\def\MATHkleiner{\leavevmode
        \hbox to 0pt{\hbox to 5.24995pt{\hss$\langle$\hss}\hss}}
\def\MATHkleiner{<}
\def\MATHgroesser{\leavevmode
        \hbox to 0pt{\hbox to 5.24995pt{\hss$\rangle$\hss}\hss}}
\def\MATHgroesser{>}
\def\MATHhoch{\leavevmode
        \hbox to 0pt{\hbox to 5.24995pt{\hss$^\land$\hss}\hss}}
\def\MATHtief{\leavevmode
        \hbox to 0pt{\hbox to 5.24995pt{\hss\vrule height0pt depth.8pt width3pt\hss}\hss}}

\def\ringerl(#1 #2){\move(#1 #2)\fcir f:0 r:.1}
\def\hdSchritt{\bsegment
  \lpatt(.05 .13)
  \rlvec(0.866025403784439 -.5) 
  \savepos(0.866025403784439 -.5)(*ex *ey)
        \esegment
  \move(*ex *ey)
        }
\def\vdSchritt{\bsegment
  \lpatt(.05 .13)
  \rlvec(0 -1) 
  \savepos(0 -1)(*ex *ey)
        \esegment
  \move(*ex *ey)
        }
\def\ldreieck{\bsegment
  \rlvec(0.866025403784439 .5) \rlvec(0 -1)
  \rlvec(-0.866025403784439 .5)
  \savepos(0.866025403784439 -.5)(*ex *ey)
        \esegment
  \move(*ex *ey)
        }
\def\rdreieck{\bsegment
  \rlvec(0.866025403784439 -.5) \rlvec(-0.866025403784439 -.5)  \rlvec(0 1)
  \savepos(0 -1)(*ex *ey)
        \esegment
  \move(*ex *ey)
        }
\def\rhombus{\bsegment
  \rlvec(0.866025403784439 .5) \rlvec(0.866025403784439 -.5)
  \rlvec(-0.866025403784439 -.5)  \rlvec(0 1)
  \rmove(0 -1)  \rlvec(-0.866025403784439 .5)
  \savepos(0.866025403784439 -.5)(*ex *ey)
        \esegment
  \move(*ex *ey)
        }
\def\RhombusA{\bsegment
  \rlvec(0.866025403784439 .5) \rlvec(0.866025403784439 -.5)
  \rlvec(-0.866025403784439 -.5) \rlvec(-0.866025403784439 .5)
  \savepos(0.866025403784439 -.5)(*ex *ey)
        \esegment
  \move(*ex *ey)
        }
\def\RhombusB{\bsegment
  \rlvec(0.866025403784439 .5) \rlvec(0 -1)
  \rlvec(-0.866025403784439 -.5) \rlvec(0 1)
  \savepos(0 -1)(*ex *ey)
        \esegment
  \move(*ex *ey)
        }
\def\RhombusC{\bsegment
  \rlvec(0.866025403784439 -.5) \rlvec(0 -1)
  \rlvec(-0.866025403784439 .5) \rlvec(0 1)
  \savepos(0.866025403784439 -.5)(*ex *ey)
        \esegment
  \move(*ex *ey)
        }

\def\al{\alpha}
\def\be{\beta}
\def\ga{\gamma}
\def\de{\delta}

\def\la{\lambda}
\def\rh{\rho}
\def\si{\sigma}

\def\om{\omega}
\def\Ga{\Gamma}

\def\R{{\mathbb R}}

\def\Z{{\mathbb Z}}

\def\today{\ifcase\month\or
 January\or February\or March\or April\or May\or June\or
 July\or August\or September\or October\or November\or December\fi
 \space\number\day, \number\year}

\def\({\left(}
\def\){\right)}
\def\[{\left[}
\def\]{\right]}

\def\Comp{\operatorname{Comp}}
\def\Part{\operatorname{Part}}

\def\inv{\operatorname{inv}}
\def\maj{\operatorname{maj}}
\def\fmaj{\operatorname{fmaj}}
\def\nmaj{\operatorname{nmaj}}
\def\neg{\operatorname{neg}}
\def\sneg{\operatorname{sneg}}

\def\zbk{\operatorname{zbk}}
\def\nzbk{\operatorname{nzbk}}
\def\NC{\operatorname{NC}}
\def\match{\operatorname{NCmatch}}

\def\sgn{\operatorname{sgn}}

\def\3{\ss}

\def\fl#1{\left\lfloor#1\right\rfloor}
\def\cl#1{\left\lceil#1\right\rceil}
\def\coef#1{\left\langle#1\right\rangle}

\def\NC{\operatorname{NC}}

\def\stat{\operatorname{stat}}
\def\match{\operatorname{NCmatch}}
\def\Pf{\operatorname{Pf}}
\DeclareMathOperator{\h}{H}
\def\QED{\raise-15pt\hbox{{\Huge$\square$}}}
\def\qbinom#1#2{\left[\begin{smallmatrix} #1\\#2\end{smallmatrix}\right]_q}

\begin{document}

\newbox\Adr
\setbox\Adr\vbox{
\centerline{\sc C.~Krattenthaler$^\dagger$}
\vskip18pt
\centerline{Institut Camille Jordan, Universit\'e Claude Bernard
Lyon-I,}
\centerline{21, avenue Claude Bernard, F-69622 Villeurbanne Cedex,
France.}
\centerline{E-mail: {\tt\footnotesize kratt@euler.univ-lyon1.fr}}
\centerline{WWW: \footnotesize\tt http://igd.univ-lyon1.fr/\~{}kratt}
}

\title{Advanced Determinant Calculus: A Complement}

\author[C.~Krattenthaler]{\box\Adr}

\address{Institut Girard 
Desargues, Universit\'e Claude Bernard Lyon-I, 
21, avenue Claude Bernard, F-69622 Villeurbanne Cedex, France.}
\email{kratt@euler.univ-lyon1.fr}

\thanks{$^\dagger$ Research partially supported by EC's IHRP Programme,
grant HPRN-CT-2001-00272, ``Algebraic Combinatorics in Europe", and by
the ``Algebraic Combinatorics" Programme during Spring 2005 
of the Institut Mittag--Leffler of the Royal Swedish Academy of Sciences}


\begin{abstract}
This is a complement to my previous article {``Advanced
Determinant Calculus"} ({\it S\'eminaire Lotharingien Combin.}\ {\bf
42} (1999), Article~B42q, 67~pp.). In the present article, I 
share with the reader my experience of applying the
methods described in 
the previous article in order to solve a particular problem 
from number theory (G.~Almkvist, J.~Petersson and the author,
{\it Experiment.\ Math.}\ {\bf 12} (2003), 441--456).
Moreover, I add a list of determinant evaluations which I consider as 
interesting, which have been found since
the appearance of the previous article, or which I failed to mention
there, including several conjectures and open problems.
\end{abstract}

\maketitle

\section{Introduction} 
In the previous article 
\machSeite{KratBN}\cite{KratBN}, 
I described several methods to evaluate determinants, and 
I provided a long list of known determinant evaluations.
The present article is meant as a complement to 
\machSeite{KratBN}\cite{KratBN}.
Its purpose is three-fold: first, I want to shed light on
the problem of evaluating determinants from a slightly different
angle, by sharing with the reader my experience of applying the
methods from 
\machSeite{KratBN}\cite{KratBN} in order to solve a particular problem 
from number theory (see Sections~\ref{sec:det} and \ref{sec:eval}); 
second, I shall address the question
why it is apparently in the first case combinatorialists 
(such as myself) who are so interested in determinant evaluations
and get so easily excited about them (see Section~\ref{sec:comb}); 
and, finally third,
I add a list of determinant evaluations, which I consider 
as interesting, which have been found since
the appearance of 
\machSeite{KratBN}\cite{KratBN}, or which I failed to mention in
the list given in Section~3 of 
\machSeite{KratBN}\cite{KratBN} (see
Section~\ref{sec:detlist}), including several conjectures and open
problems.

\section{Enumerative combinatorics, nice formulae, and determinants}
\label{sec:comb}

Why are combinatorialists so fascinated by determinant
evaluations?

A simplistic answer to this question goes
as follows. Clearly, binomial coefficients $\binom nk$ or 
Stirling numbers (of the second kind) $S(n,k)$ are basic
objects in (enumerative) combinatorics; after all they count the
subsets of cardinality $k$ of a set with $n$ elements,
respectively the ways of partitioning such a set of $n$
elements into $k$ pairwise disjoint non-empty subsets. Thus, if one
sees an identity such as\footnote{For more information on this
determinant see Theorems~\ref{thm:MM} and \ref{thm:MM2} in this section
and 
\machSeite{KratBN}\cite[Sections~2.2, 2.3 and 2.5]{KratBN}.}
\begin{equation} \label{eq:M1}
\det_{1\le i,j\le n}\(\binom {a+b}{a-i+j} \)=
\prod _{i=1} ^{n}\frac {(a+b+i-1)!\,(i-1)!} {(a+i-1)!\,(b+i-1)!},
\end{equation}
or\footnote{This determinant evaluation follows easily from the matrix
factorisation
$$\(S(i+j,i) \)_{1\le i,j\le n}=\((-1)^kk^i/(k!\,(i-k)!)\)_{1\le
i,k\le n}\cdot\(k^j\)_{1\le k,j\le n},$$
application of 
\machSeite{KratBN}\cite[Theorem~26, (3.14)]{KratBN} to the first
determinant, and
application of the Vandermonde determinant evaluation to the
second.}
\begin{equation} \label{eq:S2}
\det_{1\le i,j\le n}\(S(i+j,i) \)=\prod _{i=1} ^{n}i^i
\end{equation}
(and there are many more of that kind; see 
\machSeite{KratBN}\cite{KratBN} and Section~\ref{sec:detlist}),
there is an obvious excitement that one cannot escape. 

Although this is indeed an explanation which applies in many cases,
there is also an answer on a more substantial level, which brings us
to the reason why {\it I\/} like (and need) determinant evaluations.

The {\it favourite question} for an enumerative combinatorialist (such as
myself) is
$$\text {\it How many $\langle\dots\rangle$ are there?}$$
Here, $\langle\dots\rangle$ can be permutations with certain properties, 
certain
partitions, certain paths, certain trees, etc. The {\it
favourite theorem} then is:
\begin{Theorem} \label{thm:fav}
The number of $\langle\dots\rangle$ of size $n$ is equal to
$$NICE(n).$$
\end{Theorem}

I have already explained the meaning of $\langle\dots\rangle$. What does
$NICE(n)$ stand for?
Typical examples for $NICE(n)$ are formulae
such as 
\begin{equation} \label{eq:Cat}
\frac {1} {n+1}\binom {2n}n
\end{equation}
(\!{\it Catalan numbers}; cf.\
\machSeite{StanBI}\cite[Ex.~6.19]{StanBI}) or
\begin{equation} \label{eq:alt}
\prod _{i=0} ^{n-1}\frac {(3i+1)!}{(n+i)!}
\end{equation} 
(the number of $n\times n$ 
alternating sign matrices and several other combinatorial objects;
cf.\ 
\machSeite{BresAO}\cite{BresAO}). Let us be more precise.
\begin{DefinitionA} \label{def:1}
The symbol $NICE(n)$ is a formula of the type
\begin{equation} \label{eq:NICE}
\xi^n\cdot\text{\em Rat}(n)\cdot
\prod _{i=1} ^{k}\frac {(a_in+b_i)!} {(c_in+d_i)!},
\end{equation}
where $\text{\em Rat}(n)$ is a rational function in $n$,
and where $a_i,c_i\in \Z$ for $i=1,2,\dots,k$, 
$\Z$ denoting the set of integers.
The parameters $b_i,c_i,\xi$ can be arbitrary real or complex
numbers. {\em(}If necessary, $(a_in+b_i)!$ has to be interpreted as
$\Ga(a_in+b_i+1)$, where $\Ga(x)$ is the Euler gamma function,
and similarly for $(c_in+d_i)!$.{\em)}
\end{DefinitionA}

Clearly, the formulae \eqref{eq:Cat} and \eqref{eq:alt} fit this
``Definition".\footnote{The writing $NICE(n)$ is borrowed from 
Doron Zeilberger 
\machSeite{ZeilAP}\cite[Recitation~III]{ZeilAP}.
The technical term for a formula of the type \eqref{eq:NICE} is
{\it``hypergeometric term"}, see 
\machSeite{PeWZAA}\cite[Sec.~3.2]{PeWZAA}, whereas, most often,
the colloquial terms {\it``closed form"} or {\it``nice formula"} 
are used for it, see 
\machSeite{ZeilAP}\cite[Recitation~II]{ZeilAP}.
More recently, some authors call sequences given by
formulae of that type sequences of ``round" numbers, see 
\machSeite{KupeAH}\cite[Sec.~6]{KupeAH}.} 

If one is working on a particular problem,
how can one recognise that one is looking at a sequence of numbers
given by $NICE(n)$? The key observation is that, if we factorise
$(an+b)!$ into its prime factors, where $a$ and $b$ are integers,
then, as $n$ runs through the positive integers, the numbers
$(an+b)!$ explode quickly, whereas the prime factors
occurring in the factorisation will grow only moderately, more
precisely, they will grow roughly linearly. Thus, if we encounter a
sequence the prime factorisation of which has this property, we can
be sure that there is a formula $NICE(n)$ for this sequence.
Even better, as I explain in Appendix~A of 
\machSeite{KratBN}\cite{KratBN}, the program
{\tt Rate}\footnote{\label{foot:Rate}{\tt Rate} is available from {\tt
http://igd.univ-lyon1.fr/\~{}kratt}. It is
based on a rather simple algorithm which involves rational interpolation.
In contrast to what I read, with great surprise, in 
\machSeite{CoBWAA}\cite{CoBWAA}, 
the explanations of how {\tt Rate} works in Appendix~A of
\machSeite{KratBN}\cite{KratBN} can be read and understood without any
knowledge about
determinants and, in particular, without any knowledge of the fifty
or so pages that precede Appendix~A in 
\machSeite{KratBN}\cite{KratBN}.} 
will (normally\footnote{\label{foot:normally}{\tt Rate} 
will {\it always} be able to guess
a formula of the type \eqref{eq:NICE} if there are enough initial
terms of the sequence available. However, there is a larger class of
sequences which have the property that the size of the primes in the 
prime factorisation of the
terms of the sequence grows only slowly with $n$. These are sequences
given by formulae containing ``Abelian" factors, such as
$n^n$. Unfortunately, {\tt Rate} does not know how to handle 
such factors. Recently, Rubey
\machSeite{RubeAD}\cite{RubeAD} proposed an algorithm for covering
Abelian factors as well. His implementation {\tt Guess} is written 
in {\sl Axiom} and is available 
at {\tt http://www.mat.univie.ac.at/\~{}rubey/martin.html}.})
be able to guess the formula.

To illustrate this, let us look at a particular example.
Let us suppose that the first few values of our sequence are the
following:
\begin{multline*}
1, 2, 5, 14, 42, 132, 429, 1430, 4862, 16796, 58786, 208012, 742900, 
  2674440, \\
9694845, 35357670, 129644790, 477638700, 1767263190,
6564120420.
\end{multline*}
The prime factorisation of the second-to-last number is
(we are using {\sl
Mathematica} here)

\MATH
\goodbreakpoint%
In[1]:= FactorInteger[477638700]
\goodbreakpoint%
Out[1]= %
\MATHlbrace %
\MATHlbrace 2, 2%
\MATHrbrace , %
\MATHlbrace 3, 1%
\MATHrbrace , %
\MATHlbrace 5, 2%
\MATHrbrace , %
\MATHlbrace 7, 1%
\MATHrbrace , %
\MATHlbrace 11, 1%
\MATHrbrace , %
\MATHlbrace 23, 1%
\MATHrbrace , %
\MATHlbrace 29, 1%
\MATHrbrace , %
\MATHlbrace 31, 1%
\MATHrbrace %
\MATHrbrace 
\goodbreakpoint%
\endMATH
whereas the prime factorisations of the next-to-last and the last
number in this sequence are

\MATH
\goodbreakpoint%
In[2]:= FactorInteger[1767263190]
\goodbreakpoint%
Out[2]= %
\MATHlbrace %
\MATHlbrace 2, 1%
\MATHrbrace , %
\MATHlbrace 3, 1%
\MATHrbrace , %
\MATHlbrace 5, 1%
\MATHrbrace , %
\MATHlbrace 7, 1%
\MATHrbrace , %
\MATHlbrace 11, 1%
\MATHrbrace , %
\MATHlbrace 23, 1%
\MATHrbrace , %
\MATHlbrace 29, 1%
\MATHrbrace , %
\MATHlbrace 31, 1%
\MATHrbrace , 
 
\MATHgroesser    %
\MATHlbrace 37, 1%
\MATHrbrace %
\MATHrbrace 
\goodbreakpoint%
In[3]:= FactorInteger[6564120420]
\goodbreakpoint%
Out[3]= %
\MATHlbrace %
\MATHlbrace 2, 2%
\MATHrbrace , %
\MATHlbrace 3, 1%
\MATHrbrace , %
\MATHlbrace 5, 1%
\MATHrbrace , %
\MATHlbrace 11, 1%
\MATHrbrace , %
\MATHlbrace 13, 1%
\MATHrbrace , %
\MATHlbrace 23, 1%
\MATHrbrace , %
\MATHlbrace 29, 1%
\MATHrbrace , %
\MATHlbrace 31, 1%
\MATHrbrace , 
 
\MATHgroesser    %
\MATHlbrace 37, 1%
\MATHrbrace %
\MATHrbrace 
\goodbreakpoint%
\endMATH

\noindent
(To decipher this for the reader unfamiliar with {\sl Mathematica}:
the prime factorisation of the last number is
$ 2^ 2 3^ 1 5^ 1 11^ 1 13^ 1 23^ 1 29^ 1 31^ 1  37^ 1$.)
One observes, first of all, that the occurring prime factors are
rather small in comparison to the numbers of which they are factors,
and, second, that the size of the prime factors grows only very
slowly (from 31 to 37). Thus, we {\it can be sure} that there is a
``nice" formula $NICE(n)$ for this sequence. Indeed, {\tt Rate} needs only
the first five members of the sequence to come up with a guess for
$NICE(n)$:

\MATH
\goodbreakpoint%
In[4]:= \MATHkleiner{\MATHkleiner}rate.m
\goodbreakpoint%
In[5]:= Rate[1,2,5,14,42]
\goodbreakpoint%
\leavevmode%
           i0       1
\leavevmode%
          4   Gamma[- + i0]
\leavevmode%
                    2
Out[5]= ----------------------
\leavevmode%
        Sqrt[Pi] Gamma[2 + i0]
\goodbreakpoint%
\endMATH
As the reader will have guessed, {\tt Rate} uses the parameter $i_0$ 
instead of $n$. In fact, the formula is a fancy way to write 
$\frac {1} {i_0+1}\binom {2i_0}{i_0}$, that is, we were looking at
the sequence of Catalan numbers \eqref{eq:Cat}.

To see the sharp contrast, here are the first few terms of another sequence:
$$1,2,9,272,589185.$$
(Also these are combinatorial numbers. They count the perfect
matchings of the $n$-dimensional hypercube; cf.\
\machSeite{PropAH}\cite[Problem~19]{PropAH}.) 
Let us factorise the last two numbers:

\MATH
\goodbreakpoint%
In[6]:= FactorInteger[272]
\goodbreakpoint%
Out[6]= %
\MATHlbrace %
\MATHlbrace 2, 4%
\MATHrbrace , %
\MATHlbrace 17, 1%
\MATHrbrace %
\MATHrbrace 
\goodbreakpoint%
In[7]:= FactorInteger[589185]
\goodbreakpoint%
Out[7]= %
\MATHlbrace %
\MATHlbrace 3, 2%
\MATHrbrace , %
\MATHlbrace 5, 1%
\MATHrbrace , %
\MATHlbrace 13093, 1%
\MATHrbrace %
\MATHrbrace 
\goodbreakpoint%
\endMATH
The presence of the big prime factor 13093 in the last factorisation is 
a sure sign that we cannot expect a formula $NICE(n)$ as described
in the ``Definition" for this sequence of numbers. (There may well be a
simple formula of a different kind. 
It is not very likely, though. In any case, such a formula has not
been found up to this date.)

\medskip
Now, that I have sufficiently explained all the ingredients in the
``prototype theorem'' Theorem~\ref{thm:fav}, I can explain why theorems of this
form are so attractive (at least to me): the objects (i.e., the
permutations, partitions, paths, trees, etc.) that it deals
with are usually very simple to explain, the statement is very simple
and can be understood by anybody, the result $NICE(n)$ has a very
elegant form, and yet, very often it is not easy at all to give a
proof (not to mention a true {\it explanation} why such an elegant
result occurs.)

Here are two examples. They concern
{\it rhombus tilings}, by which I mean tilings
of a region by rhombi with side lengths 1 and angles of
$60^\circ$ and $120^\circ$. The first one is a one century old theorem
due to MacMahon 
{\machSeite{MacMAA}\cite[Sec.~429, $q\rightarrow 1$; proof in
Sec.~494]{MacMAA}}.\footnote{To be correct, MacMahon did not know
anything about rhombus tilings, they did not exist in enumerative
combinatorics at the time. The objects that he considered were
{\it plane partitions}. However, there is a very simple bijection
between plane partitions contained in an $a\times b\times c$ box 
and rhombus tilings of a hexagon
with side lengths $a,b,c,a,b,c$, as explained for example in 
\machSeite{DT}\cite{DT}.}

\begin{Theorem} 
\label{thm:MM}
The number of rhombus tilings of a hexagon with 
side lengths $a,b,c,a,b,c$ whose angles are $120^\circ$
{\em(}see Figure~\ref{fig:1}.a for an example of such a hexagon, and
Figure~\ref{fig:1}.b for an example of a rhombus tiling{\em)} is equal to
\begin{equation} \label{eq:M2}
\prod _{i=1} ^{c}\frac {(a+b+i-1)!\,(i-1)!} {(a+i-1)!\,(b+i-1)!}.
\end{equation}
\quad \quad \qed
\end{Theorem}

\begin{figure}[h]
\centertexdraw{
  \drawdim truecm  \linewd.02
  \rhombus \rhombus \rhombus \rhombus \ldreieck
  \move (-0.866025403784439 -.5)
  \rhombus \rhombus \rhombus \rhombus \rhombus \ldreieck
  \move (-1.732050807568877 -1)
  \rhombus \rhombus \rhombus \rhombus \rhombus \rhombus \ldreieck
  \move (-1.732050807568877 -1)
  \rdreieck
  \rhombus \rhombus \rhombus \rhombus \rhombus \rhombus \ldreieck
  \move (-1.732050807568877 -2)
  \rdreieck
  \rhombus \rhombus \rhombus \rhombus \rhombus \rhombus \ldreieck
  \move (-1.732050807568877 -3)
  \rdreieck
  \rhombus \rhombus \rhombus \rhombus \rhombus \rhombus 
  \move (-1.732050807568877 -4)
  \rdreieck
  \rhombus \rhombus \rhombus \rhombus \rhombus 
  \move (-1.732050807568877 -5)
  \rdreieck
  \rhombus \rhombus \rhombus \rhombus 
\move(8 0)
\bsegment
  \drawdim truecm  \linewd.02
  \rhombus \rhombus \rhombus \rhombus \ldreieck
  \move (-0.866025403784439 -.5)
  \rhombus \rhombus \rhombus \rhombus \rhombus \ldreieck
  \move (-1.732050807568877 -1)
  \rhombus \rhombus \rhombus \rhombus \rhombus \rhombus \ldreieck
  \move (-1.732050807568877 -1)
  \rdreieck
  \rhombus \rhombus \rhombus \rhombus \rhombus \rhombus \ldreieck
  \move (-1.732050807568877 -2)
  \rdreieck
  \rhombus \rhombus \rhombus \rhombus \rhombus \rhombus \ldreieck
  \move (-1.732050807568877 -3)
  \rdreieck
  \rhombus \rhombus \rhombus \rhombus \rhombus \rhombus 
  \move (-1.732050807568877 -4)
  \rdreieck
  \rhombus \rhombus \rhombus \rhombus \rhombus 
  \move (-1.732050807568877 -5)
  \rdreieck
  \rhombus \rhombus \rhombus \rhombus 
  \linewd.12
  \move(0 0)
  \RhombusA \RhombusB \RhombusB 
  \RhombusA \RhombusA \RhombusB \RhombusA \RhombusB \RhombusB
  \move (-0.866025403784439 -.5)
  \RhombusA \RhombusB \RhombusB \RhombusB \RhombusB
  \RhombusA \RhombusA \RhombusB \RhombusA 
  \move (-1.732050807568877 -1)
  \RhombusB \RhombusB \RhombusA \RhombusB \RhombusB \RhombusA
  \RhombusB \RhombusA \RhombusA 
  \move (1.732050807568877 0)
  \RhombusC \RhombusC \RhombusC 
  \move (1.732050807568877 -1)
  \RhombusC \RhombusC \RhombusC 
  \move (3.464101615137755 -3)
  \RhombusC 
  \move (-0.866025403784439 -.5)
  \RhombusC
  \move (-0.866025403784439 -1.5)
  \RhombusC
  \move (0.866025403784439 -2.5)
  \RhombusC \RhombusC 
  \move (0.866025403784439 -3.5)
  \RhombusC \RhombusC \RhombusC 
  \move (2.598076211353316 -5.5)
  \RhombusC 
  \move (0.866025403784439 -5.5)
  \RhombusC 
  \move (-1.732050807568877 -3)
  \RhombusC 
  \move (-1.732050807568877 -4)
  \RhombusC 
  \move (-1.732050807568877 -5)
  \RhombusC \RhombusC 
\esegment
\htext (-1.5 -9){\small a. A hexagon with sides $a,b,c,a,b,c$,}
\htext (-1.5 -9.5){\small \hphantom{a. }where $a=3$, $b=4$, $c=5$}
\htext (6.8 -9){\small b. A rhombus tiling of a hexagon}
\htext (6.8 -9.5){\small \hphantom{b. }with sides $a,b,c,a,b,c$}
\rtext td:0 (4.3 -4.1){$\sideset {} c 
    {\left.\vbox{\vskip2.6cm}\right\}}$}
\rtext td:60 (2.6 -.55){$\sideset {} {} 
    {\left.\vbox{\vskip2cm}\right\}}$}
\rtext td:120 (-.34 -.2){$\sideset {}  {}  
    {\left.\vbox{\vskip1.7cm}\right\}}$}
\rtext td:0 (-2.4 -3.6){$\sideset {c}  {} 
    {\left\{\vbox{\vskip2.6cm}\right.}$}
\rtext td:240 (-0.1 -6.9){$\sideset {}  {}  
    {\left.\vbox{\vskip2cm}\right\}}$}
\rtext td:300 (2.9 -7.3){$\sideset {}  {}  
    {\left.\vbox{\vskip1.7cm}\right\}}$}
\htext (-.9 0.2){$a$}
\htext (2.8 -.1){$b$}
\htext (3.2 -7.9){$a$}
\htext (-0.4 -7.65){$b$}
}
\caption{}
\label{fig:1}
\end{figure}

The second one is more recent, and is due to Ciucu, Eisenk\"olbl, Zare
and the author {\machSeite{CiEKAA}\cite[Theorem~1]{CiEKAA}}.

\begin{figure}
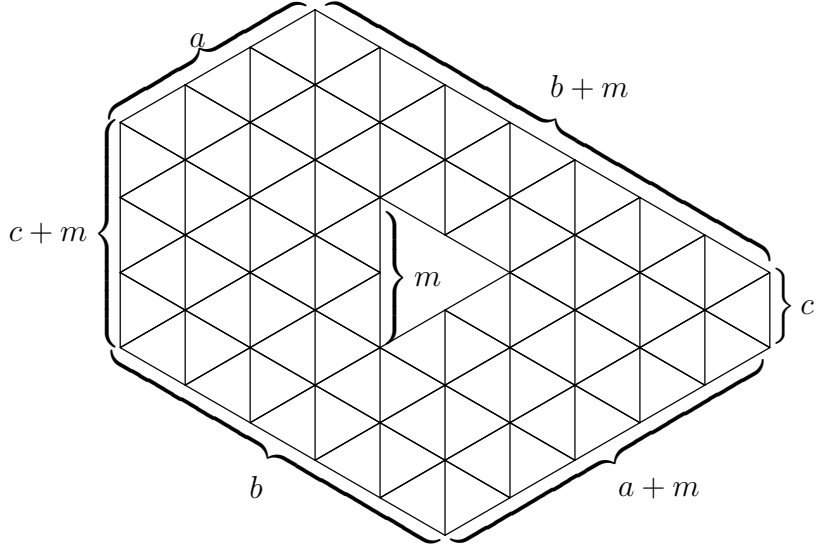
 
\centertexdraw{
\drawdim truecm \linewd.02
\rhombus \rhombus \rhombus \rhombus \rhombus \rhombus \rhombus
\ldreieck
\move(-.866025 -.5)
 \rhombus \rhombus \rhombus \rhombus \rhombus \rhombus \rhombus \rhombus 
\move(-1.7305 -1)
\rhombus \rhombus \rhombus \ldreieck \rmove(.866025 -.5)
\rhombus \rhombus \rhombus 
\move(-1.7305 -1)
\rdreieck \rhombus \rhombus \rhombus \ldreieck \rhombus \rhombus \rhombus 
\move(-1.7305 -2)
\rdreieck \rhombus \rhombus \rhombus \rhombus \rhombus \rhombus 
\move(-1.7305 -3)
\rdreieck \rhombus \rhombus \rhombus \rhombus \rhombus 
\htext(-.8 0){$a$}
\htext(4 -.7){$b+m$}
\htext(6.9 -3.95){$\left. \vbox{\vskip.6cm} \right\} c$}
\htext(-3.2 -4){$c+m \left\{ \vbox{\vskip1.6cm} \right.$}
\htext(0 -6){$b$}
\htext(4.9 -6){$a+m$}
\htext(1.7 -3.95){$\left. \vbox{\vskip1cm} \right\}m$}
\rtext td:60 (4 -1.3) {$\left. \vbox{\vskip3.6cm} \right\} $}
\rtext td:-60 (-.8 0){$\left\{ \vbox{\vskip1.6cm} \right. $}
\rtext td:-60 (4.6 -5.3) {$\left. \vbox{\vskip2.5cm} \right\} $}
\rtext td:60 (0.3 -5.6){$\left\{ \vbox{\vskip2.7cm}\right. $}
}
\caption{\protect\small A hexagon with triangular hole}
\label{hex}
\end{figure}

\begin{Theorem} \label{enum}
If $a,b,c$ have the same parity, then
the number of lozenge tilings of a hexagon with side lengths
$a,b+m,c,a+m,b,c+m$, with an equilateral triangle of
side length $m$ removed from its centre
{\em(}see Figure~\ref{hex} for an example{\em)}
is given by
\begin{multline} \label{eq:enum}
\frac {\h(a + m)\h(b + m)\h(c + m)\h(a + b + c + m)
} 
{\h(a + b + m)\h(a + c + m)\h(b + c + m)
}
\frac {\h(m + \left \lceil {\frac{a + b + c}{2}} \right \rceil)
\h(m + \left \lfloor {\frac{a + b + c}{2}} \right \rfloor)
} {\h({\frac{a + b}{2}} + m)     \h({\frac{a + c}{2}} + m)\h({\frac{b + c}{2}} + m)
}
\\
\times\frac {\h(\left \lceil {\frac{a}{2}} \right \rceil)
\h(\left \lceil {\frac{b}{2}} \right \rceil)
     \h(\left \lceil {\frac{c}{2}} \right \rceil)
     \h(\left \lfloor {\frac{a}{2}} \right \rfloor)\,
     \h(\left \lfloor {\frac{b}{2}} \right \rfloor)\,
     \h(\left \lfloor {\frac{c}{2}} \right \rfloor)\,
} 
{\h({\frac{m}{2}} + \left \lceil {\frac{a}{2}} \right \rceil)\,
     \h({\frac{m}{2}} + \left \lceil {\frac{b}{2}} \right \rceil)\,
     \h({\frac{m}{2}} + \left \lceil {\frac{c}{2}} \right \rceil)\,
\h({\frac{m}{2}} + \left \lfloor {\frac{a}{2}} \right \rfloor)\,
     \h({\frac{m}{2}} + \left \lfloor {\frac{b}{2}} \right \rfloor)\,
     \h({\frac{m}{2}} + \left \lfloor {\frac{c}{2}} \right \rfloor)\,
     }\\
\times
\frac {\h(\frac{m}{2})^2 \h({\frac{a + b + m}{2}})^2 
\h({\frac{a + c + m}{2}})^2 \h({\frac{b + c +
m}{2}})^2
} 
{\h({\frac{m}{2}} + \left \lceil {\frac{a + b + c}{2}} \right \rceil)
\h({\frac{m}{2}} + \left \lfloor {\frac{a + b + c}{2}} \right \rfloor)
\h({\frac{a + b}{2}})\h({\frac{a + c}{2}})\h({\frac{b + c}{2}})
               },
\end{multline}
where
\begin{equation} \label{eq:hyperfac} 
\h(n):=\begin{cases}
\prod _{k=0} ^{n-1}{\Gamma(k+1)}\quad &\text {for $n$ an integer,}\\
\prod _{k=0} ^{n-\frac {1} {2}}{\Gamma(k+\frac {1} {2})} \quad &\text 
{for $n$ a half-integer}.
\end{cases}
\end{equation}
\quad \quad \qed
\end{Theorem}

(There is a similar theorem if the parities of $a,b,c$ should not be the
same, see 
\machSeite{CiEKAA}\cite[Theorem~2]{CiEKAA}. Together, the two theorems
generalise MacMahon's Theorem~\ref{thm:MM}.\footnote{Bijective proofs
of Theorem~\ref{thm:MM} which ``explain" the ``nice" formula are
known 
\machSeite{KratAY}%
\machSeite{KratBK}%
\cite{KratAY,KratBK}. I do not ask for a bijective proof of
Theorem~\ref{enum} because I consider the task of finding one as daunting.})

The reader should notice that the right-hand side of \eqref{eq:M2} is
indeed of the form $NICE(a)$, while the right-hand side of
\eqref{eq:enum} is of the form $NICE(m/2)$. 

\medskip
Where is the connexion to determinants? As it turns out, these two
theorems are in fact {\it determinant evaluation theorems}. More
precisely, Theorem~\ref{thm:MM} is equivalent to the following
theorem.

\begin{Theorem} \label{thm:MM2}
\begin{equation} \label{eq:M3}
\det_{1\le i,j\le c}\(\binom {a+b}{a-i+j} \)=
\prod _{i=1} ^{c}\frac {(a+b+i-1)!\,(i-1)!} {(a+i-1)!\,(b+i-1)!}.
\end{equation}
\quad \quad \qed
\end{Theorem}

(The reader should notice that this is exactly \eqref{eq:M1} with $n$
replaced by $c$.) On the
other hand, Theorem~\ref{enum} is equivalent to the theorem
below.\footnote{To be correct, this is a little bit oversimplified.
The truth is that equivalence holds only if $m$ is even. An
additional argument is necessary for proving the result for the case
that $m$ is odd. We refer the reader who is interested in these
details to 
\machSeite{CiEKAA}\cite[Sec.~2]{CiEKAA}.}

\begin{Theorem} \label{enum2}
If $m$ is even, the determinant
\begin{equation} \label{mat1}
\det_{1\le i,j\le a+m} \begin{pmatrix} \dbinom{b+c+m}{b-i+j}&
\text {\scriptsize $1\le i\le a$}\\
\dbinom{\frac {b+c} {2}}{\frac {b+a} {2}-i+j}&
\text {\scriptsize $a+1\le i \le a+m$}
\end{pmatrix}
\end{equation}
is equal to \eqref{eq:enum}.\quad \quad \qed
\end{Theorem}

\begin{figure} 
\centertexdraw{
\drawdim cm \setunitscale.7
\linewd.01
\rhombus \rhombus \rhombus \rhombus \rhombus \rhombus \rhombus
\ldreieck
\move(-.866025 -.5)
 \rhombus \rhombus \rhombus \rhombus \rhombus \rhombus \rhombus \rhombus 
\move(-1.73205 -1)
\rhombus \rhombus \rhombus \ldreieck \rmove(.866025 -.5)
\rhombus \rhombus \rhombus 
\move(-1.73205 -1)
\rdreieck \rhombus \rhombus \rhombus \ldreieck \rhombus \rhombus \rhombus 
\move(-1.73205 -2)
\rdreieck \rhombus \rhombus \rhombus \rhombus \rhombus \rhombus 
\move(-1.73205 -3)
\rdreieck \rhombus \rhombus \rhombus \rhombus \rhombus

\linewd.08
\move(0 0)
\RhombusA  \RhombusA \RhombusA \RhombusA \RhombusA \RhombusA \RhombusA
\RhombusB 
\move(-.866025 -.5)
\RhombusA \RhombusA \RhombusA \RhombusA \RhombusA \RhombusA 
\RhombusB \RhombusA 
\move(-1.73205 -1)
\RhombusA \RhombusA \RhombusB \RhombusB \RhombusA \RhombusB \RhombusA 
\RhombusA 
\move(2.598 -3.5)
\RhombusA \RhombusB \RhombusA 
\move(1.73205 -4)
\RhombusB \RhombusA \RhombusA
\move(-1.73205 -1)
 \RhombusC \RhombusC
\move(-1.73205 -2) \RhombusC
\move(-1.73205 -3)\RhombusC \RhombusC \RhombusC
\move(.866025 -1.5) \RhombusC
\move(.866025 -2.5) \RhombusC
\htext(-4 -8){\small 
a. A lozenge tiling of the cored hexagon in Figure~\ref{hex}}

\htext(-1 0.2){$a$}
\htext(4.4 -1){$b+m$}
\htext(6.9 -3.9){$\left. \vbox{\vskip.35cm} \right\} c$}
\htext(-3.8 -4){$c+m \left\{ \vbox{\vskip1.2cm} \right.$}
\htext(0 -6.2){$b$}
\htext(5 -6){$a+m$}
\htext(1.7 -4){$\left. \vbox{\vskip0.77cm} \right\}m$}
\rtext td:60 (4 -1.3) {$\left. \vbox{\vskip2.52cm} \right\} $}
\rtext td:-60 (-.8 0.2){$\left\{ \vbox{\vskip1.02cm} \right. $}
\rtext td:-60 (4.6 -5.3) {$\left. \vbox{\vskip1.75cm} \right\} $}
\rtext td:60 (0.3 -5.8){$\left\{ \vbox{\vskip1.89cm}\right. $}

\move(11 0)
\bsegment
\linewd.05
\move(0 0)
\RhombusA  \RhombusA \RhombusA \RhombusA \RhombusA \RhombusA \RhombusA
\RhombusB 
\move(-.866025 -.5)
\RhombusA \RhombusA \RhombusA \RhombusA \RhombusA \RhombusA 
\RhombusB \RhombusA 
\move(-1.73205 -1)
\RhombusA \RhombusA \RhombusB \RhombusB \RhombusA \RhombusB \RhombusA 
\RhombusA 
\move(2.598 -3.5)
\RhombusA \RhombusB \RhombusA 
\move(1.73205 -4)
\RhombusB \RhombusA \RhombusA
\move(-1.73205 -1)
 \RhombusC \RhombusC
\move(-1.73205 -2) \RhombusC
\move(-1.73205 -3)\RhombusC \RhombusC \RhombusC
\move(.866025 -1.5) \RhombusC
\move(.866025 -2.5) \RhombusC

\ringerl(.433012 .25) 
\hdSchritt \hdSchritt  \hdSchritt \hdSchritt \hdSchritt \hdSchritt \hdSchritt  \vdSchritt  
\ringerl(-.433012 -.25)
\hdSchritt \hdSchritt \hdSchritt \hdSchritt \hdSchritt \hdSchritt  \vdSchritt \hdSchritt  
\ringerl(-1.299037 -.75) \hdSchritt \hdSchritt 
\vdSchritt \vdSchritt  \hdSchritt  \vdSchritt  \hdSchritt \hdSchritt 
\ringerl(3.031 -3.25) \hdSchritt \vdSchritt \hdSchritt
\ringerl(2.165 -3.75) \vdSchritt \hdSchritt \hdSchritt
\ringerl(6.4952 -4.25)
\ringerl(5.6292 -4.75)
\ringerl(4.7632 -5.25) 
\ringerl(3.8971 -5.75)
\ringerl(3.031  -6.25)
\htext(-2 -8){\small b. The corresponding path family}
\esegment
\htext(4 -18){\small c. The path family made orthogonal}}
\vskip-7cm
$$
\Einheit=.7cm
\Gitter(11,7)(0,0)
\Koordinatenachsen(11,7)(0,0)
\Pfad(0,3),11\endPfad
\Pfad(2,1),22\endPfad
\Pfad(2,1),1\endPfad
\Pfad(3,0),2\endPfad
\Pfad(3,0),11\endPfad
\Pfad(1,4),111111\endPfad
\Pfad(7,3),2\endPfad
\Pfad(7,3),1\endPfad
\Pfad(2,5),1111111\endPfad
\Pfad(9,4),2\endPfad
\Pfad(4,1),2\endPfad
\Pfad(4,1),11\endPfad
\Pfad(5,3),1\endPfad
\Pfad(6,2),2\endPfad
\Pfad(6,2),1\endPfad
\Kreis(0,3.02) \Kreis(5,0)
\Label\lo{A_1}(0,3) \Label\ru{E_1}(5,0)
\Kreis(1,4.02) \Kreis(8,3)
\Label\lo{A_2}(1,4) \Label\ru{E_4}(8,3)
\Kreis(2,5.02) \Kreis(9,4)
\Label\lo{A_3}(2,5) \Label\ru{E_5}(9,4)
\Kreis(4,2.02) \Kreis(6,1)
\Label\lo{A_4}(4,2) \Label\ru{E_2}(6,1)
\Kreis(5,3.02) \Kreis(7,2)
\Label\lo{A_5}(5,3) \Label\ru{E_3}(7,2)
\hskip6cm
$$
\vskip1cm
\caption{}
\label{tiling}
\end{figure}

The link between rhombus tilings (and equivalent objects such as
plane partitions, semistandard tableaux, etc.) and determinants 
which explains the above two equivalence statements is {\it
non-intersecting lattice paths}.\footnote{There exists in fact a
second link between rhombus tilings and determinants which is not
less interesting or less important. It is a well-known fact that
rhombus tilings are in bijection with perfect matchings of certain
hexagonal graphs. (See for example
\machSeite{KupeAG}\cite[Figures~13 and 14]{KupeAG}.) In view of this fact,
this second link is given by Kasteleyn's
theorem \machSeite{KastAA}\cite{KastAA} saying that the number of
perfect matchings of a planar graph is given by the Pfaffian of
a slight perturbation of the adjacency matrix of the graph.
See \machSeite{KupeAG}\cite{KupeAG} for an exposition of 
Kasteleyn's result, including historical notes, and for 
adaptations taking symmetries of the graph into
account.} The latter are families of paths in a
lattice with the property that no two paths in the family have a
point in common. Indeed, rhombus tilings are (usually) in bijection
with families of non-intersecting paths in the integer lattice
$\Z^2$ which consist of unit horizontal and vertical steps.
(Figure~\ref{tiling} illustrates 
the bijection for the rhombus tilings which appear in
Theorem~\ref{enum} in an example. 
In that bijection, all horizontal steps of the paths are in the positive
direction, and all vertical steps are in the negative direction. 
See the explanations that accompany 
\machSeite{CiEKAA}\cite[Figure~8]{CiEKAA} 
for a detailed description.
Since, as I explained, Theorem~\ref{enum}
essentially is a generalisation of Theorem~\ref{thm:MM}, this gives
also an idea for the bijection for the rhombus tilings which appear in
the latter theorem. For other instances of bijections between rhombus
tilings and non-intersecting lattice paths see
\machSeite{CiKrAA}%
\machSeite{CiKrAC}%
\machSeite{CiKrAD}%
\machSeite{EisTAA}%
\machSeite{EisTAB}%
\machSeite{EisTAF}%
\machSeite{FiscAA}%
\machSeite{KratBY}%
\machSeite{OkKrAA}%
\cite{CiKrAA,CiKrAC,CiKrAD,EisTAA,EisTAB,EisTAF,FiscAA,KratBY,OkKrAA}.). 
In the case that the starting points and the end points of the lattice
paths are fixed, the following
many-author-theorem applies.\footnote{%
This result was discovered and rediscovered several times. 
In a probabilistic form, it occurs for the
first time in work by Karlin and McGregor 
\machSeite{KaMGAB}%
\machSeite{KaMGAC}%
\cite{KaMGAB,KaMGAC}.
In matroid theory, it is discovered in its discrete form by Lindstr\"om 
\machSeite{LindAA}\cite[Lemma~1]{LindAA}. Then, in the 1980s the theorem is
rediscovered at about the same time in three different
communities, not knowing from each other at the time:
in statistical physics by Fisher 
\machSeite{FishAA}\cite[Sec.~5.3]{FishAA} 
in order to apply it to
the analysis of vicious walkers as a model of wetting and melting, 
in combinatorial chemistry by John and Sachs 
\machSeite{JoSaAB}\cite{JoSaAB} and
Gro\-nau, Just, Schade, Scheffler and Wojciechowski 
\machSeite{GrJSAA}\cite{GrJSAA}
in order to compute Pauling's bond order
in benzenoid hydrocarbon molecules, and in enumerative combinatorics
by Gessel and Viennot 
\machSeite{GeViAA}%
\machSeite{GeViAB}%
\cite{GeViAA,GeViAB} in order to count
tableaux and plane partitions. Since only Lindstr\"om, and then
Gessel and Viennot state the result in its most general form (not
reproduced here), I call this theorem most often
the ``Lindstr\"om--Gessel--Viennot theorem." It must be also
mentioned that the so-called ``Slater determinant"
in quantum mechanics (cf.\ 
\machSeite{SlatZY}\cite{SlatZY} and 
\machSeite{SlatZZ}\cite[Ch.~11]{SlatZZ}) 
may qualify as an ``ancestor" of
the Lindstr\"om--Gessel--Viennot determinant.}

\begin{Theorem}[\sc Karlin--McGregor, Lindstr\"om, Gessel--Viennot,
Fisher,\break John--Sachs, Gronau--Just--Schade--Scheffler--Wojciechowski]
\label{thm:nonint}
Let $A_1,A_2,\break \dots,A_n$ and $E_1,E_2,\dots,E_n$ be lattice points
such that for $i<j$ and $k<l$ any lattice path between $A_i$ and $E_l$ has a
common point with any lattice path between $A_j$ and $E_k$.
Then the number of all families $(P_1,P_2,\dots,P_n)$ of
non-intersecting lattice paths, $P_i$ running from $A_i$ to $E_i$,
$i=1,2,\dots,n$, is given by
$$\det_{1\le i,j\le n}\big(P(A_j\to E_i)\big),$$
where $P(A\to E)$ denotes the number of all lattice paths from $A$ to
$E$.\quad \quad \qed
\end{Theorem}

It goes beyond the scope of this article to include the proof of this
theorem here. However, I cannot help telling that it is an extremely
beautiful and simple proof that {\it every} mathematician 
should have seen once, even if (s)he
does not have any use for it in her/his own research. I refer the
reader to 
\machSeite{GeViAA}%
\machSeite{GeViAB}%
\machSeite{StemAE}%
\cite{GeViAA,GeViAB,StemAE}.

\medskip
Now the origin of the determinants becomes evident. In particular,
since, for rhombus tilings, we have to deal with lattice paths 
in the integer lattice consisting of unit horizontal and vertical 
steps, and since the number of such lattice paths
which connect two lattice points is given by a binomial coefficient,
we see that the enumeration of rhombus tilings must be a rich
source for binomial determinants. This is indeed the case, and there
are several instances in which such determinants can be evaluated in the
form $NICE(.)$ (see
\machSeite{CiucAH}%
\machSeite{CiEKAA}%
\machSeite{CiKrAA}%
\machSeite{CiKrA}%
\machSeite{CiKrAD}%
\machSeite{EisTAA}%
\machSeite{EisTAB}%
\machSeite{EisTAF}%
\machSeite{FiscAA}%
\machSeite{FuKrAC}%
\machSeite{KratBN}%
\cite{CiucAH,CiEKAA,CiKrAA,CiKrAC,CiKrAD,EisTAA,EisTAB,EisTAF,FiscAA,FuKrAC,KratBN}
and Section~\ref{sec:detlist}).
Often the evaluation part is highly non-trivial.

The evaluation of the determinant \eqref{eq:M3} is not very difficult
(see 
\machSeite{KratBN}\cite[Sections~2.2, 2.3, 2.5]{KratBN} for 3 different ways
to evaluate it). 
On the other hand, the
evaluation of the determinant \eqref{mat1} requires some effort (see
\machSeite{CiEKAA}\cite[Sec.~7]{CiEKAA}).

\medskip
To conclude this section, I state another determinant evaluation,
to which I shall come back later. Its origin lies as well in the
enumeration of rhombus tilings and plane partitions (see 
\machSeite{KratBD}\cite[Theorem~10]{KratBD} and 
\machSeite{CiKrAB}\cite[Theorem~2.1]{CiKrAB}).

\begin{Theorem} \label{thm:xy}
For any complex numbers $x$ and $y$ there holds
\begin{multline} \label{eq:Krat}
\det_{0\le i,j\le n-1}\(\frac {(x+y+i+j-1)!}
{(x+2i-j)!\,(y+2j-i)!}\)\\
=\prod _{i=0} ^{n-1}\frac {i!\,(x+y+i-1)!\,(2x+y+2i)_i\,(x+2y+2i)_i}
{(x+2i)!\,(y+2i)!},
\end{multline}
where the {\em shifted factorials} or 
{\em Pochhammer symbols} $(a)_k$ are defined by
$(a)_k:=a(a+1)\cdots(a+k-1)$, $k\ge1$, and $(a)_0:=1$.
{\em(}In this  formula,
a factorial $m!$ has to be interpreted as $\Ga(m+1)$ if $m$ is not a
non-negative integer.{\em)}
\end{Theorem}

\section{A determinant from number theory}\label{sec:det}

However, determinants do not only arise in combinatorics, they also
arise in other fields. In this section, I want to present 
a determinant which arose in number theory, explain in some detail
its origin, and then outline the steps which led to its evaluation,
thereby giving the reader an opportunity to look ``behind the scenes"
while one tries to make the determinant evaluation methods described
in 
\machSeite{KratBN}\cite{KratBN} work.

The story begins with the following two series expansions for $\pi$.
The first one is due to Bill Gosper 
\machSeite{Gosper}\cite{Gosper},
\begin{equation} \label{eq:Gosper}
\pi =\sum_{n=0}^\infty \frac{50n-6}{\binom{3n}n2^n} ,
\end{equation}
and was used by Fabrice Bellard 
\machSeite{Bellard}\cite[file {\tt pi1.c}]{Bellard} 
to find an algorithm for computing the $n$-th decimal of $\pi $ without
computing the earlier ones, thus improving an earlier algorithm due to
Simon Plouffe 
\machSeite{Plouffe}\cite{Plouffe}. 
The second one,
\begin{equation} \label{eq:Bellard} 
\pi=\frac {1} {740025}\(\sum _{n=1} ^{\infty}\frac {3P(n)} {\binom
{7n}{2n}2^{n-1}}-20379280\),
\end{equation}
where
\begin{multline*}
P(n)=-885673181n^5+3125347237n^4-2942969225n^3\\
+1031962795n^2-196882274n+10996648,
\end{multline*}
is due to Fabrice Bellard 
\machSeite{Bellard}\cite{Bellard}, and was used by him in his
world record setting
computation of the 1000 billionth {\it binary} digit of $\pi$, being based
on the algorithm in 
\machSeite{BaBoPl}\cite{BaBoPl}.

Going beyond that, my co-authors from 
\machSeite{AlKPAA}\cite{AlKPAA},
Gert Almkvist and Joakim Petersson, asked themselves the following
question:

{\em Are there more expansions of the type
$$\pi=\sum_{n=0}^\infty \frac {S(n)}{\binom{mn}{pn}a^n},$$
where $S(n)$ is some polynomial in $n$ {\em(}depending on $m,p,a${\em)}?}

How can one go about to get some intuition about this question? 
One chooses some specific
$m,p,a$, goes to the computer, computes
$$p(k)=\sum _{n=0} ^{\infty}\frac {n^k} {\binom {mn}{pn}a^n}$$
to many, many digits for $k=0,1,2,\dots$, puts
$$\pi,p(0),p(1),p(2),\dots$$
into the LLL-algorithm (which comes, for example, with the {\sl Maple}
computer algebra package), and one waits whether the algorithm
comes up with an integral linear combination of 
$\pi,p(0),p(1),p(2),\dots$.\footnote{For readers unfamiliar with the
LLL-algorithm: in this particular application,
it takes as an input rational numbers $r_1,r_2,\dots,r_m$ (which, in
our case, will be the numbers $1$ and the rational approximations of $\pi$,
$p(0)$, $p(1)$, \dots \ which we computed), and, if
successful, outputs {\it small\/} integers $c_1,c_2,\dots,c_m$ such
that $c_1r_1+c_2r_2+\dots+c_mr_m$ is {\it very small}. Thus, if $r_i$
was a good approximation for the real number $x_i$, $i=1,2,\dots,m$,
one can expect that actually $c_1x_1+c_2x_2+\dots+c_mx_m=0$.
See
\machSeite{LeLLAA}\cite[Sec.~1, in particular 
the last paragraph]{LeLLAA} and 
\machSeite{CohHAA}\cite[Ch.~2]{CohHAA} 
for the description of and more information on this important algorithm. 
In particular, also here, the output of the
algorithm (if there is) is just a (very guided) {\it guess}. Thus, a
proof is still needed, although the probability that the guess is
wrong is infinitesimal. As a matter of fact, it is very likely that
Bellard had no proof of his \hbox{formula \eqref{eq:Bellard} \dots}}
Indeed, Table~\ref{tab:2} shows the parameter values, where 
the LLL-algorithm gave a result.

\begin{table}[h]
\begin{tabular}{c|c|l|cl}
$m$ & $p$ & $\hphantom{-}a$ & $\deg(S)$ &  \\ 
\cline{1-4}
\hphantom{1}3 & \hphantom{1}1 & $\hphantom{-}2$ & \hphantom{1}1 & (Gosper) 
\\ 
\hphantom{1}7 & \hphantom{1}2 & $\hphantom{-}2$ & \hphantom{1}5 & (Bellard) 
\\ 
\hphantom{1}8 & \hphantom{1}4 & $-4$ & \hphantom{1}4 &  \\ 
10 & \hphantom{1}4 & $\hphantom{-}4$ & \hphantom{1}8 &  \\ 
12 & \hphantom{1}4 & $-4$ & \hphantom{1}8 &  \\ 
16 & \hphantom{1}8 & $\hphantom{-}16$ & \hphantom{1}8 &  \\ 
24 & 12 & $-64$ & 12 &  \\ 
32 & 16 & $\hphantom{-}256$ & 16 &  \\ 
40 & 20 & $-4^5$ & 20 &  \\ 
48 & 24 & $\hphantom{-}4^6$ & 24 &  \\ 
56 & 28 & $-4^7$ & 28 &  \\ 
64 & 32 & $\hphantom{-}4^8$ & 32 &  \\ 
72 & 36 & $-4^9$ & 36 &  \\ 
80 & 40 & $\hphantom{-}4^{10}$ & 40 & 
\end{tabular}
\vskip10pt
\caption{}
\label{tab:2}
\end{table}

For example, it found
$$
\pi =\frac 1r\sum_{n=0}^\infty \frac{S(n)}{\binom{16n}{8n}16^n}, 
$$
where
$$
r=3^65^37^211^213^2 
$$
and
\begin{multline*} 
S(n)=-869897157255-3524219363487888n+112466777263118189n^2 
\\
-1242789726208374386n^3+6693196178751930680n^4-19768094496651298112n^5 
\\
+32808347163463348736n^6-28892659596072587264n^7+10530503748472012800n^8,
\end{multline*}
and
$$
\pi =\frac 1r\sum_{n=0}^\infty \frac{S(n)}{\binom{32n}{16n}256^n} ,
$$
where
$$
r=2^33^{10}5^67^311^1 13^217^219^223^229^231^2 
$$
and
{\allowdisplaybreaks
\begin{align*}
S(n)=&-2062111884756347479085709280875 
\\
&+1505491740302839023753569717261882091900n 
\\
&-112401149404087658213839386716211975291975n^2 
\\
&+3257881651942682891818557726225840674110002n^3 
\\
&-51677309510890630500607898599463036267961280n^4 
\\
&+517337977987354819322786909541179043148522720n^5 
\\
&-3526396494329560718758086392841258152390245120n^6 
\\
&+171145766235995166227501216110074805943799363584n^7 
\\
&-60739416613228219940886539658145904402068029440n^8 
\\
&+159935882563435860391195903248596461569183580160n^9 
\\
&-313951952615028230229958218839819183812205608960n^{10} 
\\
&+457341091673257198565533286493831205566468325376n^{11} 
\\
&-486846784774707448105420279985074159657397780480n^{12} 
\\
&+367314505118245777241612044490633887668208926720n^{13} 
\\
&-185647326591648164598342857319777582801297080320n^{14} 
\\
&+56224688035707015687999128994324690418467340288n^{15} 
\\
&-7687255778816557786073977795149360408612044800n^{16} .
\end{align*}}%
Of course, there could be many more.

If one looks more closely at Table~\ref{tab:2}, then, if one
disregards the first, second and fourth line, one cannot escape to
notice a pattern: {\it apparently, for each $k=1,2,\dots$, there is a
formula
$$
\pi =\sum_{n=0}^\infty \frac{S_k(n)}{\binom{8kn}{4kn}(-4)^{kn}} ,
$$
where $S_k(n)$ is some polynomial in $n$ of degree $4k$.}

In order to make progress on this observation, we have to first
see how one can prove such an identity, once it is found. In fact,
this is not difficult at all. To illustrate the idea, let us go
through a proof of Gosper's identity \eqref{eq:Gosper}.

The beta integral evaluation (cf.\ 
\machSeite{AAR}\cite[Theorem~1.1.4]{AAR}) gives
$$
\frac 1{\binom{3n}n}=(3n+1)\int_0^1x^{2n}(1-x)^ndx .
$$
Hence the right hand side of the formula will be
$$
\int_0^1\sum_{n=0}^\infty (50n-6)(3n+1)\left(\frac {x^2(1-x)}{2}\right)^ndx
.
$$
We have
\begin{equation} \label{eq:rational} 
\sum_{n=0}^\infty (50n-6)(3n+1)y^n=\frac{2(56y^2+97y-3)}{(1-y)^3} .
\end{equation}
Thus, if substituted, we obtain
\begin{align}\notag
RHS&=8\int_0^1\frac{28x^6-56x^5+28x^4-97x^3+97x^2-6}{(x^3-x^2+2)^3}dx\\
&=
\left[ \frac{4x(x-1)(x^3-28x^2+9x+8)}{(x^3-x^2+2)^2}+4\arctan (x-1)\right]
_0^1=\pi .
\label{eq:RHS}
\end{align}
(Clearly, both \eqref{eq:rational} and \eqref{eq:RHS} are routine
calculations, and therefore we did not do it by hand, but let them be
worked out by {\sl Maple}.)

Now let us fix $k\ge1$. We apply the same procedure to 
$
\sum_{n=0}^\infty {S_k(n)}\big/{\binom{8kn}{4kn}(-4)^{kn}} ,
$
where $S_k(n)$ is (hopefully) some (unknown) polynomial in $n$.
The beta integral evaluation gives
$$
\frac 1{\binom{8kn}{4kn}}=(8kn+1)\int_0^1x^{4kn}(1-x)^{4kn}dx .
$$
Hence, if $S_k(n)$ should have degree $d$ in $n$,
\begin{align}\notag
\sum _{n=0} ^{\infty}\frac {S_k(n)}{\binom{8kn}{4kn}(-4)^{kn}} 
&=\int_0^1\sum_{n=0}^\infty (8kn+1)S_k(n)\left(\frac 
{x^{4k}(1-x)^{4k}}{(-4)^k}\right)^n\,dx
\notag\\
&=\int_0^1 \frac {P_k(x)} {\(x^{4k}(1-x)^{4k}-(-4)^k\)^{d+2}}\,dx,
\label{eq:sumint}
\end{align}
where $P_k(x)$ is some polynomial in $x$. For convenience, let us
write $P$ as a short-hand for $P_k$.
Let $Q(x):=x^{4k}(1-x)^{4k}-(-4)^k$. Now we make the wild assumption
that
$$
\int \frac{P(x)}{Q(x)^{d+2}}\,dx=\frac{R(x)}{Q(x)^{d+1}}+2\arctan (x)+2\arctan
(x-1) ,
$$
for some polynomial $R(x)$ with $R(0)=R(1)=0$. Then the original sum
would indeed be equal to $\pi$. The last equality is equivalent to
$$
\frac P{Q^{d+2}}=\frac{R^{\prime }}{Q^{d+1}}-(d+1)\frac{Q^{\prime }R}{Q^{d+2}%
}+2\left(\frac 1{x^2+1}+\frac 1{x^2-2x+2}\right), 
$$
or
$$QR^{\prime }-(d+1)Q^{\prime }R=P-2Q^{d+2}\left(\frac 1{x^2+1}+
\frac 1{x^2-2x+2}\right) .
$$
In our examples, we observed that
$$
R(x)=(2x-1)\check{R}\big(x(1-x)\big) 
$$
for a polynomial $\check R$. So, let us make the substitution
$$t=x(1-x).$$
Then, after some simplification, the above differential equation
becomes
\begin{equation} \label{eq:diff}
-(1-4t)Q\frac{d\check{R}}{dt}+(2Q+4k(4k+1)(1-4t)t^{4k-1})\check{R}-P+2(3-2t)%
\frac{Q^{4k+2}}{t^2-2t+2}=0 ,
\end{equation}
where $Q(t)=t^{4k}-(-4)^{k}$.

Now, writing $N(k)=4k(4k+1)$, we make the Ansatz
\begin{align*}
\check R(t)&=\sum _{j=1} ^{N(k)-1}a(j)t^j,\\
S_k(n)&=\sum _{j=0} ^{4k}a(N(k)+j)n^j.
\end{align*}
(The reader should recall that $S_k(n)$ defines $P_k(t)=P(t)$
through \eqref{eq:sumint}.) Comparing
coefficients of powers of $t$ on both sides of \eqref{eq:diff}, we
get a system of $N(k)+4k$ linear equations for the unknowns
$a(1),a(2),\dots,a(N(k)+4k)$.

Hence: {\it If the determinant of this system
of linear equations is non-zero, then
there does indeed exist a representation
$$\pi=\sum_{n=0}^\infty \frac {S_k(n)}{\binom{8kn}{4kn}(-4)^n}.$$
}

To see whether we could indeed hope for the determinant to be
non-zero, we went again to the computer and looked at the values of
the determinant in some small
instances. (Obviously, we do not want to do this by
hand, since for $k=1$ the matrix is already a $24\times 24$ matrix!)
So, let us program the matrix. (We shall see
the mathematical definition of the matrix in just a
moment, see \eqref{eq:Det}.\footnote{To tell the truth, this is the form of the
matrix after some simplifications have already been carried out.
(In particular, we are looking at a matrix which is slightly smaller
than the original one.)
See 
\machSeite{AlKPAA}\cite[beginning of Section~4]{AlKPAA} for these details.
There, the matrix in \eqref{eq:Det} is called $M'''$.})

\MATH
\goodbreakpoint%
In[8]:= a[k\MATHtief ,j\MATHtief ]:=Module[%
\MATHlbrace Var=j/(4k)%
\MATHrbrace ,
\leavevmode%
                    (-1)\MATHhoch (Var-1)*8k(4k+1)(-4)\MATHhoch (k*(Var+1))*
\leavevmode%
                    Product[4k*l-1,%
\MATHlbrace l,1,4k-Var%
\MATHrbrace ]*Product[4k*l+1,%
\MATHlbrace l,1,Var-1%
\MATHrbrace ]
\leavevmode%
                    ]

In[9]:= A[k\MATHtief , i\MATHtief , j\MATHtief ] := Module[%
\MATHlbrace Var%
\MATHrbrace , 
\leavevmode%
        Var = %
\MATHlbrace Floor[(i - 2)/(4*k - 1)], 
\leavevmode%
            Floor[(j - 1)/(4*k)], Mod[i - 2, 4*k - 1], 
\leavevmode%
            Mod[j - 1, 4*k]%
\MATHrbrace ; 
\leavevmode%
        If[i == 1, 
\leavevmode%
          If[Mod[j, 4*k] === 0, a[k, j], 0], 
\leavevmode%
          If[Var[[1]] - Var[[2]] == 0, 
\leavevmode%
            Switch[Var[[3]] - Var[[4]], 0, f1[k, Var[[3]]+1, j], -1, 
\leavevmode%
              f0[k, Var[[3]]+1, j], \MATHtief , 0], 
\leavevmode%
            If[Var[[1]] - Var[[2]] == 1, 
\leavevmode%
              Switch[Var[[3]] - Var[[4]], 0, g1[k, Var[[3]]+1, j], -1, 
\leavevmode%
                g0[k, Var[[3]]+1, j], \MATHtief , 0], 0]]]]
\goodbreakpoint%
In[10]:= A[k\MATHtief ] := Table[A[k, i, j], %
\MATHlbrace i, 1, 16*k\MATHhoch 2%
\MATHrbrace , %
\MATHlbrace j, 1, 16*k\MATHhoch 2%
\MATHrbrace ]
\goodbreakpoint%
In[11]:= f0[k\MATHtief , t\MATHtief , j\MATHtief ] := j*(-4)\MATHhoch k;
\leavevmode%
         f1[k\MATHtief , t\MATHtief , j\MATHtief ] := -(2 + 4*j)*(-4)\MATHhoch k;
\leavevmode%
         g0[k\MATHtief , t\MATHtief , j\MATHtief ] := (4*k*(4*k + 1) - j);
\leavevmode%
         g1[k\MATHtief , t\MATHtief , j\MATHtief ] := (-4*4*k*(4*k + 1) + 2 + 4*j)
\goodbreakpoint%
\endMATH

We shall not try to digest this at this point. Let us accept the
program as a black box, and let us compute the determinant for $k=2$.

\MATH
\goodbreakpoint%
In[12]:= Det[A[2]]
\goodbreakpoint%
Out[12]= -601576375580370166777074138698518196031142518971568946712\MATHbackslash 
 
\MATHgroesser     2204136674781038302774231725971306459064075121023092662279814\MATHbackslash 
 
\MATHgroesser     015195545600000000000
\goodbreakpoint%
\endMATH
Magnificent! This is certainly {\it not\/} zero. However, what are we going to
do with this gigantic number? Remembering our discussion about
``nice" numbers and ``nice" formulae in the preceding section, 
let us factorise it in its prime factors.

\MATH
\goodbreakpoint%
In[13]:= FactorInteger[\%]
\goodbreakpoint%
Out[13]= %
\MATHlbrace %
\MATHlbrace -1, 1%
\MATHrbrace , %
\MATHlbrace 2, 325%
\MATHrbrace , %
\MATHlbrace 3, 39%
\MATHrbrace , %
\MATHlbrace 5, 11%
\MATHrbrace , %
\MATHlbrace 7, 11%
\MATHrbrace , %
\MATHlbrace 11, 3%
\MATHrbrace , %
\MATHlbrace 13, 2%
\MATHrbrace %
\MATHrbrace 
\goodbreakpoint%
\endMATH
I would say that this is sensational: a number with 139 digits, 
and the biggest prime factor is 13! As a matter of fact, this is not just 
a rare exception. Table~\ref{tab:1} shows the factorisations of the first five
determinants. (We could not go further because of the exploding size
of the matrix of which the determinant is taken.)

\begin{table}[h]
\vskip10pt
\begin{tabular}{l|l}
\hphantom{1}$k$ & $\hphantom{-}\det(A(k))$ \\ 
\hline\\[-8pt]
\hphantom{1}1 & $\hphantom{-}2^{59}3^55^67^1$ \\ 
\hphantom{1}2 & $-2^{325}3^{39}5^{11}7^{11}11^313^2$ \\ 
\hphantom{1}3 & $\hphantom{-}2^{772}3^{146}5^{28}7^{17}11^{17}13^{18}17^419^323^1$ \\ 
\hphantom{1}4 & $-2^{1913}3^{111}5^{58}7^{38}11^{21}13^{22}17^{24}19^723^529^231^1$ \\ 
\hphantom{1}5 & $\hphantom{-}
2^{2932}3^{202}5^{306}7^{69}11^{29}13^{27}17^{28}19^{29}23^{9}29^631^537^2$%
\end{tabular}
\vskip10pt
\caption{}
\label{tab:1}
\end{table}

Thus, these
experimental results {\it make us sure} that there must be a ``nice"
formula for the determinant. Indeed, we prove in 
\machSeite{AlKPAA}\cite{AlKPAA} 
that\footnote{Strictly speaking, this is not a formula 
$NICE(k)$ according to my ``Definition" in the preceding section,
because of the presence of the ``Abelian" factors $k^{8k^2+2k}$ and
$(4k+1)!^{4k}$, see Footnote~\ref{foot:normally}. Nevertheless, the
reader will certainly admit that this is a {\it nice} and {\it
closed\/} formula.}
\begin{equation} \label{eq:det(A(k))}
\det (A(k))=(-1)^{k-1}
2^{16k^3+20k^2+6k}k^{8k^2+2k}(4k+1)!^{4k}
\prod_{j=1}^{4k}\frac{(2j)!}{j!^2}.
\end{equation}

Hence the desired theorem follows.

\begin{Theorem} \label{T1}
For all $k\geq 1$ there is a formula
$$
\pi =\sum_{n=0}^\infty \frac{S_k(n)}{\binom{8kn}{4kn}(-4)^{kn}},
$$
where $S_k(n)$ is a polynomial in $n$ 
of degree $4k$ with rational coefficients. The polynomial $S_k(n)$ can be
found by solving the previously described system of linear equations.
\end{Theorem}

I must admit that we were extremely lucky that
it was indeed possible to {\it evaluate} the determinant {\it
explicitly}. To recall,  
``all" we needed to prove our theorem (Theorem~\ref{T1}) 
was to show that the determinant was
{\it non-zero}. To be honest, I would not have the slightest idea how to
do this here without finding the exact value of the determinant.

\medskip
Now, after all this somewhat ``dry" discussion, let me present the
determinant. 
We had to determine the
determinant of the $16k^2\times 16k^2$ matrix
\begin{equation} \label{eq:Det}
\begin{pmatrix} 
0\dots0\,*&
0\dots0\,*&
0\dots0\,*&
\dots&
\dots&
\dots&
0\dots0\,*\\
\hbox{\Large$F_1$}&\hbox{\Large$0$}&\hbox{\Large$0$}&
\dots&\dots&\dots&\hbox{\Large$0$}\\
\hbox{\Large$G_1$}&\hbox{\Large$F_2$}& \hbox{\Large$0$}&
\dots&\dots&\dots&\hbox{\Large$0$}\\
\hbox{\Large$0$}&\hbox{\Large$G_2$}&\hbox{\Large$F_3$}& &&&\vdots\\
\hbox{\Large$0$}&\hbox{\Large$0$}&\hbox{\Large$G_3$}& &
\ddots& &\vdots\\
\vdots&\ddots&\ddots&\ddots&\ddots&\ddots&\vdots\\
\vdots&&\ddots&\ddots&\ddots&\hbox{\Large$F_{4k-1}$}&\hbox{\Large$0$}\\
\vdots&&&\hbox{\Large$0$}&\hbox{\Large$0$}&
\hbox{\Large$G_{4k-1}$}&\hbox{\Large$F_{4k}$}\\
\hbox{\Large$0$}&\dots&\dots&\dots&\hbox{\Large$0$}&\hbox{\Large$0$}&\hbox{\Large$G_{4k}$}
\end{pmatrix},
\end{equation}
where the $\ell$-th non-zero entry in the first row
(these are marked by $*$) is
$$
(-1)^{\ell -1}(-4)^{(\ell +1)k}8k(4k+1)\left(\prod _{i=1} ^{4k-\ell }(4ik-1)\right)
\left(\prod _{i=1} ^{\ell -1}(4ik+1)\right),
$$
and where each block $F_t$ and $G_t$ is a $(4k-1)\times(4k)$ matrix (that
is, these are {\it rectangular} blocks!) with non-zero entries only
on its (two) main diagonals,
$$
F_t =\left(\smallmatrix f_1(4(t -1)k+1)&f_0(4(t -1)k+2)&0&\dots\\
0&f_1(4(t -1)k+2)&f_0(4(t -1)k+3)&0&\dots\\
&\ddots&\ddots&\\
&&\ddots&\ddots&\\
&&0&f_1(4t k-2)&f_0(4t k-1)&0\\
&&&0&f_1(4t k-1)&f_0(4t k)\\
\endsmallmatrix\right)
$$
and
$$
G_t =\left(\smallmatrix g_1(4(t -1)k+1)&g_0(4(t -1)k+2)&0&\dots\\
0&g_1(4(t -1)k+2)&g_0(4(t -1)k+3)&0&\dots\\
&\ddots&\ddots&\\
&&\ddots&\ddots&\\
&&0&g_1(4t k-2)&g_0(4t k-1)&0\\
&&&0&g_1(4t k-1)&g_0(4t k)\\
\endsmallmatrix\right).
$$
We have almost worked our way through the definition of the
determinant. The only missing piece is the definition of the
functions $f_0,f_1,g_0,g_1$ in the blocks $F_t$ and $G_t$. Here
it is:
\begin{align} 
\notag
f_0(j)&=j(-4)^k,\\
\notag
f_1(j)&=-(4j+2)(-4)^k,\\
\notag
g_0(j)&=(N(k)-j),\\
\label{eq:fg} 
g_1(j)&=-(4N(k)-4j-2),
\end{align}
where, as before, we write $N(k)=4k(4k+1)$ for short.

\section{The evaluation of the determinant}
\label{sec:eval}

I now describe how the determinant of \eqref{eq:Det} was evaluated by
applying to it the methods described in 
\machSeite{KratBN}\cite{KratBN}. To make this
section as self-contained as possible, for each of them
I briefly recall how it works before putting it into action.

\medskip
{\it ``Method" 0}: {\it Do row and column operations until the determinant
reduces to something manageable.}

In fact, at a first glance, this does not look too bad. 
Our matrix \eqref{eq:Det}, of which we want to compute
the determinant and show that it is non-zero, is a very sparse
matrix. Moreover, it looks almost like a two-diagonal matrix.
It seems that one should be able to do a few row and
column manipulations and thus reduce the matrix to a matrix of a
simpler form of which we can evaluate the determinant. 

Well, we tried that. Unfortunately, the above impression is
deceiving. First of all, the diagonals of the blocks do not really
fit together to form diagonals which run from one end of the matrix
to the other. Second, there remains still the first row which does
not fit the pattern of the rest of the matrix.
So, whatever we did, we ended up nowhere.
Maybe we should try something more \hbox{sophisticated \dots}

\medskip
{\it Method 1} 
\machSeite{KratBN}\cite[Sec.~2.6]{KratBN}: {\it LU-factorisation}.
Suppose we are given a family of matrices $A(1), A(2),A(3),\dots$ of
which we want to compute the determinants. Suppose further that we
can write
$$A(k)\cdot U(k)=L(k),$$
where $U(k)$ is an upper triangular matrix with 1s on the diagonal,
and where $L(k)$ is a lower triangular matrix. Then, clearly,
$$\det(A(k))=\text {product of the diagonal entries of $L(k)$}.$$

But how do we find $U(k)$ and $L(k)$? We go to the computer, crank
out $U(k)$ and $L(k)$ for $k=1,2,3,\dots$, until we are able to make
a guess. Afterwards we prove the guess by proving the corresponding
identities.

Well, we programmed that, we stared at the output on the computer
screen, but we could not make any sense of it.

\medskip
{\it Method 2} 
\machSeite{KratBN}\cite[Sec.~2.3]{KratBN}: {\it Condensation}.
This is based on a determinant formula due to Jacobi (see 
\machSeite{BresAO}\cite[Ch.~4]{BresAO}
and 
\machSeite{KnutAF}\cite[Sec.~3]{KnutAF}). Let $A$ be an $n\times n$ matrix. Let
$A_{i_1,i_2,\dots,i_\ell}^{j_1,j_2,\dots,j_\ell}$ denote the
submatrix of $A$ in which 
rows $i_1,i_2,\dots,i_\ell$ and columns $j_1,j_2,\dots,j_\ell$ are 
omitted. Then there holds
\begin{equation} \label{eq:cond}
\hfill \det A\cdot \det A_{1,n}^{1,n}=\det A_{1}^{1}\cdot \det A_n^n-
\det A_1^n\cdot \det A_n^1.
\end{equation}
If we consider a family of matrices $A(1),A(2),\dots$, and if all the
consecutive minors of $A(n)$ belong to the same family, then this
allows one to give an inductive proof of a conjectured determinant
evaluation for $A(n)$.

Let me illustrate this by reproducing Amdeberhan's condensation proof
\machSeite{AmdeAD}\cite{AmdeAD} of \eqref{eq:Krat}. Let $M_n(x,y)$ denote the
determinant in \eqref{eq:Krat}. Then we have
\begin{align} \notag
\big(M_n(x,y)\big)_n^n&=M_{n-1}(x,y), \\
\notag
\big(M_n(x,y)\big)_1^1&=M_{n-1}(x+1,y+1), \\
\notag
\big(M_n(x,y)\big)_n^1&=M_{n-1}(x-1,y+2), \\
\notag
\big(M_n(x,y)\big)_1^n&=M_{n-1}(x+2,y-1), \\
\label{eq:minors}
\big(M_n(x,y)\big)_{1,n}^{1,n}&=M_{n-2}(x+1,y+1).
\end{align}
Thus, we know that Equation \eqref{eq:cond} is satisfied
with $A$ replaced by $M_n(x,y)$,
where the minors appearing in \eqref{eq:cond} are given by
\eqref{eq:minors}. This can be interpreted as a recurrence for
the sequence $\big(M_n(x,y)\big)_{n\ge0}$. 
Indeed, given $M_0(x,y)$ and $M_1(x,y)$, the
equation \eqref{eq:cond} determines $M_n(x,y)$ uniquely
for all $n\ge0$ (given that
$M_n(x,y)$ never vanishes).
Thus, since the right-hand side of \eqref{eq:Krat} is indeed never
zero, for the proof of \eqref{eq:Krat}
it suffices to check \eqref{eq:Krat} for
$n=0$ and $n=1$, and that the right-hand side of \eqref{eq:Krat}
also satisfies \eqref{eq:cond},
all of which is a routine task.

\medskip

Now, a short glance at the definition of our matrix \eqref{eq:Det} 
will convince us
quickly that application of 
this method to it is rather hopeless. For example, omission of
the first row already brings us outside of our family of matrices.
So, also this method is not much help to solve our problem, which is
really a pity, because it is the most painless of \hbox{all \dots}

\medskip
{\it Method 3} 
\machSeite{KratBN}\cite[Sec.~2.4]{KratBN}: {\it Identification of
factors}. In order to sketch the idea, let us quickly go through a
(standard) proof of the {\it Vandermonde determinant evaluation},
\begin{equation} \label{eq:Vandermonde}
\det_{1\le i,j\le n}\(X_i^{j-1}\)=\prod _{1\le i<j\le n}
^{}(X_j-X_i).
\end{equation}

\begin{proof}If $X_{i_1}=X_{i_2}$ with $i_1\ne i_2$, then the Vandermonde
determinant \eqref{eq:Vandermonde} certainly vanishes because in that
case two rows of the determinant are identical. Hence,
$(X_{i_1}-X_{i_2})$ divides the determinant as a polynomial 
in the $X_i$'s. But that means that the complete product $\prod
_{1\le i<j\le n} (X_j-X_i)$ (which is exactly the right-hand side
of \eqref{eq:Vandermonde}) must divide the determinant.

On the other hand, the determinant is a polynomial in the $X_i$'s of
degree at most $\binom n2$. Combined with the previous observation,
this implies that the determinant equals the right-hand side product
times, possibly, some constant. To compute the constant, compare
coefficients of $X_1^0X_2^1\cdots X_n^{n-1}$ on both sides of 
\eqref{eq:Vandermonde}. This completes the proof of
\eqref{eq:Vandermonde}.
\end{proof}

At this point, let us extract the essence of this proof. 
The basic steps are:
{\em 
\begin{enumerate}
\item[(S1)] Identification of factors
\item[(S2)] Determination of degree bound
\item[(S3)] Computation of the multiplicative constant.
\end{enumerate}
}

As I report in 
\machSeite{KratBN}\cite{KratBN}, this turns out to be an extremely
powerful method which has numerous applications. To given an idea of
the flavour of the method, I show a few steps when it is applied
to the determinant in
\eqref{eq:Krat} (ignoring the fact that we have already found a very
simple proof of its evaluation; see 
\machSeite{KratBD}\cite[proof of Theorem~10]{KratBD}
for the complete proof using the ``identification of factors" method). 

\medskip
To get started, we have to transform the assertion \eqref{eq:Krat} into
an assertion about polynomials. This is easily done, we just have to
factor 
$$(x+y+i-1)!/(x+2i)!/(y+2n-i-2)!$$
out of the $i$-th row of the determinant. If we subsequently cancel
common factors on both sides of \eqref{eq:Krat}, we arrive at the
equivalent assertion
\begin{multline} \label{eq:Krat1}
\det_{0\le i,j\le
n-1}\big((x+y+i)_{j}\,(x+2i-j+1)_j\,(y+2j-i+1)_{2n-2j-2}\big)\\
=\prod _{i=0}
^{n-1}\big(i!\,(y+2i+1)_{n-i-1}\,(2x+y+2i)_i\,(x+2y+2i)_i\big),
\end{multline}
where, as before, 
$(\alpha)_k$ is the standard notation for shifted factorials
(Pochhammer symbols) explained in the statement of Theorem~\ref{thm:xy}.

In order to apply the same idea as in the above evaluation of the
Vandermonde determinant, as a first step we have to show that the
right-hand side of \eqref{eq:Krat1} divides the determinant on the
left-hand side as a polynomial in $x$ and $y$. For example, we would
have to prove that $(x+2y+2n-2)$ 
(actually,
$(x+2y+2n-2)^{\fl{(n+1)/3}}$, we will come to that in a moment)
divides the determinant. Equivalently,
if we set $x=-2y-2n+2$ in the determinant, then it should vanish. How
could we prove that? Well, if it vanishes then there must be a linear
combination of the columns, or of the rows, that vanishes. Equivalently, for
$x=-2y-2n+2$ we find a vector in the kernel of the matrix in \eqref{eq:Krat1}, 
respectively of its transpose. More generally (and this addresses the
fact that
we actually want to prove that $(x+2y+2n-2)^{\fl{(n+1)/3}}$ divides the
determinant):

\medskip
{\em \leftskip1cm
\rightskip1cm
\noindent
For proving that $(x+2y+2n-2)^E$ divides the
determinant, we find $E$ linear independent vectors in the
kernel.
\par}
\medskip
\noindent
(For a formal justification that this does indeed suffice,
see Section~2 of 
\machSeite{KratBI}\cite{KratBI}, and in particular the Lemma in that
section.)

Okay, how is this done in practice? You go to your computer, crank
out these vectors in the kernel, for
$n=1,2,3,\dots$, and try to make a guess what they are in general.
To see how this works, let us do it in our example.
First of all, we program the kernel of the matrix in \eqref{eq:Krat1}
with $x=-2y-2n+2$ (again, we are using {\sl Mathematica}
here).\footnote{In the program, {\tt V[n]} represents the kernel,
which 
is clearly a vector space. In the computer output, it is given in
parametric form, the parameters being the {\tt c[i]}'s.} 

\MATH
\goodbreakpoint%
In[14]:= p=Pochhammer;
\leavevmode%
         m[i\MATHtief ,j\MATHtief ,n\MATHtief ]:=p[x+y+i, j]*p[y+2*j+1-i, 2*n-2*j-2]*p[x-j+1+2i,j];
\leavevmode%
         V[n\MATHtief ]:=(x=-2y-2n+2;
\leavevmode%
                 Var=Sum[c[j]*Table[m[i,j,n],\MATHlbrace i,0,n-1\MATHrbrace ],\MATHlbrace j,0,n-1\MATHrbrace ];
\leavevmode%
                 Var=Solve[Var==Table[0,\MATHlbrace n\MATHrbrace ],Table[c[i],\MATHlbrace i,0,n-1\MATHrbrace ]];
\leavevmode%
                 Factor[Table[c[i],\MATHlbrace i,0,n-1\MATHrbrace ]/.Var])
\goodbreakpoint%
\endMATH
What the computer gives is the following:

\MATH
\goodbreakpoint%
In[15]:= V[2]
 
Out[15]= \MATHlbrace \MATHlbrace -2 c[1], c[1]\MATHrbrace \MATHrbrace 
\goodbreakpoint%
In[16]:= V[3]
 
Out[16]= \MATHlbrace \MATHlbrace -2 c[2], -c[2], c[2]\MATHrbrace \MATHrbrace 
\goodbreakpoint%
In[17]:= V[4]
 
Out[17]= \MATHlbrace \MATHlbrace -2 c[3], -3 c[3], 0, c[3]\MATHrbrace \MATHrbrace 
\goodbreakpoint%
In[18]:= V[5]
 
Out[18]= \MATHlbrace \MATHlbrace -2 c[4], -5 c[4], -2 c[3] - c[4], c[3], c[4]\MATHrbrace \MATHrbrace 
\goodbreakpoint%
In[19]:= V[6]
 
Out[19]= \MATHlbrace \MATHlbrace -2 c[5], -7 c[5], -2 (c[4] + 2 c[5]), -c[4], c[4], c[5]\MATHrbrace \MATHrbrace 
\goodbreakpoint%
In[20]:= V[7]
 
Out[20]= \MATHlbrace \MATHlbrace -2 c[6], -9 c[6], -2 c[5] - 9 c[6], %
-3 c[5] - c[6], 0, %

>     c[5], c[6]\MATHrbrace \MATHrbrace 
\goodbreakpoint%
\endMATH
At this point, the computations become somewhat slow.
So we should help our computer. Indeed, on the basis of what we have
obtained so far, it is ``obvious" that, somewhat unexpectedly, $y$ does not
appear in the result. Therefore we simply set $y$ equal to some random
number, and then the computer can go much further without any effort.

\MATH
\goodbreakpoint%
In[21]:= y=101
\goodbreakpoint%
In[22]:= V[8]
 
Out[22]= \MATHlbrace \MATHlbrace -2 c[7], -11 c[7], -2 (c[6] + 8 c[7]), -5 (c[6] + c[7]),
  
>     -2 c[5] - c[6], c[5], c[6], c[7]\MATHrbrace \MATHrbrace 
\goodbreakpoint%
In[23]:= V[9]
 
Out[23]= \MATHlbrace \MATHlbrace -2 c[8], -13 c[8], -2 c[7] - 25 c[8], -7 (c[7] + 2 c[8]),
  
>     -2 c[6] - 4 c[7] - c[8], -c[6], c[6], c[7], c[8]\MATHrbrace \MATHrbrace 
\goodbreakpoint%
In[24]:= V[10]
 
Out[24]= \MATHlbrace \MATHlbrace -2 c[9], -15 c[9], -2 (c[8] + 18 c[9]), -3 (3 c[8] + 10 c[9]),
  
>     -2 c[7] - 9 c[8] - 6 c[9], -3 c[7] - c[8], 0, c[7], c[8], c[9]\MATHrbrace \MATHrbrace 
\goodbreakpoint%
In[25]:= V[11]
 
Out[25]= \MATHlbrace \MATHlbrace -2 c[10], -17 c[10], -2 c[9] - 49 c[10], -11 (c[9] + 5 c[10]),
  
>     -2 (c[8] + 8 c[9] + 10 c[10]), -5 c[8] - 5 c[9] - c[10],
  
>     -2 c[7] - c[8], c[7], c[8], c[9], c[10]\MATHrbrace \MATHrbrace 
\goodbreakpoint%
%
%
%
%
%
%
%
\endMATH

Let us extract some information out of these data. 
For convenience, we write $M_n$ for the
matrix in \eqref{eq:Krat1} in the sequel. For example, by
just looking at the coefficients of $c[n-1]$ appearing in $V[n]$, we
extract that

\begin{enumerate}
\item [] the vector $(-2,1)$ is in the kernel of $M_2$, 
\item [] the vector $(-2,-1,1)$ is in the kernel of $M_3$, 
\item [] the vector $(-2,-3,0,1)$ is in the kernel of $M_4$,
\item [] the vector $(-2,-5,-1,0,1)$ is in the kernel of $M_5$, 
\item [] the vector $(-2,-7,-4,0,0,1)$ is in the kernel of $M_6$, 
\item [] the vector $(-2,-9,-9,-1,0,0,1)$ is in the kernel of $M_7$, 
\item [] the vector $(-2,-11,-16,-5,0,0,0,1)$ is in the kernel of $M_8$, 
\item [] the vector $(-2,-13,-25,-14,-1,0,0,0,1)$ is in the kernel of $M_9$, 
\item [] the vector $(-2,-15,-36,-30,-6,0,0,0,0,1)$ is in the kernel of $M_{10}$, 
\item [] the vector $(-2,-17,-49,-55,-20,-1,0,0,0,0,1)$ is in the kernel of $M_{11}$. 
\end{enumerate}

Okay, now we have to make sense out of this.
Our vectors in the kernel have the following structure: first, there
are some negative numbers, then follow a few zeroes, and finally there
is a trailing 1. I believe that we do not have any problem to
guess what the zeroeth\footnote{The indexing convention in the matrix
in \eqref{eq:Krat1} of which the determinant is taken is that rows and
columns are indexed by $0,1,\dots,n-1$. We keep this convention here.}
or the first coordinate of our vector is. Since
the second coordinates are always negatives of 
squares, there is also no problem
there. What about the third coordinates? Starting with the vector for
$M_7$, these are $-1,-5,-14,-30,-55,\dots$. I guess, rather than
thinking hard, we should consult {\tt Rate} (see Footnote~\ref{foot:Rate}):

\MATH
\goodbreakpoint%
In[26]:= Rate[-1,-5,-14,-30,-55]

\leavevmode%
          -(i0 (1 + i0) (1 + 2 i0))
Out[26]= \MATHlbrace -------------------------\MATHrbrace
\leavevmode%
                      6
\goodbreakpoint%
\endMATH
After replacing {\tt i0} by $n-6$ (as we should), this becomes
$-(n-6)(n-5)(2n-11)/6$. An interesting feature of this formula is that
it also works well for $n=6$ and $n=5$. Equipped with this experience,
we let {\tt Rate} work out the fourth coordinate:

\MATH
\goodbreakpoint%
In[27]:= Rate[0,0,0,-1,-6,-20]

\leavevmode%
                               2
\leavevmode%
          -((-3 + i0) (-2 + i0)  (-1 + i0))
Out[27]= \MATHlbrace ---------------------------------\MATHrbrace
\leavevmode%
                         12
\goodbreakpoint%
\endMATH
After replacement of {\tt i0} by $n-5$, this is
$-(n-8)(n-7)^2(n-6)/12$. Let us summarise our results so far: the
first five coordinates of our vector in the kernel of $M_n$ are
\begin{multline*}
-2,\ -(2n-5),\ -\frac {(n-4)(2n-8)} {2}, 
\ -\frac {(n-6)(n-5)(2n-11)} {6},
\\ 
-\frac {(n-8)(n-7)(n-6)(2n-14)} {12}.
\end{multline*}
I would say, there is a clear pattern emerging: the $s$-th coordinate
is equal to
$$-\frac {(n-2s)_{s-1}\,(2n-3s-2)} {s!}=
-\frac {(2n-3s-2)} {(n-s-1)}\frac {(n-2s)_{s}} {s!}.$$
Denoting the above expression by $f(n,s)$, the vector
$$(f(n,0),f(n,1),\dots,f(n,n-2),1)$$
is apparently in the kernel of $M_n$ for $n\ge2$. To prove this, we
have to show that
\begin{multline*} 
\sum _{s=0} ^{n-2}\frac
{(2n-3s-2)} {(n-s-1)}\frac {(n-2s)_{s}} {(s)!}\\
\cdot(-y-2n+i+2)_{s}\,(-2y-2n+2i-s+3)_{s}\,
(y+2s-i+1)_{2n-2s-2}\\
=(-y-2n+i+2)_{n-1}\,(-2y-3n+2i+4)_{n-1}.
\end{multline*}
In 
\machSeite{KratBD}\cite{KratBD} it was argued that this identity follows from
a certain hypergeometric identity due to Singh 
\machSeite{SinVAA}\cite{SinVAA}. However,
for just having {\it some} proof of this identity, this careful
literature search was not necessary. In fact, nowadays, {\it once you
write down a binomial or hypergeometric identity, it is already
proved!} One simply puts the binomial/hypergeometric sum into the
{\it Gosper--Zeilberger algorithm}
(see 
\machSeite{PeWZAA}%
\machSeite{ZeilAP}%
\machSeite{ZeilAM}%
\machSeite{ZeilAV}%
\cite{PeWZAA,ZeilAP,ZeilAM,ZeilAV}), 
which outputs a recurrence for it, and
then the only task is to verify that the (conjectured) right-hand side
also satisfies the same recurrence, and to check the identity 
for sufficiently many
initial values (which one has already done anyway while producing the
conjecture).\footnote{As you may have suspected, this is again a little
bit oversimplified. But not much. The Gosper--Zeilberger algorithm
applies {\it always} to hypergeometric sums, and there are only very few
binomial sums where it does not apply. (For the sake of completeness, I
mention that there are also several algorithms available to deal with
multi-sums, see 
\machSeite{ChSaAA}%
\machSeite{WiZeAC}%
\cite{ChSaAA,WiZeAC}. These do, however, rather quickly
challenge the resources of today's computers.) {\em Maple} 
implementations written by Doron
Zeilberger are available from {\tt
http://www.math.temple.edu/\~{}zeilberg}, 
those written by Fr\'ed\'eric Chyzak are available from 
{\tt http://algo.inria.fr/chyzak/mgfun.html},
{\sl Mathematica} implementations 
written by Peter Paule, Axel Riese, Markus Schorn, Kurt Wegschaider
are available 
from {\tt
http://www.risc.uni-linz.ac.at/research/combinat/risc/software}.}

As I mentioned earlier, actually
we need more vectors in the kernel. However, this is not difficult.
Take a closer look, and you will see that the pattern
persists (set $c[n-1]=0$ in the vector for $V[n]$, etc.). 
It will take you no time to work
out a full-fledged conjecture for $\fl{(n+1)/3}$ linear independent
vectors in the kernel of $M_n$.

I do not want to go through Steps (S2) and (S3), that is, the degree
calculation and the computation of the constant. As it turns out, to
do this conveniently you need to introduce more variables in the
determinant in \eqref{eq:Krat1}. Once you do this, everything works
out very smoothly. I refer the reader to 
\machSeite{KratBD}\cite{KratBD} for these details.

\medskip
Now, let us come back to our determinant, the 
determinant of \eqref{eq:Det}, and apply
``identification of factors" to it.%
\footnote{What I describe in the sequel is, except for very few 
excursions that ended up in a dead end, and which are therefore omitted here, 
the way how the determinant evaluation was found.}
To begin with, here is bad news: ``identification of factors"
crucially requires the existence of indeterminates.
But, where are they in \eqref{eq:Det}? If we look at the definition
of the matrix \eqref{eq:Det}, which, in the end, depends on the
auxiliary functions $f_0,f_1,g_0,g_1$ defined in \eqref{eq:fg}, 
then we see that there are no
indeterminates at all. Everything is (integral) numbers. So,
to get even started, we need
to introduce indeterminates in a way such that the more general
determinant would still factor ``nicely." We do not have much
guidance. Maybe, since we already made the abbreviation
$N(k)=4k(4k+1)$, we should replace $N(k)$ by $X$? Okay, let us try
this, that is, let us put
\begin{align} 
\notag
f_0(j)&=j(-4)^k,\\
\notag
f_1(j)&=-(4j+2)(-4)^k,\\
\notag
g_0(j)&=(X-j),\\
\label{eq:fg1} 
g_1(j)&=-(4X-4j-2)
\end{align}
instead of \eqref{eq:fg}. Let us compute the new determinant for
$k=2$. We program the new functions $f_0,f_1,g_0,g_1$,

\MATH
\goodbreakpoint%
In[28]:= f0[k\MATHtief , t\MATHtief , j\MATHtief ] := j*(-4)\MATHhoch k;
\leavevmode%
         f1[k\MATHtief , t\MATHtief , j\MATHtief ] := -(2 + 4*j)*(-4)\MATHhoch k;
\leavevmode%
         g0[k\MATHtief , t\MATHtief , j\MATHtief ] := (X - j);
\leavevmode%
         g1[k\MATHtief , t\MATHtief , j\MATHtief ] := (-4*X + 2 + 4*j)
\goodbreakpoint%
\endMATH
we enter the new determinant for $k=2$,

\MATH
\goodbreakpoint%
In[29]:= Factor[Det[A[2]]]
\goodbreakpoint%
\endMATH
and, after a waiting time of more than 15
minutes,\footnote{\label{foot:kompl}which I use
to explain why our computer needs so long to calculate this
determinant of a very sparse matrix of size $16\cdot 2^2=64$: isn't it
true that, nowadays, determinants of matrices with several hundreds of
rows and columns can be calculated without the slightest difficulty
(particularly if they are very sparse)?
Well, we should not forget that this is true for determinants of
matrices with {\it numerical\/} entries. However, our matrix
\eqref{eq:Det} with the modified definitions \eqref{eq:fg1} of
$f_0,f_1,g_0,g_1$ has now entries which are {\it polynomials} in
$X$. Hence, when our computer algebra program applies (internally)
some elimination algorithm to compute the determinant, huge rational
expressions will slowly build up and will slow down the calculations
(and, at times, will make our computer crash \dots). As I learn from
Dave Saunders, {\sl Maple} and {\sl Mathematica} do currently in fact
not use the best known algorithms for dealing with determinants of
matrices with polynomial entries. (This may have to do with the fact
that the developers try to optimise the algorithms for numerical
determinants in the first case.) It is known how to avoid the
expression swell and compute polynomial matrix determinants in time 
about $mn^3$, where $n$ is the dimension of the matrix 
and $m$ is the bit length of the determinant 
(roughly, in univariate case, $m$ is degree times
maximum coefficient length).}
we obtain

\MATH
\goodbreakpoint%
Out[29]= -1406399608474882323154910525986578515918369681041517636\MATHbackslash 

\leavevmode%
                                                              2
\MATHgroesser     11783762359972003840000000 (-64 + X) (-48 + X) (-40 + X) 

\leavevmode%
               3          4          5         6  7
\MATHgroesser     (-32 + X)  (-24 + X)  (-16 + X)  (-8 + X)  X 

\leavevmode%
                                                              2
\MATHgroesser     (9653078694297600 - 916000657637376 X + 36130368757760 X 

\leavevmode%
                       3               4             5           6
\MATHgroesser      - 758218948608 X  + 8928558848 X  - 55938432 X  + 145673 X )
\goodbreakpoint%
\endMATH
Not bad. There are many factors which are linear in $X$. (This is what
we were after.) However, the
irreducible polynomial of degree 6 gives us some headache. (The
degrees of the irreducible part of the polynomial 
will grow quickly with $k$.) How are we going to guess what this
factor could be, and, even more daunting, even if we should be
able to come up with a guess, how would we go about to prove it?

So, maybe we should modify our choice of how to introduce
indeterminates into the matrix. In fact, we overlooked something:
maybe, in a hidden manner, the variable $X$ is also there at other
places in \eqref{eq:fg}, that is, when $X$ is specialised to
$N(k)=4k(4k+1)$ at these places it becomes invisible. 
More specifically, maybe we should insert the
difference $X-4k(4k+1)$ in the definitions of $f_0$ and $f_1$ (which
would disappear for $X=4k(4k+1)$). So, maybe we should try:
\begin{align*} 
\notag
f_0(j)&=(4k(4k + 1) - X + j)(-4)^k,\\
\notag
f_1(j)&=-(16k(4k + 1) - 4X + 2 + 4j)(-4)^k,\\
\notag
g_0(j)&=(X-j),\\
g_1(j)&=-(4X-4j-2),
\end{align*}
Okay, let us modify our computer program accordingly,

\MATH
\goodbreakpoint%
In[30]:= f0[k\MATHtief , t\MATHtief , j\MATHtief ] := (4*k*(4*k + 1) - X + j)*(-4)\MATHhoch k;
\leavevmode%
         f1[k\MATHtief , t\MATHtief , j\MATHtief ] := -(4*4*k*(4*k + 1) - 4*X + 2 + 4*j)*(-4)\MATHhoch k;
\leavevmode%
         g0[k\MATHtief , t\MATHtief , j\MATHtief ] := (X - j);
\leavevmode%
         g1[k\MATHtief , t\MATHtief , j\MATHtief ] := (-4*X + 2 + 4*j)
\goodbreakpoint%
\endMATH
and let us compute the new determinant for $k=2$: 

\MATH
\goodbreakpoint%
In[31]:= Factor[Det[A[2]]]
\goodbreakpoint%
\endMATH
This makes us wait for
another 15 minutes, after which we are rewarded with:

\MATH
\goodbreakpoint%
Out[31]= -296777975397624679901369809794412104454134763494070841\MATHbackslash 

\MATHgroesser     1155365196124754770317472271790417634937439881166252558632\MATHbackslash 

\MATHgroesser     616674197504000000000 (-141 + 2 X) (-139 + 2 X) (-137 + 2 X) 

\MATHgroesser     (-135 + 2 X) (-133 + 2 X) (-131 + 2 X) (-129 + 2 X)
\goodbreakpoint%
\endMATH
Excellent! There is no big irreducible polynomial anymore. Everything
is linear factors in $X$. But, wait, there is still a problem: in the
end (recall Step~(S2)!) 
we will have to compare the degrees of the determinant and of the 
right-hand side as polynomials in $X$. If we expand the determinant
according to its definition, then the conclusion is that the degree
of the determinant is bounded above by $16k^2-1$, which, for $k=2$ is
equal to $31$. The right-hand side polynomial however which we
computed above has degree 7. This is a big gap!

I skip some other things (ending up in dead ends \dots) that we tried 
at this point. 
Altogether they pointed to the
fact that, apparently, {\it one} indeterminate is not sufficient. Perhaps it
is a good idea to ``diversify" the variable $X$, that is, to make two
variables, $X_1$ and $X_2$, out of $X$:
\begin{align*} 
\notag
f_0(j)&=(4k(4k + 1) - X_2 + j)(-4)^k,\\
\notag
f_1(j)&=-(16k(4k + 1) - 4X_1 + 2 + 4j)(-4)^k,\\
\notag
g_0(j)&=(X_2-j),\\
g_1(j)&=-(4X_1-4j-2).
\end{align*}
We program this,

\MATH
\goodbreakpoint%
In[32]:= f0[k\MATHtief , t\MATHtief , j\MATHtief ] := (4*k*(4*k + 1) - X[2] + j)*(-4)\MATHhoch k;
\leavevmode%
         f1[k\MATHtief , t\MATHtief , j\MATHtief ] := -(4*4*k*(4*k + 1) - 4*X[1] + 2 + 4*j)*(-4)\MATHhoch k;
\leavevmode%
         g0[k\MATHtief , t\MATHtief , j\MATHtief ] := (X[2] - j);
\leavevmode%
         g1[k\MATHtief , t\MATHtief , j\MATHtief ] := (-4*X[1] + 2 + 4*j)
\goodbreakpoint%
\endMATH
and, in order to avoid overstraining our computer, compute this
time the new determinant for $k=1$:

\MATH
\goodbreakpoint%
In[33]:= Factor[Det[A[1]]]
\goodbreakpoint%
\endMATH
After some minutes there appears

\MATH
\goodbreakpoint%
Out[33]= 3242591731706757120000 (-37 + 2 X[1]) (-35 + 2 X[1]) 

\leavevmode%
                                           3                      2
\MATHgroesser      (-33 + 2 X[1]) (1 + 2 X[1] - 2 X[2])  (3 + 2 X[1] - 2 X[2])  

\MATHgroesser      (5 + 2 X[1] - 2 X[2])
\goodbreakpoint%
\endMATH
on the computer screen.
On the positive side: the determinant still factors completely into
linear factors, something which we could not expect a priori.
Moreover, the degree (in $X_1$ and $X_2$) 
has increased, it is now equal to 9 although we were
only computing the determinant for $k=1$. However, a gap remains, the
degree should be $16k^2-1=15$ if $k=1$.

Thus, it may be wise to introduce another genuine variable, $Y$.
For example, we may think of simply homogenising the definitions of
$f_0,f_1,g_0,g_1$:
\begin{align*} 
\notag
f_0(j)&=(4k(4k + 1)Y - X_2 + jY)(-4)^k,\\
\notag
f_1(j)&=-(16k(4k + 1)Y - 4X_1 + (2 + 4j)Y)(-4)^k,\\
\notag
g_0(j)&=(X_2-jY),\\
g_1(j)&=-(4X_1-(4j+2)Y).
\end{align*}
We program this,

\MATH
\goodbreakpoint%
In[34]:= f0[k\MATHtief , t\MATHtief , j\MATHtief ] := (4*k*(4*k + 1)*Y - X[2] + j*Y)*(-4)\MATHhoch k;
\leavevmode%
         f1[k\MATHtief , t\MATHtief , j\MATHtief ] := -(4*4*k*(4*k + 1)*Y - 4*X[1] + 
\leavevmode%
                               (2 + 4*j)*Y)*(-4)\MATHhoch k;
\leavevmode%
         g0[k\MATHtief , t\MATHtief , j\MATHtief ] := (X[2] - j*Y);
\leavevmode%
         g1[k\MATHtief , t\MATHtief , j\MATHtief ] := (-4*X[1] + (2 + 4*j)*Y)
\goodbreakpoint%
In[35]:= Factor[Det[A[1]]]
\goodbreakpoint%
\endMATH
we wait for some more minutes, 
and we obtain

\MATH
\leavevmode%
                                  6
Out[35]= -3242591731706757120000 Y  (33 Y - 2 X[1]) (35 Y - 2 X[1]) 

\leavevmode%
                                            3                        2
\MATHgroesser      (37 Y - 2 X[1]) (Y + 2 X[1] - 2 X[2])  (3 Y + 2 X[1] - 2 X[2]) 

\MATHgroesser      (5 Y + 2 X[1] - 2 X[2])
\goodbreakpoint%
\endMATH
Great! The degree in $X_1,X_2,Y$ is 15, as it should
be!

At this point, one becomes greedy. The more variables we have, the
easier will be the proof. We ``diversify" the variables $X_1,X_2,Y$,
that is, we make them $X_{1,t},X_{2,t},Y_t$ if they appear in the
blocks $F_t$ or $G_t$, respectively, $t=1,2,\dots,4k$ (cf.\
\eqref{eq:Det} and the {\sl Mathematica} code for the precise meaning
of this definition):
\begin{align} 
\notag
f_0(j)&=(4k(4k + 1)Y_t - X_{2,t} + jY_t)(-4)^k,\\
\notag
f_1(j)&=-(16k(4k + 1)Y_t - 4X_{1,t} + (2 + 4j)Y_t)(-4)^k,\\
\notag
g_0(j)&=(X_{2,t}-jY_t),\\
\label{eq:fg5} 
g_1(j)&=-(4X_{1,t}-(4j+2)Y_t).
\end{align}

Now there are so many variables so that there is no way to do the
factorisation of the new determinant for $k=1$ on the computer unless
one plays tricks (which we\break did).%
\footnote{\label{foot:tricks}See
Footnote~\ref{foot:kompl} for the explanation of the complexity problem.
``Playing tricks" would mean to compute the
determinant for various special choices of the variables
$X_{1,t},X_{2,t},Y_t$, 
and then reconstruct the general result by interpolation.
This is possible because we know an a priori degree bound 
(namely $15$) for the polynomial.
However, this would become infeasible for $k=3$, for example. ``Playing
tricks" then would mean to be content with an 
``almost sure" guess, the latter being based on features of the (unknown)
general result that are already visible in the earlier results, 
and on calculations done for special values of the variables.
For example, if we encounter determinants $\det M(k)$, where the $M(k)$'s
are some square matrices, $k=1,2,\dots$, and the results for
$k=1,2,\dots,k_0-1$
show that $x-y$ must be a factor of $\det M(k)$ to some power, then
one would specialise $y$ to some value that would make $x-y$ distinct
from any other linear factors containing $x$, and, supposing that 
$y=17$ is such a choice, compute $\det M(k_0)$ with $y=17$. The exact
power of $x-y$ in the unspecialised determinant $\det M(k_0)$ 
can then be read off from
the exponent of $x-17$ in the specialised one. 
If it should happen that it is also infeasible
to calculate $\det M(k_0)$ with $x$ still unspecialised, then there is
still a way out. In that case, one specialises $y$ {\it and\/} $x$, in
such a way that $x-y$ would be a prime $p$ that one expects not to occur
as a prime factor in any other factor of the determinant $\det
M(k_0)$. The exact
power of $x-y$ in the unspecialised determinant $\det M(k_0)$
can then be read off from the exponent of $p$ in the prime
factorisation of the specialised determinant. See
Subsections~\ref{sec:signed} and \ref{sec:poset}, and in particular
Footnote~\ref{foot:maj} for further instances where this trick was applied.} 
But let us pretend that we are able to do it:

\MATH
\goodbreakpoint%
In[36]:= f0[k\MATHtief , t\MATHtief , j\MATHtief ] := (4*k*(4*k + 1)*Y[t] - X[2, t] 
\leavevmode%
                                       + j*Y[t])*(-4)\MATHhoch k;
\leavevmode%
         f1[k\MATHtief , t\MATHtief , j\MATHtief ] := -(4*4*k (4*k + 1)*Y[t] - 4*X[1, t] + 
\leavevmode%
                                 (2 + 4*j)*Y[t])*(-4)\MATHhoch k;
\leavevmode%
         g0[k\MATHtief , t\MATHtief , j\MATHtief ] := (X[2, t] - j*Y[t]);
\leavevmode%
         g1[k\MATHtief , t\MATHtief , j\MATHtief ] := (-4*X[1, t] + (2 + 4*j)*Y[t])
\goodbreakpoint%
In[37]:= Factor[Det[A[1]]]
\goodbreakpoint%
Out[37]= 3242591731706757120000 (2 X[1, 1] - 33 Y[1]) Y[1] 

\MATHgroesser      (2 X[1, 1] - 2 X[2, 1] + Y[1]) (2 X[1, 2] - 35 Y[2]) Y[2] 

\MATHgroesser      (2 X[1, 2] - 2 X[2, 2] + Y[2]) (-2 X[2, 2] Y[1] + 2 X[1, 1] 

\MATHgroesser       Y[2] + 3 Y[1] Y[2]) (2 X[1, 3] - 37 Y[3]) Y[3] (2 X[1, 3] - 

\MATHgroesser       2 X[2, 3] + Y[3]) (-2 X[2, 3] Y[1] + 2 X[1, 1] Y[3] + 

\MATHgroesser       5 Y[1] Y[3]) (-2 X[2, 3] Y[2] + 2 X[1, 2] Y[3] + 3 Y[2] Y[3])
\goodbreakpoint%
\endMATH

By staring a little bit at this result (and the one that we computed
for $k=2$), we extracted that,
apparently, we have
\begin{multline} \label{eq:CK4}
\det A^X=(-1)^{k-1}4^{2k(4k^2+7k+2)}k^{2k(4k+1)}
\prod _{i=1} ^{4k}(i+1)_{4k-i+1}\\
\times
\prod _{a=1}
^{4k-1}\left(2X_{1,a}-(32k^2+2a-1)Y_a\right)\\
\times
\prod _{1\le a\le b\le 4k-1} ^{}
(2X_{2,b}Y_a-2X_{1,a}Y_b-(2b-2a+1)Y_aY_b),
\end{multline}
where $A^X$ denotes the new general matrix given through
\eqref{eq:Det} and \eqref{eq:fg5}, and
where, as before, 
$(\alpha)_k$ is the standard notation for shifted factorials
(Pochhammer symbols) explained in the statement of Theorem~\ref{thm:xy}.
The special case that we need in the end to prove our Theorem~\ref{T1} is
$X_{1,t}=X_{2,t}=N(k)$ and $Y_t=1$.

\medskip
Now we are in business. Here is the {\it Sketch of the proof of
\eqref{eq:CK4}}:

\medskip
Re (S1): For each factor of the (conjectured) result \eqref{eq:CK4}, we
find a linear combination of the rows which vanishes if the factor
vanishes. (In other terms: if the indeterminates in the matrix are
specialised so that a particular factor vanishes, we find a vector in
the kernel of the transpose of the specialised matrix.) 
For example, to explain the factor
$(2X_{1,1}-(32k^2+1)Y_1)$, we found:

If $X_{1,1}=\frac {32k^2+1} {2}Y_1$, then
\begin{multline} \label{eq:combin} 
\frac {2(X_{2,4k-1}-(N(k)-1)Y_{4k-1})} {(-4)^{k(4k+1)+1}(16k^2+1)\prod
_{\ell=1} ^{4k-1}(4\ell k+1)}\cdot(\text {row 0 of $A^X$})
\\+
\sum _{s =0} ^{4k}\sum _{t =0} ^{4k-2}\Bigg(\frac {(-1)^{s (k-1)}2^t} {4^{s k}}
\prod _{\ell=0} ^{s -1}\frac {4k-1+4\ell k} {16k^2+1-4\ell k}
 \prod _{\ell=4k-t } ^{4k-1}\frac
{2X_{1,\ell}-(32k^2+2\ell-1)Y_\ell} 
{X_{2,\ell-1}-(16k^2+\ell-1)Y_{\ell-1}}\Bigg)\\
\cdot(\text {row $(16k^2-(4k-1)s -t -1)$ of $A^X$})=0,
\end{multline}
as is easy to verify. (Since the coefficients of the various rows in
\eqref{eq:combin} are rational functions in the indeterminates
$X_{1,t},X_{2,t},Y_t$, they
are rather easy to work out from computer data. 
One does not even need {\tt Rate} \dots)

\medskip
Re (S2): The total degree in the $X_{1,t}$'s, $X_{2,t}$'s, $Y_{t}$'s of the
product on the right-hand side of \eqref{eq:CK4} is $16k^2-1$. As we
already remarked earlier, the degree of the determinant is at most
$16k^2-1$. Hence, the determinant is equal to the product times,
possibly, a constant.

\medskip
Re (S3): For the evaluation of the constant, we compare coefficients of
$$
X_{1,1}^{4k}X_{1,2}^{4k-1}\cdots X_{1,4k-1}^2Y_1^1Y_2^2\cdots
Y_{4k-1}^{4k-1}.
$$
After some reflection,
it turns out that the constant is equal to a determinant of the same
form, that is, of the form \eqref{eq:Det}, but with auxiliary
functions 
\begin{align} 
\notag
f_0(j)&=(N(k)+j)(-4)^k,\\
\notag
f_1(j)&=4(-4)^k,\\
\notag
g_0(j)&=-j,\\
\label{eq:fg6} 
g_1(j)&=-4.
\end{align}

\medskip
What a set-back!
It seems that we are in the same situation as at the very beginning. 
We started with the determinant
of the matrix
\eqref{eq:Det} with auxiliary functions \eqref{eq:fg}, and we ended up
with the same type of determinant, with auxiliary functions
\eqref{eq:fg6}. There is little hope though: the functions in
\eqref{eq:fg6} are somewhat simpler as those in \eqref{eq:fg}.
Nevertheless, we have to play the same game again; that is, if we
want to apply the method of identification of factors, then we have
to introduce indeterminates. Skipping the experimental part, we came up
with
\begin{align*} f_0(j)&=(Z_t+j)(-4)^k,\\
f_1(j)&=4(-4)^kX_t,\\
g_0(j)&=-j,\\
g_1(j)&=-4X_t,
\end{align*}
where $t$ has the same meaning as before in \eqref{eq:fg5}.
Denoting the new matrix by $A^Z$, computer calculations suggested
that apparently
\begin{equation} \label{eq:CK9}
\det A^Z=(-1)^{k-1}2^{16k^3+20k^2+14k-1}k^{4k}(4k+1)!
\prod _{a=1} ^{4k-1}\Bigg(X_a^{4k+1-a}
\prod _{b=0} ^{a-1}(Z_a-4bk)\Bigg).
\end{equation}
The special case that we need is $Z_t=N(k)$ and $X_t=1$.

So, we apply again the method of identification of factors.
Everything runs smoothly (except that the details of the
verification of the factors are somewhat more unpleasant here). 
When we come finally to the point that we want to determine the
constant, it turns out that the constant is
equal to --- no surprise anymore --- 
the determinant of a matrix of the form \eqref{eq:Det} with
auxiliary functions
\begin{align*} f_0(j)&=(-4)^k,\\
f_1(j)&=4(-4)^k,\\
g_0(j)&=0,\\
g_1(j)&=-4.
\end{align*}
Now, is this good or bad news? In other words, while painfully working
through the steps of ``identification of factors," will we 
forever continue 
producing new determinants of the form \eqref{eq:Det}, which we must
again handle by the same method? To give it away: this is indeed {\it
very good} news. The function $g_0(j)$ vanishes identically!
It makes it possible that now Method~0 (= do some row and column
manipulations) works. (See 
\machSeite{AlKPAA}\cite{AlKPAA} for the details.) 
We are --- finally --- done with the proof of
\eqref{eq:CK4}, and, since the right-hand side {\it does not\/} vanish
for $X_{1,t}=X_{2,t}=N(k)$ and $Y_t=1$, with the proof of
Theorem~\ref{T1}!
%
\QED

\section{More determinant evaluations}\label{sec:detlist}

This section complements the list of known determinant evaluations
given in Section~3 of 
\machSeite{KratBN}\cite{KratBN}. I list here several determinant
evaluations which I believe are interesting or attractive
(and, in the ideal case, both), that have appeared
since 
\machSeite{KratBN}\cite{KratBN}, or that I failed to mention in 
\machSeite{KratBN}\cite{KratBN}.
I also include several conjectures and open problems, some of them old,
some of them new. As in 
\machSeite{KratBN}\cite{KratBN}, each evaluation is accompanied by some
remarks providing information on the context in which it arose.
Again, the selection of determinant evaluations presented 
reflects totally my taste, which must be blamed in the case
of any shortcomings. The order of presentation follows loosely the
order of presentation of determinants in 
\machSeite{KratBN}\cite{KratBN}.

\subsection{More basic determinant evaluations}
I begin with two determinant evaluations belonging to the category ``standard
determinants" (see Section~2.1 in 
\machSeite{KratBN}\cite{KratBN}). They are among
those which I missed to state in 
\machSeite{KratBN}\cite{KratBN}.
The reminder for inclusion here is the paper 
\machSeite{AmZeAB}\cite{AmZeAB}.
There, Amdeberhan and Zeilberger propose an {\it automated approach} towards
determinant evaluations via the condensation method (see ``Method~2" in
Section~\ref{sec:eval}). They provide a list of examples which
can be obtained in that way. As they remark at the end of the paper,
all of these are special cases of Lemma~5 in 
\machSeite{KratBN}\cite{KratBN}, with the
exception of three, namely Eqs.~(8)--(10) in 
\machSeite{AmZeAB}\cite{AmZeAB}. 
In their turn, two of them, namely (8) and (9), are special cases of
the following evaluation. (For (10), see Lemma~\ref{prop:AmZe} below.)

\begin{Lemma} \label{lem:AmZe}
Let $P(Z)$ be a polynomial in $Z$ of degree $n-1$ 
with leading coefficient $L$. Then
\begin{equation} \label{eq:D1}
\det_{1\le i,j\le n}\(P(X_i+Y_j)\)=L^n
\prod _{i=1} ^{n}\binom {n-1}i\prod _{1\le i<j\le n}
^{}(X_i-X_j)(Y_j-Y_i).
\end{equation}
\quad \quad \qed
\end{Lemma}

This lemma is easily proved along the lines of 
the standard proof of the Vandermonde determinant evaluation which
we recalled in Section~\ref{sec:eval} (see the proof of 
\eqref{eq:Vandermonde}) or by condensation.

A multiplicative version of Lemma~\ref{lem:AmZe} is the following.

\begin{Lemma} \label{lem:AmZe2}
Let $P(Z)=p_{n-1}Z^{n-1}+p_{n-2}Z^{n-2}+\dots+p_0$. Then
\begin{equation} \label{eq:D2}
\det_{1\le i,j\le n}\(P(X_iY_j)\)=
\prod _{i=0} ^{n-1}p_i\prod _{1\le i<j\le n}
^{}(X_i-X_j)(Y_i-Y_j).
\end{equation}
\quad \quad \qed
\end{Lemma}

On the other hand, identity~(10) from 
\machSeite{AmZeAB}\cite{AmZeAB} can be generalised to the following
Cauchy-type determinant evaluation. As all the identities from
\machSeite{AmZeAB}\cite{AmZeAB}, it can also be proved by the
condensation method.

\begin{Lemma} \label{prop:AmZe}
Let $a_0,a_1,\dots,a_{n-1}$, $c_0,c_1,\dots,c_{n-1}$, $b$, $x$ and $y$
be indeterminates. Then, for any positive integer $n$, there holds
\begin{multline} \label{eq:AmZe}
\det_{0\le i,j\le n-1}\(\frac{(x+a_i+c_j)(y+bi+c_j)}
    { (x+a_i+bi+c_j)}\)\\=b^{n-1}\,(n-1)!
     \(\binom n2b+(n-1)x+y+\sum_{i=1}^{n-1}a_i +
      \sum_{i=0}^{n-1}c_i\)\\
\times
\frac{\displaystyle\prod_{0\le i<j\le n-1}^{}(c_j-c_i)\
     \prod_{i=1}^{n-1}(y-x-a_i)
     \prod_{1\le i<j\le n-1}^{}((j-i)b-a_i+a_j)}
{\displaystyle\prod_{i=1}^{n-1}\prod_{j=0}^{n-1}
     {(x+a_i+bi+c_j)}}.
\end{multline}
\quad \quad \qed
\end{Lemma}

Speaking of Cauchy-type determinant evaluations, this brings us to a
whole family of such evaluations which were instrumental in
Kuperberg's recent advance
\machSeite{KupeAH}\cite{KupeAH}
on the enumeration of (symmetry classes of) {\it alternating sign matrices}.
The reason that determinants, and also Pfaffians, play an important
role in this context is due to Propp's discovery 
(described for the first time in
\machSeite{ElKLAB}\cite[Sec.~7]{ElKLAB} and exploited in 
\machSeite{KupeAD}%
\machSeite{KupeAH}\cite{KupeAD,KupeAH}) that alternating sign matrices 
are in bijection with
{\it configurations in the six vertex model}, and due to
determinant and Pfaffian formulae due to Izergin
\machSeite{IzerAA}\cite{IzerAA} and Kuperberg
\machSeite{KupeAH}\cite{KupeAH} for certain multivariable
partition functions of the six vertex model under various
boundary conditions. In many cases,
this leads to determinants which are, or are similar to,
{\it Cauchy's evaluation of the double alternant\/} (see
\machSeite{MuirAB}%
\cite[vol.~III, p.~311]{MuirAB} and \eqref{eq:Cauchy} below) 
or {\it Schur's Pfaffian version} 
\machSeite{SchuAA}\cite[pp.~226/227]{SchuAA}
of it (see \eqref{eq:Schur} below).

Let me recall that
the {\it Pfaffian} $\Pf(A)$ of a
skew-symmetric $(2n)\times(2n)$ matrix $A$ is defined by
\begin{equation} \label{eq:Pfaff} 
\Pf(A)=\sum _{\pi} ^{}(-1)^{c(\pi)}\prod _{(ij)\in \pi} ^{}A_{ij},
\end{equation}
where the sum is over all perfect matchings $\pi$ of the complete
graph on $2n$ vertices, where $c(\pi)$ is the {\em crossing
number} of $\pi$, and where the product is over all edges $(ij)$,
$i<j$, in the matching $\pi$ (see e.g.\ 
\machSeite{StemAE}\cite[Sec.~2]{StemAE}). 
What links Pfaffians
so closely to determinants is (aside from similarity of definitions)
the fact that the Pfaffian of a skew-symmetric matrix is, up to sign,
the square root of its determinant. That is,
$\det(A)=\Pf(A)^2$ for any skew-symmetric $(2n)\times(2n)$ matrix $A$
(cf.\ 
\machSeite{StemAE}\cite[Prop.~2.2]{StemAE}). See the corresponding remarks and
additional references in 
\machSeite{KratBN}\cite[Sec.~2.8]{KratBN}.

The following three theorems present the relevant evaluations.
They are Theorems~15--17 from
\machSeite{KupeAH}\cite{KupeAH}. All of them are proved using
identification of factors (see ``Method~3" in Section~\ref{sec:eval}). 
The results in Theorem~\ref{thm:Kup1}
contain whole sets of indeterminates, whereas the results in
Theorems~\ref{thm:Kup2} and \ref{thm:Kup3} 
only have two indeterminates $p$ and $q$, respectively three
indeterminates $p$, $q$ and $r$, in them.
Identity \eqref{eq:Kup4} is originally due to
Laksov, Lascoux and Thorup 
\machSeite{LaLTAA}\cite{LaLTAA} and Stembridge
\machSeite{StemAE}\cite{StemAE}, independently.
The reader must be warned that the statements in
\machSeite{KupeAH}\cite[Theorems~15--17]{KupeAH} are often blurred
by typos.

\begin{Theorem} \label{thm:Kup1}
Let $x_1,x_2,\dots$ and $y_1,y_2,\dots$ be indeterminates. Then, for
any positive integer $n$, there hold
\begin{equation}
\label{eq:Cauchy}
\det_{1\leq i,j\leq n}\left(
		\frac{1}{x_i+y_j} 
\right)  =
	\frac{\displaystyle
		\prod_{1\leq i<j\leq n}(x_i-x_j)(y_i-y_j)
	}{\displaystyle
		\prod_{1\leq i, j\leq n}(x_i+y_j)
	},
\end{equation}
\begin{multline}
\label{eq:Kup2}
\det_{1\le i,j\le n}\(\frac1{x_i+y_j} - \frac1{1+x_iy_j}\) 
= \frac{\displaystyle\prod_{1\le i<j\le n}
(1-x_ix_j)(1-y_iy_j)(x_j-x_i)(y_j-y_i)}
{\displaystyle\prod_{1\le i,j\le n} (x_i+y_j)(1+x_iy_j)}\\
           \times \prod_{i=1}^n (1-x_i)(1-y_i) ,
\end{multline}
\begin{equation} \label{eq:Schur}
\underset{1\le i,j\le 2n}\Pf\(\frac {x_i-x_j} {x_i+x_j}\)=
\prod _{1\le i<j\le 2n} ^{}\frac {x_i-x_j} {x_i+x_j}.
\end{equation}
\begin{equation}
\label{eq:Kup4}
\underset{1\le i,j\le 2n}\Pf\(\frac{x_i-x_j}{1-x_ix_j}\)
= \prod_{1\le i<j\le 2n} \frac{x_i-x_j}{1-x_ix_j}.
\end{equation}
\quad \quad \qed
\end{Theorem}

\begin{Theorem} \label{thm:Kup2}
Let $p$ and $q$ be indeterminates. Then, for
any positive integer $n$, there hold
\begin{equation} \label{eq:Kup5} 
\det_{1\le i,j\le n}\(\frac{q^{n+j-i}-q^{-(n+j-i)}}
{p^{n+j-i}-p^{-(n+j-i)}}\)
= \frac{\displaystyle\prod_{1\le i \ne j\le n} (p^{j-i}-p^{-(j-i)})
\prod_{1\le i,j\le n} (qp^{j-i}-q^{-1}p^{-(j-i)})}
    {\displaystyle\prod_{1\le i,j\le n} (p^{n+j-i}-p^{-(n+j-i)})},
\end{equation}
\begin{equation} \label{eq:Kup6} 
\det_{1\le i,j\le n}\(\frac{q^{j-i}+q^{-(j-i)}}{p^{j-i}+p^{-(j-i)}}\)
= (-1)^{\binom{n}2}
    \frac{\displaystyle
2^n\prod_{\substack{1\le i \ne j\le n \\ 2\mid j-i}}(p^{j-i}-p^{-j+i})
    \prod_{\substack{1\le i,j\le n \\ 2\nmid j-i}}(qp^{j-i}-q^{-1}p^{-j+i})}
    {\displaystyle\prod_{1\le i,j\le n} (p^{j-i}+p^{-j+i})} ,
\end{equation}
\begin{multline} \label{eq:Kup7} 
\det_{1\le i,i\le n}\(\frac{q^{n+j+i}-q^{-(n+j+i)}}{p^{n+j+i}-p^{-(n+j+i)}} -
\frac{q^{n+j-i}-q^{-(n+j-i)}}{p^{n+j-i}-p^{-(n+j-i)}}\) \\=
\frac{\prod_{1\le i<j\le 2n}(p^{j-i}-p^{-(j-i)})
\prod_{\displaystyle\substack{1\le i,j \le 2n+1 \\ 2|j}}
(qp^{j-i}-q^{-1}p^{-(j-i)})}
{\displaystyle\prod_{1\le i,j\le n} (p^{n+j-i}-p^{-(n+j-i)})
(p^{n+j+i}-p^{-(n+j+i)})},
\end{multline}
\begin{multline} \label{eq:Kup8} 
\det_{1\le i,j\le n}\(\frac{q^{j+i}+q^{-(j+i)}}{p^{j+i}+p^{-(j+i)}} -
\frac{\displaystyle q^{j-i}+q^{-j+i}}{p^{j-i}+p^{-(j-i)}}\)\\
=(-1)^{\binom n2}
 \frac{\displaystyle 2^n\prod_{1\le i<j \le n} (p^{2(j-i)}-p^{-2(j-i)})^2
\prod_{\substack{1\le i,j \le 2n+1 \\ 2\nmid i,\,2|j}} (qp^{j-i}-q^{-1}p^{-(j-i)})}
{\displaystyle\prod_{1\le i,j\le n} (p^{j-i}+p^{-(j-i)})(p^{j+i}+p^{-(j+i)})}.
\end{multline}
\quad \quad \qed
\end{Theorem}

\begin{Theorem} \label{thm:Kup3}
Let $p$, $q$, and $r$ be indeterminates. Then, for
any positive integer $n$, there hold
\begin{multline} \label{eq:Kup9}
\underset{1\le i,j\le 2n}\Pf \(\frac{(q^{j-i}-q^{-(j-i)})(r^{j-i}-r^{-(j-i)})}
{(p^{j-i}-p^{-(j-i)})}\) \\
= \frac{\displaystyle\prod_{1\le i<j\le n} (p^{j-i}-p^{-(j-i)})^2
\prod_{1\le i,j\le n}
    (qp^{j-i}-q^{-1}p^{-(j-i)})(rp^{j-i}-r^{-1}p^{-(j-i)})}
{\displaystyle\prod_{1\le i,j \le n} (p^{n+j-i}-p^{-(n+j-i)})} 
\end{multline}
\begin{multline} \label{eq:Kup10}
\underset{1\le i,j\le 2n}\Pf \Bigg((p^{j+i}-p^{-(j+i)})(p^{j-i}-p^{-(j-i)})
\biggl(\frac{q^{j+i}-q^{-(j+i)}}{p^{j+i}-p^{-(j+i)}} -
\frac{q^{j-i}-q^{-(j-i)}}{p^{j-i}-p^{-(j-i)}}\biggr)\\
\cdot 
 \biggl(\frac{r^{j+i}-r^{-(j+i)}}
{p^{j+i}-p^{-(j+i)}} - 
\frac{r^{j-i}-r^{-(j-i)}}{p^{j-i}-p^{-(j-i)}}\biggr)\Bigg) \\
= \frac{\displaystyle\prod_{1\le i<j \le 2n}(p^{j-i}-p^{-(j-i)})
\prod_{\substack{1\le i,j \le 2n+1 \\ 2|j}}
(qp^{j-i}-q^{-1}p^{-(j-i)})(rp^{j-i}-r^{-1}p^{-(j-i)})}
{\displaystyle\prod_{1\le i<j \le 2n}(p^{j+i}-p^{-(j+i)})}.
\end{multline}
\quad \quad \qed
\end{Theorem}

Subsequent to Kuperberg's work,
Okada \machSeite{OkadAJ}\cite{OkadAJ}
related Kuperberg's determinants and Pfaffians to characters of
classical groups, by coming up with rather complex, but still
beautiful determinant identities. In particular, this allowed him to
settle one more of the conjectured enumeration formulae on symmetry
classes of alternating sign matrices.
Generalising even further, Ishikawa, Okada, Tagawa and
Zeng
\machSeite{IsOTAA}\cite{IsOTAA} have found more such determinant 
identities. Putting them into the framework of certain special
representations of the symmetric group, Lascoux 
\machSeite{LascAT}\cite{LascAT} has clarified the mechanism which
gives rise to these identities.

\medskip
The next six determinant lemmas are corollaries of
{\it elliptic determinant evaluations}
due to Rosengren and Schlosser
\machSeite{RoScAC}\cite{RoScAC}. 
(The latter will be addressed later in
Subsection~\ref{sec:ell}.) They partly extend the fundamental
determinant lemmas in
\machSeite{KratBN}\cite[Sec.~2.2]{KratBN}.
For the statements of the lemmas, we need the notion
of a {\it norm} of a polynomial
$a_0+a_1z+\dots+a_kz^k$, which we define to be the reciprocal of the
product of its roots, or, more explicitly,
as $(-1)^ka_k/a_0$.

If we specialise $p=0$ in Lemma~\ref{wp}, \eqref{awpi}, then we obtain a
determinant identity which generalises at the same time 
the Vandermonde determinant evaluation, Lemma~\ref{lem:AmZe} and
Lemma~\ref{lem:AmZe2}. 

\begin{Lemma} \label{lem:Van1}
Let $P_1,P_2,\dots,P_n$ be polynomials of degree $n$ and norm
$t$, given by 
$$P_j(x)=(-1)^nta_{j,0}x^n+
\sum _{k=0} ^{n-1}a_{j,k}x^k.$$
Then
\begin{equation} \label{eq:Van1} 
\det_{1\leq i,j\leq n}\left(P_j(x_i)\right)=(1-tx_1\dotsm x_n)
\bigg(\prod_{1\leq i<j\leq n}(x_j-x_i)\bigg)
\underset{0\le j\le n-1}{\det_{1\le i\le n}}(a_{i,j}).
\end{equation}
\quad \quad \qed
\end{Lemma}

Further determinant identities which generalise other {\it Weyl
denominator formulae} (cf.\ 
\machSeite{KratBN}\cite[Lemma~2]{KratBN}) could be obtained from
the special case $p=0$ of the other determinant evaluations in 
Lemma~\ref{wp}.

A generalisation of Lemma~6 from
\machSeite{KratBN}\cite{KratBN} in the same spirit 
can be obtained by setting $p=0$ in
Theorem~\ref{adet}. It is given as Corollary~5.1 in
\machSeite{RoScAC}\cite{RoScAC}.

\begin{Lemma} \label{lem:RS1}
Let $x_1,\dots,x_n$, $a_1,\dots,a_n$, and $t$ be indeterminates.
For each $j=1,\dots,n$, let $P_j$ be a polynomial
 of degree $j$ and norm $ta_1\dotsm a_j$. Then there holds
\begin{multline} \label{eq:RS1}
\det_{1\le i,j\le n}\left(P_j(x_i)
\prod_{k=j+1}^n(1-a_kx_i)\right)\\
=\frac{1-ta_1\dotsm a_nx_1\dotsm x_n}{1-t}
\prod_{i=1}^nP_i(1/a_i)
\prod_{1\le i<j\le n}a_j(x_j-x_i).
\end{multline}
\quad \quad \qed
\end{Lemma}

We continue with a consequence of Theorem~\ref{adetcor} (see Corollary~5.3
in \machSeite{RoScAC}\cite{RoScAC}).
The special case $P_{j-1}(x)=1$, $j=1,\dots,n$,
is Lemma~A.1 of 
\machSeite{SchlAB}\cite{SchlAB}, which was needed in order to obtain an
{\it $A_n$ matrix inversion} that played a crucial role in the derivation
of {\it multiple basic hypergeometric series identities}.
A slight generalisation was given in 
\machSeite{SchlAF}\cite[Lemma~A.1]{SchlAF}.

\begin{Lemma}\label{adetcorr}
Let $x_1,\dots,x_n$ and $b$ be indeterminates.
For each $j=1,\dots,n$, let $P_{j-1}(x)$ be a polynomial  in
$x$ of degree at most $j-1$ with constant term $1$,
and let $Q(x)=(1-y_1x)\dotsm (1-y_{n+1}x)$.
Then there holds
\begin{multline}\label{adetcorrid}
Q(b)\;\det_{1\le i,j\le n}\left(x_i^{n+1-j}P_{j-1}(x_i)
-b^{n+1-j}P_{j-1}(b)\frac{Q(x_i)}{Q(b)}\right)\\
=(1-bx_1\cdots x_ny_1\dotsm y_{n+1})
\prod_{i=1}^n(x_i-b)\prod_{1\le i<j\le n}(x_i-x_j).
\end{multline}
\quad \quad \qed
\end{Lemma}


Pairing the $(i,j)$-entry in the determinant in \eqref{eq:RS1} with
itself, but with $x_i$ replaced by $1/x_i$, one can construct another
determinant which evaluates in closed form. The result given below is
Corollary~5.5 in
\machSeite{RoScAC}\cite{RoScAC}. It is the special case
$p=0$ of Theorem~\ref{cdet}. 

\begin{Lemma}\label{cdetr}
Let $x_1,\dots,x_n$, $a_1,\dots,a_n$, and $c_1,\dots,c_{n+2}$ be
indeterminates. For each $j=1,\dots,n$, let $P_j$ be a polynomial
  of degree $j$ with norm
$(c_1\dotsm c_{n+2}a_{j+1}\dotsm a_n)^{-1}$.
Then there holds
\begin{multline}\label{cdetrid}
\det_{1\leq i,j\leq n}\left(x_i^{-n-1}
\prod_{k=1}^{n+2}(1-c_kx_i)\,
P_j(x_i)\prod_{k=j+1}^n(1-a_kx_i)\right.\\
\left.-x_i^{n+1}\prod_{k=1}^{n+2}(1-c_kx_i^{-1})\,
P_j(x_i^{-1})\prod_{k=j+1}^n(1-a_kx_i^{-1})\right)\\
=\frac{a_1\dotsm a_n}
{x_1\dotsm x_n\,(1-c_1\dotsm c_{n+2}a_1\dotsm a_n)}
\prod_{i=1}^nP_i(1/a_i)\\\times
\prod_{1\leq i<j\leq n+2}(1-c_ic_j)\prod_{i=1}^n(1-x_i^2)
\prod_{1\le i<j\le n}a_j(x_i-x_j)(1-1/x_ix_j).
\end{multline}
\quad \quad \qed
\end{Lemma}

It is worthwhile to state the limit case $c_{n+2}\to\infty$ of this
lemma separately, in which case the norm constraint on the polynomials
$P_j$ drops out, but, in return, the degree of $P_j$ gets lowered by
one (see 
\machSeite{RoScAC}\cite[Cor.~5.8]{RoScAC}).

\begin{Lemma}\label{cdetr1}
Let $x_1,\dots,x_n$, $a_2,\dots,a_n$, and $c_1,\dots,c_{n+1}$ be
indeterminates. For each $j=1,\dots,n$, let $P_{j-1}$ be a polynomial
  of degree at most $j-1$. Then there holds
\begin{multline}\label{cdetr1id}
\det_{1\leq i,j\leq n}\left(x_i^{-n}
\prod_{k=1}^{n+1}(1-c_kx_i)\,
P_{j-1}(x_i)\prod_{k=j+1}^n(1-a_kx_i)\right.\\
\left.-x_i^n\prod_{k=1}^{n+1}(1-c_kx_i^{-1})\,
P_{j-1}(x_i^{-1})\prod_{k=j+1}^n(1-a_kx_i^{-1})\right)\\
=\prod_{i=1}^nP_{i-1}(1/a_i)\prod_{1\leq i<j\leq n+1}(1-c_ic_j)\\\times
\prod_{i=1}^nx_i^{-1}(1-x_i^2)
\prod_{1\le i<j\le n}a_j(x_i-x_j)(1-1/x_ix_j).
\end{multline}
\quad \quad \qed
\end{Lemma}

Dividing both sides of \eqref{cdetr1id} by $\prod_{i=2}^n a_i^{i-1}$
and then letting $a_i$ tend to $\infty$, $i=2,3,\dots,n$, 
we arrive at the determinant evaluation below (see 
\machSeite{RoScAC}\cite[Cor.~5.11]{RoScAC}).
Its special case $P_{j-1}(x)=1$, $j=1,\dots,n$,
is Lemma~A.11 of 
\machSeite{SchlAB}\cite{SchlAB}, needed there in order to obtain a
{\it $C_n$ matrix
inversion}, which was later applied in 
\machSeite{SchlAG}\cite{SchlAG} to derive
{\it multiple q-Abel and q-Rothe summations}.

\begin{Lemma}\label{cdetr1cor}
Let $x_1,\dots,x_n$, and $c_1,\dots,c_{n+1}$ be
indeterminates. For each $j=1,\dots,n$, let $P_{j-1}$ be a polynomial
 of degree at most $j-1$. Then there holds
\begin{multline}\label{cdetr1corid}
\det_{1\leq i,j\leq n}\left(x_i^{-j}
\prod_{k=1}^{n+1}(1-c_kx_i)\,P_{j-1}(x_i)
-x_i^j\prod_{k=1}^{n+1}(1-c_kx_i^{-1})\,P_{j-1}(x_i^{-1})\right)\\
=\prod_{i=1}^n P_{i-1}(0)
\prod_{1\leq i<j\leq n+1}(1-c_ic_j)
\prod_{i=1}^nx_i^{-1}(1-x_i^2)
\prod_{1\le i<j\le n}(x_j-x_i)(1-1/x_ix_j).
\end{multline}
\quad \quad \qed
\end{Lemma}

It is an attractive feature of this determinant identity that it
contains, at the same time, the {\it Weyl denominator formulae} for the
classical root systems $B_n$, $C_n$ and $D_n$ as special cases (cf.\ 
\machSeite{KratBN}\cite[Lemma~2]{KratBN}). This is seen by setting
$P_j(x)=1$ for all $j$, $c_1=c_2=\dots=c_{n-1}=0$, 
and then $c_n=0$, $c_{n+1}=-1$ for the type $B_n$ case, 
$c_n=c_{n+1}=0$ for the type $C_n$ case, 
and $c_n=1$, $c_{n+1}=-1$ for the type $D_n$ case, respectively. 

\medskip
A determinant which is of completely different type, but 
which also belongs to the category of basic determinant evaluations, is the
determinant of a matrix where only two (circular) diagonals are filled
with non-zero elements. It was applied with advantage in 
\machSeite{HaKrAA}\cite{HaKrAA}
to evaluate {\it Scott-type permanents}.

\begin{Lemma} \label{prop:2diag}
Let $n$ and $r$ be positive integers, $r\le n$, and
$x_1,x_2,\dots,x_n$, $y_1,y_2,\dots,y_n$ be indeterminates. Then, with
$d=\gcd(r,n)$, we have
\begin{multline}
\det\begin{pmatrix}
x_1&0&\dots&0&y_{n-r+1}&0&\\
0&x_2&0&&0&y_{n-r+2}&0\\
&&\ddots&&&&\ddots&0\\
0&&&&&&0&y_n\\
y_1&0&\\
0&y_2&0\\
&0&\ddots&0&&&\ddots&0\\
&&0&y_{n-r}&0&&0&x_n\end{pmatrix}\\
\hskip3.6cm=\prod _{i=1} ^{d}\bigg(\prod _{j=1} ^{n/d}x_{i+(j-1)d}-
(-1)^{n/d}\prod _{j=1} ^{n/d}y_{i+(j-1)d}\bigg).
\end{multline}
{\em(}I.e., in the matrix there are only nonzero entries along two
diagonals, one of which is a broken diagonal.{\em)}
\quad \quad \qed
\end{Lemma}

A further basic determinant evaluation which I missed to state in
\machSeite{KratBN}\cite{KratBN} is the evaluation of the
determinant of a {\it skew circulant matrix} attributed to Scott 
\machSeite{ScotAB}\cite{ScotAB} in
\machSeite{MuirAB}\cite[p.~356]{MuirAB}. It was in fact
recently used by Fulmek in 
\machSeite{FulmAF}\cite{FulmAF}
to find a closed form formula for the number of {\it non-intersecting
lattice paths with equally spaced starting and end points living
on a cylinder}, improving on earlier results by Forrester 
\machSeite{ForrAC}\cite{ForrAC} on the {\it vicious walker model\/}
in {\it statistical mechanics}, see 
\machSeite{FulmAF}\cite[Lemma~9]{FulmAF}.

\begin{Theorem} \label{thm:circulant1}
Let $n$ by a fixed positive integer, and let 
$a_0,a_1,\dots,a_{n-1}$ be indeterminates. Then
\begin{multline} \label{eq:circulant1}
\det\begin{pmatrix} a_0&a_1&a_2&\dots&a_{n-2}&a_{n-1}\\
-a_{n-1}&a_0&a_1&\dots&a_{n-3}&a_{n-2}\\
-a_{n-2}&-a_{n-1}&a_0&\dots&a_{n-4}&a_{n-3}\\
\hdotsfor6\\
-a_{1}&-a_2&-a_3&\dots&-a_{n-1}&a_{0}
\end{pmatrix}\\=\prod _{i=0}
^{n-1}(a_0+\om^{2i+1}a_1+\om^{2(2i+1)}a_2+\dots+\om^{(n-1)(2i+1)}a_{n-1}),
\end{multline}
where $\om$ is a primitive $(2n)$-th root of unity.\quad \quad \qed
\end{Theorem}

\subsection{More confluent determinants}
Here I continue the discussion from the beginning of Section~3 in 
\machSeite{KratBN}\cite[Theorems~20--24]{KratBN}. There I presented
determinant evaluations of matrices which, essentially, 
consist of several vertical strips, each of which is formed by
taking a certain column vector and gluing it together with its
derivative, its second derivative, etc., respectively by a similar
construction where the derivative is replaced by a difference or
$q$-difference operator. 

Since most of this subsection will be under the influence of the
so-called {\it ``$q$-disease'',}\footnote{\label{foot:q}The 
distinctive symptom of this
disease is to invariably raise the question ``Is there also a
$q$-analogue?" My epidemiological research on {\textsf MathSciNet} revealed
that, while basically non-existent during the 1970s, this disease slowly
spread during the 1980s, and then had a sharp increase around 1990,
when it jumped from about 20 papers per year published with the word
``$q$-analog$*$" in it to over 80 in 1995, and since then it has been
roughly stable at 60--70 papers per year. In its
simplest form, somebody who is infected by this disease takes a
combinatorial identity, and replaces every occurrence of a positive
integer $n$
by its {\it ``$q$-analogue"} $1+q+q^2+\dots+q^{n-1}$, inserts some powers of
$q$ at the appropriate places, and hopes that the result of these
manipulations would be again
an identity, thus constituting a ``$q$-analogue" of the original
equation. I refer the reader to the bible 
\machSeite{GaRaAA}\cite{GaRaAA} for a rich
source of $q$-identities, and for the right way to look at (most)
combinatorial $q$-identities. In another form, given a certain set of
objects of which one knows the exact number, one defines a {\it statistics}
stat on these objects and now tries to evaluate 
$\sum _{O\text{ an object}} ^{}q^{\stat(O)}$. For a very instructive
text following these lines see 
\machSeite{FoHaAL}\cite{FoHaAL}, with emphasis on the
objects being permutations. There is also an important third form of
the disease in which one works in the ring of polynomials in variables
$x,y,\dots$ with coefficients being rational functions in $q$, 
and in which some pairs of variables satisfy commutation
relations of the type $xy=qxy$. The study of such polynomial rings
and algebras is often motivated by {\it quantum groups} and {\it 
quantum algebras}.
The reader may want to consult 
\machSeite{KoorAG}\cite{KoorAG} to learn more about this
direction.
While my description did not make this clear,
the three described forms of the $q$-disease are indeed strongly
inter-related.}
we shall need the standard $q$-notations
$(a;q)_k$, denoting the {\em $q$-shifted factorial\/} and being given by 
$(a;q)_0:=1$ and
$$(a;q)_k:=(1-a)(1-aq)\cdots(1-aq^{k-1})$$ 
if $k$ is a positive integer, as well as
$\left[\begin{smallmatrix}\al\\k\end{smallmatrix}\right]_q$, denoting 
the {\em $q$-binomial coefficient\/} and being defined by
$\left[\begin{smallmatrix}\al\\k\end{smallmatrix}\right]_q=0$ if $k<0$,
$\left[\begin{smallmatrix}\al\\0\end{smallmatrix}\right]_q=1$, and
$$\begin{bmatrix} \al\\k\end{bmatrix}_q:=
\frac {(1-q^\al)(1-q^{\al-1})\cdots(1-q^{\al-k+1})}
{(1-q^k)(1-q^{k-1})\cdots(1-q)}$$
if $k$ is a positive integer.
Clearly we have $\lim_{q\to1}\[\smallmatrix \al\\k\endsmallmatrix\]_q=
\binom \al k$.

The first result that I present is a $q$-extension of the evaluation
of the {\it
confluent alternant\/} due to Schendel
\machSeite{ScheAA}\cite{ScheAA} (cf.\ 
\machSeite{KratBN}\cite[paragraph before Theorem~20]{KratBN}). In fact,
Theorem~23 of \machSeite{KratBN}\cite{KratBN} already provided a
$q$-extension of (a generalisation of) Schendel's formula.
However, in 
\machSeite{JohWAF}\cite[Theorem~1]{JohWAF}, Johnson found a
different $q$-extension. The theorem below is a slight generalisation
of it. (The theorem below reduces to Johnson's theorem if one puts
$C=0$. For $q=1$, the theorem below and 
\machSeite{KratBN}\cite[Theorem~23]{KratBN} become equivalent. To go
from one determinant to the other in this special case, one would have
to take a certain factor out of each column.)

\begin{Theorem} \label{thm:Johnson1}
Let $n$ be a non-negative integer, and
let $A_m(X)$ denote the $n\times m$ matrix
$$\(\begin{bmatrix} C+i\\i-j\end{bmatrix}_q
(X;q)_{i-j}\)_{0\le i\le n-1,\,0\le j\le m-1}.$$
Given a composition of $n$, $n=m_1+\dots+m_\ell$, there holds
\begin{multline} \label{eq:Johnson1}
\det_{n\times n}\big(A_{m_1}(X_1)\,A_{m_2}(X_2)\dots
A_{m_\ell}(X_\ell)\big)\\=
q^{\sum_{1\le i<j<k\le \ell}m_im_jm_k}
\prod _{1\le i<j\le \ell} ^{}
\prod _{g=1} ^{m_i}
\prod _{h=1} ^{m_j}
\dfrac{(q^{h-1}X_i-q^{g-1}X_j)
      (1-q^{C+g+h-1+\sum_{r=1}^{i-1}m_r)}}
     {(1-q^{g+h-1+\sum_{r=1}^{i-1}m_r)}}.
\end{multline}
\quad \quad \qed
\end{Theorem}

In 
\machSeite{JohWAF}\cite[Theorem~2]{JohWAF}, Johnson provides as well a
confluent $q$-extension of the evaluation of Cauchy's double
alternant \eqref{eq:Cauchy}.
Already the case $q=1$ seems to not have appeared in the literature earlier.
Here, I was not able to introduce an additional parameter (as, for
example, the $C$ in Theorem~\ref{thm:Johnson1}).

\begin{Theorem} \label{thm:Johnson2}
Let $n$ be a non-negative integer, and
let $A_m(X)$ denote the $n\times m$ matrix
$$\(\frac {1} {(Y_i-X)(Y_i-qX)(Y_i-q^2X)\cdots(Y_i-q^{j-1}X)}\)
_{1\le i\le n,\,1\le j\le m}.$$
Given a composition of $n$, $n=m_1+\dots+m_\ell$, there holds
\begin{multline} \label{eq:Johnson2}
\det_{n\times n}\big(A_{m_1}(X_1)\,A_{m_2}(X_2)\dots
A_{m_\ell}(X_\ell)\big)\\=
\frac{\displaystyle
\(\prod _{1\le i<j\le n} ^{}(Y_i-Y_j)\)
\(\prod _{1\le i<j\le \ell} ^{}
\prod _{g=1} ^{m_i}
\prod _{h=1} ^{m_j}
(q^{h-1}X_j-q^{g-1}X_i)
\)}
{\displaystyle
\prod _{i=1} ^{n}
\prod _{j=1} ^{\ell}(Y_i-X_j)(Y_i-qX_j)(Y_i-q^2X_j)\cdots(Y_i-q^{m_j-1}X_j)}.
\end{multline}
\quad \quad \qed
\end{Theorem}

A surprising mixture between the confluent alternant and a confluent
double alternant appears in
\machSeite{CiucAL}\cite[Theorem~A.1]{CiucAL}. Ciucu used it there in
order to prove a {\it Coulomb gas law} and a {\it superposition
principle} for the joint correlation of certain
collections of holes for the {\it rhombus tiling model on the triangular
lattice}. (His main result is in fact based on an even more general, and
more complex, determinant evaluation, see
\machSeite{CiucAL}\cite[Theorem~8.1]{CiucAL}.)

\begin{Theorem} \label{thm:Ciucu}
Let $s_1,s_2, \dots, s_m\geq1$ and 
$t_1, t_2,\dots, t_n\geq1$
be integers. 
Write $S=\sum_{i=1}^m s_i$, $T=\sum_{j=1}^n t_j$, and assume $S\geq T$. 
Let $x_1,x_2,\dots,x_m$ and $y_1,y_2,\dots,y_n$ be indeterminates. Define $N$ to be the $S\times S$ matrix
\begin{equation} \label{eq:Ciuc} 
N=\left[\begin{matrix} A&B \end{matrix}\right]
\end{equation}
whose blocks are given by
\begin{multline}
A=
\\
\left(\!
\begin{matrix}
{\scriptscriptstyle \frac{{\binom 0 0}}{y_1-x_1}}\!\!\!\!&
{\scriptscriptstyle \frac{{\binom 1 0}}{(y_1-x_1)^2}}\!\!\!\!&{\scriptscriptstyle \cdots}\!\!\!\!&
{\scriptscriptstyle \frac{{\binom {t_1-1} 0}}{(y_1-x_1)^{t_1}}}&
\ &
{\scriptscriptstyle \frac{{\binom 0 0}}{y_n-x_1}}\!\!\!\!&
{\scriptscriptstyle \frac{{\binom 1 0}}{(y_n-x_1)^2}}\!\!\!\!&{\scriptscriptstyle \cdots}\!\!\!\!&
{\scriptscriptstyle \frac{{\binom {t_n-1} 0}}{(y_n-x_1)^{t_n}}}
\\
{\scriptscriptstyle \frac{{\binom 1 1}}{(y_1-x_1)^2}}\!\!\!\!&
{\scriptscriptstyle \frac{{\binom 2 1}}{(y_1-x_1)^3}}\!\!\!\!&{\scriptscriptstyle \cdots}\!\!\!\!&
{\scriptscriptstyle \frac{{\binom {t_1} 1}}{(y_1-x_1)^{t_1+1}}}&
\ &
{\scriptscriptstyle \frac{{\binom 1 1}}{(y_n-x_1)^2}}\!\!\!\!&
{\scriptscriptstyle \frac{{\binom 2 1}}{(y_n-x_1)^3}}\!\!\!\!&{\scriptscriptstyle \cdots}\!\!\!\!&
{\scriptscriptstyle \frac{{\binom {t_n} 1}}{(y_n-x_1)^{t_n+1}}}
\\
{\scriptstyle\cdot}\!\!\!\!&
{\scriptstyle\cdot}\!\!\!\!&\ \!\!\!\!&
{\scriptstyle\cdot}&
\!\!\!\!\cdots\!\!\! &
{\scriptstyle\cdot}\!\!\!\!&
{\scriptstyle\cdot}\!\!\!\!&\ \!\!\!\!&
{\scriptstyle\cdot}
\\
{\scriptstyle\cdot}\!\!\!\!&
{\scriptstyle\cdot}\!\!\!\!&\ \!\!\!\!&
{\scriptstyle\cdot}&
\ &
{\scriptstyle\cdot}\!\!\!\!&
{\scriptstyle\cdot}\!\!\!\!&\ \!\!\!\!&
{\scriptstyle\cdot}
\\
{\scriptstyle\cdot}\!\!\!\!&
{\scriptstyle\cdot}\!\!\!\!&\ \!\!\!\!&
{\scriptstyle\cdot}&
\ &
{\scriptstyle\cdot}\!\!\!\!&
{\scriptstyle\cdot}\!\!\!\!&\ \!\!\!\!&
{\scriptstyle\cdot}
\\
{\scriptscriptstyle \frac{{\binom {s_1-1} {s_1-1}}}{(y_1-x_1)^{s_1}}}\!\!\!\!&
{\scriptscriptstyle \frac{{\binom {s_1} {s_1-1}}}{(y_1-x_1)^{s_1+1}}}\!\!\!\!&{\scriptscriptstyle \cdots}\!\!\!\!&
{\scriptscriptstyle \frac{{\binom {s_1+t_1-2} {s_1-1}}}{(y_1-x_1)^{s_1+t_1-1}}}&
\ &
{\scriptscriptstyle \frac{{\binom {s_1-1} {s_1-1}}}{(y_n-x_1)^{s_1}}}\!\!\!\!&
{\scriptscriptstyle \frac{{\binom {s_1} {s_1-1}}}{(y_n-x_1)^{s_1+1}}}\!\!\!\!&{\scriptscriptstyle \cdots}\!\!\!\!&
{\scriptscriptstyle \frac{{\binom {s_1+t_n-2} {s_1-1}}}{(y_n-x_1)^{s_1+t_n-1}}}
\\
\ \!\!\!\!&
\ \!\!\!\!&\cdot\!\!\!\!&
\ &
\ &
\ \!\!\!\!&
\ \!\!\!\!&\cdot\!\!\!\!&
\ 
\\
\ \!\!\!\!&
\ \!\!\!\!&\cdot\!\!\!\!&
\ &
\ &
\ \!\!\!\!&
\ \!\!\!\!&\cdot\!\!\!\!&
\ 
\\
\ \!\!\!\!&
\ \!\!\!\!&\cdot\!\!\!\!&
\ &
\ &
\ \!\!\!\!&
\ \!\!\!\!&\cdot\!\!\!\!&
\ 
\\
{\scriptscriptstyle \frac{{\binom 0 0}}{y_1-x_m}}\!\!\!\!&
{\scriptscriptstyle \frac{{\binom 1 0}}{(y_1-x_m)^2}}\!\!\!\!&{\scriptscriptstyle \cdots}\!\!\!\!&
{\scriptscriptstyle \frac{{\binom {t_1-1} 0}}{(y_1-x_m)^{t_1}}}&
\ &
{\scriptscriptstyle \frac{{\binom 0 0}}{y_n-x_m}}\!\!\!\!&
{\scriptscriptstyle \frac{{\binom 1 0}}{(y_n-x_m)^2}}\!\!\!\!&{\scriptscriptstyle \cdots}\!\!\!\!&
{\scriptscriptstyle \frac{{\binom {t_n-1} 0}}{(y_n-x_m)^{t_n}}}
\\
{\scriptscriptstyle \frac{{\binom 1 1}}{(y_1-x_m)^2}}\!\!\!\!&
{\scriptscriptstyle \frac{{\binom 2 1}}{(y_1-x_m)^3}}\!\!\!\!&{\scriptscriptstyle \cdots}\!\!\!\!&
{\scriptscriptstyle \frac{{\binom {t_1} 1}}{(y_1-x_m)^{t_1+1}}}&
\ &
{\scriptscriptstyle \frac{{\binom 1 1}}{(y_n-x_m)^2}}\!\!\!\!&
{\scriptscriptstyle \frac{{\binom 2 1}}{(y_n-x_m)^3}}\!\!\!\!&{\scriptscriptstyle \cdots}\!\!\!\!&
{\scriptscriptstyle \frac{{\binom {t_n} 1}}{(y_n-x_m)^{t_n+1}}}
\\
{\scriptstyle\cdot}\!\!\!\!&
{\scriptstyle\cdot}\!\!\!\!&\ \!\!\!\!&
{\scriptstyle\cdot}&
\!\!\!\!\cdots\!\!\! &
{\scriptstyle\cdot}\!\!\!\!&
{\scriptstyle\cdot}\!\!\!\!&\ \!\!\!\!&
{\scriptstyle\cdot}
\\
{\scriptstyle\cdot}\!\!\!\!&
{\scriptstyle\cdot}\!\!\!\!&\ \!\!\!\!&
{\scriptstyle\cdot}&
\ &
{\scriptstyle\cdot}\!\!\!\!&
{\scriptstyle\cdot}\!\!\!\!&\ \!\!\!\!&
{\scriptstyle\cdot}
\\
{\scriptstyle\cdot}\!\!\!\!&
{\scriptstyle\cdot}\!\!\!\!&\ \!\!\!\!&
{\scriptstyle\cdot}&
\ &
{\scriptstyle\cdot}\!\!\!\!&
{\scriptstyle\cdot}\!\!\!\!&\ \!\!\!\!&
{\scriptstyle\cdot}
\\
{\scriptscriptstyle \frac{{\binom {s_m-1} {s_m-1}}}{(y_1-x_m)^{s_m}}}\!\!\!\!&
{\scriptscriptstyle \frac{{\binom {s_m} {s_m-1}}}{(y_1-x_m)^{s_m+1}}}\!\!\!\!&{\scriptscriptstyle \cdots}\!\!\!\!&
{\scriptscriptstyle \frac{{\binom {s_m+t_1-2} {s_m-1}}}{(y_1-x_m)^{s_m+t_1-1}}}&
\ &
{\scriptscriptstyle \frac{{\binom {s_m-1} {s_m-1}}}{(y_n-x_m)^{s_m}}}\!\!\!\!&
{\scriptscriptstyle \frac{{\binom {s_m} {s_m-1}}}{(y_n-x_m)^{s_m+1}}}\!\!\!\!&{\scriptscriptstyle \cdots}\!\!\!\!&
{\scriptscriptstyle \frac{{\binom {s_m+t_n-2} {s_m-1}}}{(y_n-x_m)^{s_m+t_n-1}}}
\end{matrix}
\!\right)
\end{multline}
and
\begin{equation}
B=
\left(\begin{matrix}
{\scriptstyle {\binom 0 0}x_1^0}\!\!\!&
{\scriptstyle {\binom 1 0}x_1}\!\!\!&{\scriptstyle \cdots}\!\!\!&
{\scriptstyle {\binom {S-T-1} {0}}x_1^{S-T-1}}
\\
{\scriptstyle {\binom {0} {1}}x_1^{-1}}\!\!\!&
{\scriptstyle {\binom {1} {1}}x_1^0}\!\!\!&{\scriptstyle \cdots}\!\!\!&
{\scriptstyle {\binom {S-T-1} {1}}x_1^{S-T-2}}
\\
{\scriptstyle \cdot}\!\!\!&
{\scriptstyle \cdot}\!\!\!&{\scriptstyle \ }\!\!\!&
{\scriptstyle \cdot}
\\
{\scriptstyle \cdot}\!\!\!&
{\scriptstyle \cdot}\!\!\!&{\scriptstyle \ }\!\!\!&
{\scriptstyle \cdot}
\\
{\scriptstyle \cdot}\!\!\!&
{\scriptstyle \cdot}\!\!\!&{\scriptstyle \ }\!\!\!&
{\scriptstyle \cdot}
\\
{\scriptstyle {\binom {0} {s_1-1}}x_1^{1-s_1}}\!\!\!&
{\scriptstyle {\binom {1} {s_1-1}}x_1^{2-s_1}}\!\!\!&{\scriptstyle \cdots}\!\!\!&
{\scriptstyle {\binom {S-T-1} {s_1-1}}x_1^{S-T-s_1}}
\\
{\scriptstyle \ }\!\!\!&
{\scriptstyle \ }\!\!\!&\cdot \!\!\!&
{\scriptstyle \ }
\\
{\scriptstyle \ }\!\!\!&
{\scriptstyle \ }\!\!\!&\cdot \!\!\!&
{\scriptstyle \ }
\\
{\scriptstyle \ }\!\!\!&
{\scriptstyle \ }\!\!\!&\cdot \!\!\!&
{\scriptstyle \ }
\\
{\scriptstyle {\binom {0} {0}}x_m^0}\!\!\!&
{\scriptstyle {\binom {1} {0}}x_m}\!\!\!&{\scriptstyle \cdots}\!\!\!&
{\scriptstyle {\binom {S-T-1} {0}}x_m^{S-T-1}}
\\
{\scriptstyle {\binom {0} {1}}x_m^{-1}}\!\!\!&
{\scriptstyle {\binom {1} {1}}x_m^0}\!\!\!&{\scriptstyle \cdots}\!\!\!&
{\scriptstyle {\binom {S-T-1} {1}}x_m^{S-T-2}}
\\
{\scriptstyle \cdot}\!\!\!&
{\scriptstyle \cdot}\!\!\!&{\scriptstyle \ }\!\!\!&
{\scriptstyle \cdot}
\\
{\scriptstyle \cdot}\!\!\!&
{\scriptstyle \cdot}\!\!\!&{\scriptstyle \ }\!\!\!&
{\scriptstyle \cdot}
\\
{\scriptstyle \cdot}\!\!\!&
{\scriptstyle \cdot}\!\!\!&{\scriptstyle \ }\!\!\!&
{\scriptstyle \cdot}
\\
{\scriptstyle {\binom {0} {s_m-1}}x_m^{1-s_m}}\!\!\!&
{\scriptstyle {\binom {1} {s_m-1}}x_m^{2-s_m}}\!\!\!&{\scriptstyle \cdots}\!\!\!&
{\scriptstyle {\binom {S-T-1} {s_m-1}}x_m^{S-T-s_m}}
\end{matrix}\right).
\end{equation}
Then we have

\vbox{\noindent
\begin{equation}
\det N =\frac{\prod_{1\leq i<j\leq m}(x_j-x_i)^{s_is_j}\prod_{1\leq i<j\leq n}(y_i-y_j)^{t_it_j}}
{\prod_{i=1}^m\prod_{j=1}^n(y_j-x_i)^{s_it_j}}.
\end{equation}
\quad \quad \qed}
\end{Theorem}

This theorem generalises at the same time numerous previously obtained
determinant evaluations. It reduces of course to Cauchy's double
alternant when $m=n$ and $s_1=s_2=\dots=s_m=t_1=t_2=\dots=t_n=1$. 
(In that case, the submatrix $B$ is empty.)
It reduces to the confluent alternant for
$t_1=t_2=\dots=t_n=0$ (i.e., in the case where the submatrix $A$ is
empty). The case $m=rn$, $s_1=s_2=\dots=s_m=1$, $t_1=t_2=\dots=t_n=r$
is stated as an exercise in
\machSeite{MuirAD}\cite[Ex.~42, p.~360]{MuirAD}. Finally, a mixture
of the double alternant and the Vandermonde determinant appeared already
in \machSeite{HaKrAA}\cite[Theorem~(Cauchy+)]{HaKrAA} where it was
used to
evaluate {\it Scott-type permanents}. This mixture turns out to be the
special case $s_1=s_2=\dots=s_m=t_1=t_2=\dots=t_n=1$ (but not
necessarily $m=n$) of Theorem~\ref{thm:Ciucu}.

If $S=T$ (i.e., in the case where the submatrix $B$ is empty),
Theorem~\ref{thm:Ciucu} 
provides the evaluation of a confluent double alternant which
is different from the one in Theorem~\ref{thm:Johnson2} for $C=0$ and
$q=1$. While, for the general form of Theorem~\ref{thm:Ciucu}, I was
not able to find a $q$-analogue, I was able to find one for this
special case, that is, for the case where $B$ is empty. 
In view of the
fact that there are also
$q$-analogues for the other extreme
case where the submatrix $A$ is empty (namely Theorem~\ref{thm:Johnson1} and 
\machSeite{KratBN}\cite[Theorem~23]{KratBN}), I still suspect that a
$q$-analogue of the general form of Theorem~\ref{thm:Ciucu} should exist.

\begin{Theorem} \label{thm:Ciucu1}
Let $s_1, s_2,\dots, s_m\geq1$ and $t_1, t_2,\dots, t_n\geq1$
be integers such that $s_1+s_2+\dots+s_m=t_1+t_2+\dots+t_n$. 
Let $x_1,x_2,\dots,x_m$ and $y_1,y_2,\dots,y_n$ be indeterminates. 
Let $A$ be the matrix
\begin{multline}
A=
\\
\left(\!
\begin{matrix}
{\scriptscriptstyle \frac{{\qbinom 0 0}}{\coef{x_1,y_1}^1}}\!\!\!\!&
{\scriptscriptstyle \frac{{\qbinom 1 0}}{\coef{x_1,y_1}^2}}\!\!\!\!&{\scriptscriptstyle \cdots}\!\!\!\!&
{\scriptscriptstyle \frac{{\qbinom {t_1-1} 0}}{\coef{x_1,y_1}^{t_1}}}&
\ &
{\scriptscriptstyle \frac{{\qbinom 0 0}}{\coef{x_1,y_n}^1}}\!\!\!\!&
{\scriptscriptstyle \frac{{\qbinom 1 0}}{\coef{x_1,y_n}^2}}\!\!\!\!&{\scriptscriptstyle \cdots}\!\!\!\!&
{\scriptscriptstyle \frac{{\qbinom {t_n-1} 0}}{\coef{x_1,y_n}^{t_n}}}
\\
{\scriptscriptstyle \frac{{\qbinom 1 1}}{\coef{x_1,y_1}^2}}\!\!\!\!&
{\scriptscriptstyle \frac{{\qbinom 2 1}}{\coef{x_1,y_1}^3}}\!\!\!\!&{\scriptscriptstyle \cdots}\!\!\!\!&
{\scriptscriptstyle \frac{{\qbinom {t_1} 1}}{\coef{x_1,y_1}^{t_1+1}}}&
\ &
{\scriptscriptstyle \frac{{\qbinom 1 1}}{\coef{x_1,y_n}^2}}\!\!\!\!&
{\scriptscriptstyle \frac{{\qbinom 2 1}}{\coef{x_1,y_n}^3}}\!\!\!\!&{\scriptscriptstyle \cdots}\!\!\!\!&
{\scriptscriptstyle \frac{{\qbinom {t_n} 1}}{\coef{x_1,y_n}^{t_n+1}}}
\\
{\scriptstyle\cdot}\!\!\!\!&
{\scriptstyle\cdot}\!\!\!\!&\ \!\!\!\!&
{\scriptstyle\cdot}&
\!\!\!\!\cdots\!\!\! &
{\scriptstyle\cdot}\!\!\!\!&
{\scriptstyle\cdot}\!\!\!\!&\ \!\!\!\!&
{\scriptstyle\cdot}
\\
{\scriptstyle\cdot}\!\!\!\!&
{\scriptstyle\cdot}\!\!\!\!&\ \!\!\!\!&
{\scriptstyle\cdot}&
\ &
{\scriptstyle\cdot}\!\!\!\!&
{\scriptstyle\cdot}\!\!\!\!&\ \!\!\!\!&
{\scriptstyle\cdot}
\\
{\scriptstyle\cdot}\!\!\!\!&
{\scriptstyle\cdot}\!\!\!\!&\ \!\!\!\!&
{\scriptstyle\cdot}&
\ &
{\scriptstyle\cdot}\!\!\!\!&
{\scriptstyle\cdot}\!\!\!\!&\ \!\!\!\!&
{\scriptstyle\cdot}
\\
{\scriptscriptstyle \frac{{\qbinom {s_1-1} {s_1-1}}}{\coef{x_1,y_1}^{s_1}}}\!\!\!\!&
{\scriptscriptstyle \frac{{\qbinom {s_1} {s_1-1}}}{\coef{x_1,y_1}^{s_1+1}}}\!\!\!\!&{\scriptscriptstyle \cdots}\!\!\!\!&
{\scriptscriptstyle \frac{{\qbinom {s_1+t_1-2} {s_1-1}}}{\coef{x_1,y_1}^{s_1+t_1-1}}}&
\ &
{\scriptscriptstyle \frac{{\qbinom {s_1-1} {s_1-1}}}{\coef{x_1,y_n}^{s_1}}}\!\!\!\!&
{\scriptscriptstyle \frac{{\qbinom {s_1} {s_1-1}}}{\coef{x_1,y_n}^{s_1+1}}}\!\!\!\!&{\scriptscriptstyle \cdots}\!\!\!\!&
{\scriptscriptstyle \frac{{\qbinom {s_1+t_n-2} {s_1-1}}}{\coef{x_1,y_n}^{s_1+t_n-1}}}
\\
\ \!\!\!\!&
\ \!\!\!\!&\cdot\!\!\!\!&
\ &
\ &
\ \!\!\!\!&
\ \!\!\!\!&\cdot\!\!\!\!&
\ 
\\
\ \!\!\!\!&
\ \!\!\!\!&\cdot\!\!\!\!&
\ &
\ &
\ \!\!\!\!&
\ \!\!\!\!&\cdot\!\!\!\!&
\ 
\\
\ \!\!\!\!&
\ \!\!\!\!&\cdot\!\!\!\!&
\ &
\ &
\ \!\!\!\!&
\ \!\!\!\!&\cdot\!\!\!\!&
\ 
\\
{\scriptscriptstyle \frac{{\qbinom 0 0}}{\coef{x_m,y_1}^1}}\!\!\!\!&
{\scriptscriptstyle \frac{{\qbinom 1 0}}{\coef{x_m,y_1}^2}}\!\!\!\!&{\scriptscriptstyle \cdots}\!\!\!\!&
{\scriptscriptstyle \frac{{\qbinom {t_1-1} 0}}{\coef{x_m,y_1}^{t_1}}}&
\ &
{\scriptscriptstyle \frac{{\qbinom 0 0}}{\coef{x_m,y_n}^1}}\!\!\!\!&
{\scriptscriptstyle \frac{{\qbinom 1 0}}{\coef{x_m,y_n}^2}}\!\!\!\!&{\scriptscriptstyle \cdots}\!\!\!\!&
{\scriptscriptstyle \frac{{\qbinom {t_n-1} 0}}{\coef{x_m,y_n}^{t_n}}}
\\
{\scriptscriptstyle \frac{{\qbinom 1 1}}{\coef{x_m,y_1}^2}}\!\!\!\!&
{\scriptscriptstyle \frac{{\qbinom 2 1}}{\coef{x_m,y_1}^3}}\!\!\!\!&{\scriptscriptstyle \cdots}\!\!\!\!&
{\scriptscriptstyle \frac{{\qbinom {t_1} 1}}{\coef{x_m,y_1}^{t_1+1}}}&
\ &
{\scriptscriptstyle \frac{{\qbinom 1 1}}{\coef{x_m,y_n}^2}}\!\!\!\!&
{\scriptscriptstyle \frac{{\qbinom 2 1}}{\coef{x_m,y_n}^3}}\!\!\!\!&{\scriptscriptstyle \cdots}\!\!\!\!&
{\scriptscriptstyle \frac{{\qbinom {t_n} 1}}{\coef{x_m,y_n}^{t_n+1}}}
\\
{\scriptstyle\cdot}\!\!\!\!&
{\scriptstyle\cdot}\!\!\!\!&\ \!\!\!\!&
{\scriptstyle\cdot}&
\!\!\!\!\cdots\!\!\! &
{\scriptstyle\cdot}\!\!\!\!&
{\scriptstyle\cdot}\!\!\!\!&\ \!\!\!\!&
{\scriptstyle\cdot}
\\
{\scriptstyle\cdot}\!\!\!\!&
{\scriptstyle\cdot}\!\!\!\!&\ \!\!\!\!&
{\scriptstyle\cdot}&
\ &
{\scriptstyle\cdot}\!\!\!\!&
{\scriptstyle\cdot}\!\!\!\!&\ \!\!\!\!&
{\scriptstyle\cdot}
\\
{\scriptstyle\cdot}\!\!\!\!&
{\scriptstyle\cdot}\!\!\!\!&\ \!\!\!\!&
{\scriptstyle\cdot}&
\ &
{\scriptstyle\cdot}\!\!\!\!&
{\scriptstyle\cdot}\!\!\!\!&\ \!\!\!\!&
{\scriptstyle\cdot}
\\
{\scriptscriptstyle \frac{{\qbinom {s_m-1} {s_m-1}}}{\coef{x_m,y_1}^{s_m}}}\!\!\!\!&
{\scriptscriptstyle \frac{{\qbinom {s_m} {s_m-1}}}{\coef{x_m,y_1}^{s_m+1}}}\!\!\!\!&{\scriptscriptstyle \cdots}\!\!\!\!&
{\scriptscriptstyle \frac{{\qbinom {s_m+t_1-2} {s_m-1}}}{\coef{x_m,y_1}^{s_m+t_1-1}}}&
\ &
{\scriptscriptstyle \frac{{\qbinom {s_m-1} {s_m-1}}}{\coef{x_m,y_n}^{s_m}}}\!\!\!\!&
{\scriptscriptstyle \frac{{\qbinom {s_m} {s_m-1}}}{\coef{x_m,y_n}^{s_m+1}}}\!\!\!\!&{\scriptscriptstyle \cdots}\!\!\!\!&
{\scriptscriptstyle \frac{{\qbinom {s_m+t_n-2} {s_m-1}}}{\coef{x_m,y_n}^{s_m+t_n-1}}}
\end{matrix}
\!\right),
\end{multline}
where $\coef{x,y}^e:=(y-x)(qy-x)(q^2y-x)\cdots(q^{e-1}y-x)$.
Then we have

\vbox{
\begin{multline}
\det A =
q^{\frac {1} {6}
\sum _{i=1} ^{m}(s_i-1)s_i(2s_i-1)}
\(\prod _{1\le i<j\le m}
\prod _{g=1} ^{s_i} 
\prod _{h=1} ^{s_j} (q^{g-1}x_j-q^{h-1}x_i)\)\\
\times
\(\prod _{1\le i<j\le n}
\prod _{g=1} ^{t_i}
\prod _{h=1} ^{t_j}
      (q^{g-1}y_i-q^{h-1}y_j)\)
\(\prod _{j=1} ^{m}
\prod _{i=1} ^{n}
\prod _{g=1} ^{t_i}
\prod _{h=1} ^{s_j}
    \frac  1{(q^{g+h-2}y_i-x_j)}\).
\end{multline}
\quad \quad \qed}
\end{Theorem}

\subsection{More determinants containing derivatives and compositions
of series}
Inspired by formulae of Mina 
\machSeite{MinaAA}\cite{MinaAA}, Kedlaya
\machSeite{KedlAA}\cite{KedlAA} and Strehl and Wilf 
\machSeite{StWiAA}\cite{StWiAA} for determinants of
matrices the entries of which being given by (coefficients of)
{\it multiple derivatives} and {\it compositions of formal power
series} (see also 
\machSeite{KratBN}\cite[Lemma~16]{KratBN}), 
Chu embedded all these in a larger
context in the remarkable systematic study
\machSeite{ChuWBG}\cite{ChuWBG}.
He shows that, at the heart of these formulae, there is the
{\it Fa\`a di Bruno formula}\footnote{As one can read in
\machSeite{JohWAE}\cite{JohWAE}, ``Fa\`a di Bruno was neither the
first to state the formula that bears his name nor the first to prove
it." In Section~4 of that article, the author tries to trace back the
roots of the formula. It is apparently impossible to find the author
of the formula with certainty. In his book
\machSeite{ArboAA}\cite[p.~312]{ArboAA}, Arbogast describes a recursive rule
how, from the top term, to generate all other terms in the
formula. However, the explicit formula is never written down. 
(I am not able to verify the conclusions in 
\machSeite{CraiAA}\cite{CraiAA}. It seems to
me that the author mixes the knowledge that we have today with what is
really written in 
\machSeite{ArboAA}\cite{ArboAA}.)
The formula appears explicitly in Lacroix's book
\machSeite{LacrAA}\cite[p.~629]{LacrAA}, but Lacroix's precise
sources remain unknown. I refer the reader to
\machSeite{JohWAE}\cite[Sec.~4]{JohWAE} and
\machSeite{CraiAA}\cite{CraiAA}
for more detailed remarks on the history of the formula.} 
for multiple derivatives of a composition of
two formal power series. Using it, he derives the following
determinant reduction formulae 
\machSeite{ChuWBG}\cite[Theorems~4.1 and 4.2]{ChuWBG} for determinants
of matrices containing multiple derivatives of compositions of formal
power series. 

\begin{Theorem} \label{thm:Chu1}
Let $f(x)$ and $\phi_k(x)$ and $w_k(x)$, $k=0,1,\dots,n$, 
be formal power series in $x$
with coefficients in a commutative ring. Then
\begin{equation} \label{eq:Chu1}
\det_{0\le i,j,\le n}\(\frac {d^j}
{dx^j}\bigg(w_j(x)\phi_i\big(f(x)\big)\bigg)\)=
\big(f'(x)\big)^{\binom {n+1}2}
\(\prod _{k=0} ^{n}w_k(x)\)\det_{0\le i,j,\le
n}\(\phi_i^{(j)}\big(f(x)\big)\),
\end{equation}
where $\phi^{(j)}(x)$ is short for $\frac {d^j} {dx^j}\phi(x)$. If, in
addition, $w_k(x)$ is a polynomial of degree at most $k$,
$k=1,2,\dots,n$, then
\begin{equation} \label{eq:Chu2}
\det_{1\le i,j,\le n}\(\frac {d^j}
{dx^j}\bigg(w_j(x)\phi_i\big(f(x)\big)\bigg)\)=
\big(f'(x)\big)^{\binom {n+1}2}
\(\prod _{k=1} ^{n}w_k(x)\)\det_{1\le i,j,\le
n}\(\phi_i^{(j)}\big(f(x)\big)\).
\end{equation}
\quad \quad \qed
\end{Theorem}

Specialising the series $\phi_i(x)$ so that the determinants
on the right-hand sides of \eqref{eq:Chu1} or \eqref{eq:Chu2} can be
evaluated, he obtains numerous nice corollaries.
Possible choices are $\phi_i(x)=\exp(y_ix)$, 
$\phi_i(x)=\log(1+y_ix)$, 
$\phi_i(x)=x^{y_i}$, or 
$\phi_i(x)=(a_i+b_ix)/(c_i+d_ix)$. See
\machSeite{ChuWBG}\cite[Cor.~4.3 and 4.4]{ChuWBG} for the
corresponding results. 

Further reduction formulae and determinant evaluations from 
\machSeite{ChuWBG}\cite{ChuWBG} 
address determinants of matrices
formed out of coefficients of iterated compositions of formal power
series. In order to have a convenient notation, let us write
$f^{\coef{n}}(x)$ for the $n$-fold composition of $f$ with itself,
$$f^{\coef{n}}(x)=f(f(\cdots(f(x)))),$$
with $n$ occurrences of $f$ on the right-hand side. Chu shows 
\machSeite{ChuWBG}\cite[Sec.~1.4]{ChuWBG} that it
is possible to extend this $n$-fold composition to values of $n$ other
than non-negative integers. This given, Theorems~4.6 and 4.7 from
\machSeite{ChuWBG}\cite{ChuWBG} read as follows.

\begin{Theorem} \label{thm:Chu3}
Let $f(x)=x+
\sum _{m=2} ^{\infty}f_mx^m$, $g(x)=
\sum _{m=1} ^{\infty}g_mx^m$ and $w_k(x)$, $k=1,2,\dots,n$, 
be formal power series with
coefficients in some commutative ring. Then
\begin{equation} \label{eq:Chu3}
\det_{1\le i,j\le n}\Big([x^j]w_j(x)f^{\coef{y_i}}\big(g(x)\big)\Big)=
f_2^{\binom n2}g_1^{\binom {n+1}2} \(
\prod _{k=1} ^{n}w_k(0)\)\(
\prod _{1\le i<j\le n} ^{}(y_j-y_i)\),
\end{equation}
where $[x^j]h(x)$ denotes the coefficient of $x^j$ in the series $h(x)$.
If, in addition, $w_n(0)=0$, then
\begin{multline} \label{eq:Chu4}
\det_{1\le i,j\le
n}\Big([x^{j+1}]w_j(x)f^{\coef{y_i}}\big(g(x)\big)\Big)
\\=
f_2^{\binom n2}g_1^{\binom {n+1}2} 
\(\prod _{1\le i<j\le n} ^{}(y_j-y_i)\)
\det_{1\le i,j\le n}\Big([x^{1+j-i}]w_j(x)\Big).
\end{multline}
Furthermore, we have
\begin{equation} \label{eq:Chu5}
\det_{1\le i,j\le n}\Big([x^{j+1}]f^{\coef{y_i}}(x)\Big)=
f_2^{\binom {n+1}2}\(\prod _{k=1} ^{n}y_k\)\(
\prod _{1\le i<j\le n} ^{}(y_j-y_i)\).
\end{equation}
\quad \quad \qed
\end{Theorem}

\subsection{More on Hankel determinant evaluations}
Section~2.7 of 
\machSeite{KratBN}\cite{KratBN} was devoted to {\it Hankel
determinants}. There, I tried to convince the reader that, whenever
you think that a certain Hankel determinant evaluates nicely, then 
the explanation will be (sometimes more sometimes less) hidden 
in the theory of {\em continued fractions} and 
{\em orthogonal polynomials}. In retrospect,
it seems that the success of this try 
was mixed. Since readers are always right, this has to be blamed
entirely on myself, and, indeed, the purpose of the present subsection
is to rectify some shortcomings from then.

Roughly speaking, I explained in 
\machSeite{KratBN}\cite{KratBN} that, given a Hankel determinant
\begin{equation} \label{eq:Hankel1} 
\det_{0\le i,j\le n-1}(\mu_{i+j}),
\end{equation}
to find its evaluation one should expand the generating function of
the sequence of coefficients $(\mu_k)_{k\ge0}$ in terms of a continued
fraction, respectively find the sequence of orthogonal polynomials
$(p_n(x))_{n\ge0}$ with moments $\mu_k$, $k=0,1,\dots$, and then the
value of the Hankel determinant \eqref{eq:Hankel1} can be read off the
coefficients of the continued fraction, respectively from the
recursion coefficients of the orthogonal polynomials. What I missed to
state is that the knowledge of the
orthogonal polynomials makes it also possible to find the value of the
Hankel determinants which start with $\mu_1$ and $\mu_2$, respectively
(instead of $\mu_0$). In the theorem below I summarise the results
that were already discussed in 
\machSeite{KratBN}\cite{KratBN}
(for which classical references are 
\machSeite{WallCF}%
\cite[Theorem 51.1]{WallCF} or
\machSeite{VienAE}%
\cite[Cor.~6, (19), on p.~IV-17; Proposition~1, (7), 
on p.~V-5]{VienAE}), and I add the two missing ones.

\begin{Theorem}
\label{cor:cfracHankel}
Let $(\mu_k)_{k\ge0}$ be a sequence with generating
function $\sum_{k=0}^\infty{\mu_k}x^k$ written in the form
\begin{equation}
\label{eq:momentgf}
\sum_{k=0}^\infty{\mu_k}x^k=\cfrac{
	\mu_0}
		{1+a_0x-\cfrac{
			b_1x^2}
				{1+a_1x-\cfrac{
					b_2x^2}
						{1+a_2x-\cdots}}}\quad .
\end{equation}
Then 
\begin{equation} \label{eq:Hankel2} 
\det_{0\le i,j\le n-1}(\mu_{i+j})=\mu_0^nb_1^{n-1}b_2^{n-2}\cdots
b_{n-2}^2b_{n-1}. 
\end{equation}

Let $(p_n(x))_{n\ge0}$ be a sequence of monic polynomials, the
polynomial $p_n(x)$ having degree $n$, which is orthogonal with
respect to some functional $L$, that is, $L(p_m(x)p_n(x))=\de_{m,n}c_{n}$,
where the $c_n$'s are some non-zero constants and $\de_{m,n}$ is the
Kronecker delta. Let
\begin{equation}
p_{n+1}(x)=(a_{n}+x)p_{n}(x)-b_{n}p_{n-1}(x)
\label{eq:three-term2}
\end{equation}
be the corresponding three-term recurrence which is guaranteed by
Favard's theorem. Then the generating function $\sum _{k=0}
^{\infty}\mu_kx^k$ for the moments
$\mu_k=L(x^k)$ satisfies \eqref{eq:momentgf} with the $a_i$'s and
$b_i$'s being the coefficients in the three-term recurrence
\eqref{eq:three-term2}. In particular, the Hankel determinant
evaluation \eqref{eq:Hankel2} holds, with the $b_i$'s from the
three-term recurrence \eqref{eq:three-term2}.

If $(q_n)_{n\ge0}$ is the sequence recursively defined by $q_0=1$,
$q_1=-a_0$, and 
$$q_{n+1}=-a_n q_n-b_nq_{n-1},$$
then in the situation above we also have
\begin{equation} \label{eq:Hankel3} 
\det_{0\le i,j\le n-1}(\mu_{i+j+1})=\mu_0^nb_1^{n-1}b_2^{n-2}\cdots
b_{n-2}^2b_{n-1}q_n
\end{equation}
and
\begin{equation} \label{eq:Hankel4} 
\det_{0\le i,j\le n-1}(\mu_{i+j+2})=\mu_0^nb_1^{n-1}b_2^{n-2}\cdots
b_{n-2}^2b_{n-1}
\sum _{k=0} ^{n}q_k^2b_{k+1}\cdots b_{n-1}b_n. 
\end{equation}
\quad \quad \qed
\end{Theorem}

I did not find a reference for \eqref{eq:Hankel3} and
\eqref{eq:Hankel4}. These two identities follow however easily from
Viennot's combinatorial model
\machSeite{VienAE}\cite{VienAE} for orthogonal polynomials and Hankel
determinants of moments. More precisely, in this theory the moments
$\mu_k$ are certain generating functions for {\it Motzkin paths}, and,
due to Theorem~\ref{thm:nonint}, the
Hankel determinants $\det_{0\le i,j\le n-1}(\mu_{i+j+m})$ are
generating functions for families $(P_1,P_2,\dots,P_n)$
of non-intersecting Motzkin paths, $P_i$ running from $(-i,0)$ to
$(j+m,0)$. In the case $m=0$, it is explained in
\machSeite{VienAE}\cite[Ch.~IV]{VienAE} how to find the corresponding
Hankel determinant evaluation \eqref{eq:Hankel2}
using this combinatorial model. The idea is that in that case there is
a unique family of non-intersecting Motzkin paths, and its weight
gives the right-hand side of \eqref{eq:Hankel2}.
If $m=1$ or $m=2$ one can argue similarly. The paths are uniquely
determined with the exception of their portions in 
the strip $0\le x\le m$. The various possibilities that one has there
then yield the right-hand sides of \eqref{eq:Hankel3} and
\eqref{eq:Hankel4}. 

Since there are so many explicit families of orthogonal polynomials,
and, hence, so many ways to apply the above theorem, I listed only a
few standard Hankel determinant evaluations explicitly in
\machSeite{KratBN}\cite{KratBN}. I did append a long list of
references and sketched in which ways these give rise to more Hankel
determinant evaluations. Apparently, these remarks were at times too
cryptic, in particular concerning the theme {\it ``orthogonal polynomials as
moments}.'' This is treated systematically in the two papers
\machSeite{IsStAB}\machSeite{IsStAC}\cite{IsStAB,IsStAC} by Ismail and
Stanton. There it is shown that certain classical polynomials
$(r_n(x))_{n\ge0}$, such as, for example,
the {\em Laguerre polynomials}, the {\em Meixner polynomials}, or
the {\em Al-Salam--Chihara polynomials} (but there are others as well,
see
\machSeite{IsStAB}\machSeite{IsStAC}\cite{IsStAB,IsStAC}), are {\it
moments} of other families of classical orthogonal polynomials. Thus,
application of Theorem~\ref{cor:cfracHankel} with $\mu_n=r_n(x)$ 
immediately tells that
the evaluations of the corresponding Hankel determinants
\begin{equation} \label{eq:pn} 
\det_{0\le i,j\le n-1}\big(r_{i+j}(x)\big)
\end{equation}
(and also the higher ones in \eqref{eq:Hankel3} and
\eqref{eq:Hankel4}) are known. In particular, the explicit forms can
be extracted from the coefficients of the three-term recursions for
these other families of orthogonal polynomials. Thus, whenever you encounter a
determinant of the form \eqref{eq:pn}, you must check whether
$(r_n(x))_{n\ge0}$ is a family of orthogonal polynomials (which, as I
explained in
\machSeite{KratBN}\cite{KratBN}, one does by consulting the
compendium \machSeite{KoSwAA}\cite{KoSwAA} of hypergeometric
orthogonal polynomials compiled by Koekoek and Swarttouw), and if the
answer is ``yes", you will find the solution of your determinant
evaluation through the results in
\machSeite{IsStAB}\machSeite{IsStAC}\cite{IsStAB,IsStAC}
by applying Theorem~\ref{cor:cfracHankel}.

\medskip
While Theorem~\ref{cor:cfracHankel} describes in detail the connexion
between Hankel determinants and the continued fractions of the type
\eqref{eq:momentgf}, which are
often called {\it $J$-fractions} (which is short for {\it Jacobi continued
fractions}), I missed to tell in
\machSeite{KratBN}\cite{KratBN} that there is also a close relation
between Hankel determinants and so-called {\it $S$-fractions}
(which is short for {\it Stieltjes continued fractions}). I try to remedy this
by the theorem below (cf.\ for example
\machSeite{JoThAA}\cite[Theorem~7.2]{JoThAA}, where  
$S$-fractions are called {\it regular $C$-fractions}). 
In principle, since $S$-fractions are special cases of $J$-fractions
\eqref{eq:momentgf} in which the coefficients $a_i$ are all zero,
the corresponding result for the Hankel determinants is in fact implied by
Theorem~\ref{cor:cfracHankel}. Nevertheless, it is useful to state it
separately. I am not
able to give a reference for \eqref{eq:Hankel7}, but, again, it is not too
difficult to derive it from Viennot's combinatorial model
\machSeite{VienAE}\cite{VienAE} for orthogonal
polynomials and moments that was mentioned above. 

\begin{Theorem} \label{thm:cfrac2}
Let $(\mu_k)_{k\ge0}$ be a sequence with generating
function $\sum_{k=0}^\infty{\mu_k}x^k$ written in the form
\begin{equation}
\label{eq:momentgfS}
\sum_{k=0}^\infty{\mu_k}x^k=\cfrac{
	\mu_0}
		{1+\cfrac{
			a_1x}
				{1+\cfrac{
					a_2x}
						{1+\cdots}}}\quad .
\end{equation}
Then 
\begin{align} \label{eq:Hankel5} 
\det_{0\le i,j\le n-1}(\mu_{i+j})&=\mu_0^n(a_1a_2)^{n-1}(a_3a_4)^{n-2}\cdots
(a_{2n-5}a_{2n-4})^2(a_{2n-3}a_{2n-2}),\\
\label{eq:Hankel6} 
\det_{0\le i,j\le n-1}(\mu_{i+j+1})&=(-1)^n\mu_0^n
a_1^n(a_2a_3)^{n-1}(a_4a_5)^{n-2}\cdots
(a_{2n-4}a_{2n-3})^2(a_{2n-2}a_{2n-1}),
\end{align}
and

\vbox{
\begin{multline} \label{eq:Hankel7} 
\det_{0\le i,j\le n-1}(\mu_{i+j+2})=\mu_0^n
a_1^n(a_2a_3)^{n-1}(a_4a_5)^{n-2}\cdots
(a_{2n-4}a_{2n-3})^2(a_{2n-2}a_{2n-1})\\
\times
\sum _{0\le i_1-1<i_2-2<\dots<i_n-n\le n} ^{}a_{i_1}a_{i_2}\cdots a_{i_n}. 
\end{multline}
\quad \quad \qed}
\end{Theorem}

Using this theorem, Tamm 
\machSeite{TammAA}\cite[Theorem~3.1]{TammAA}
observed that from {\it Gau\3' continued fraction for the ratio of two
contiguous $_2F_1$-series} one can deduce several interesting
binomial Hankel determinant evaluations, some of them had already
been found earlier by E\u gecio\u glu, Redmond and Ryavec 
\machSeite{EgRiAA}\cite[Theorem~4]{EgRiAA}
while working on {\it polynomial Riemann hypotheses}.
Gessel and Xin
\machSeite{GeXiAB}\cite{GeXiAB} undertook a systematic analysis of
this approach, and they arrived at a set of 18 Hankel determinant
evaluations, which I list as
\eqref{eq:3nA}--\eqref{eq:3nR} in the theorem below. They are preceded by the
Hankel determinant evaluation \eqref{eq:3nAAAA}, which appears only in 
\machSeite{EgRiAA}\cite[Theorem~4]{EgRiAA}.



\begin{Theorem} \label{thm:Tamm}
For any positive integer $n$, there hold 
\begin{equation}
\label{eq:3nAAAA}
\det_{0\le i,j\le n-1}\(\binom {3i+3j+2}{i+j}\) 
=
\prod _{i=0} ^{n-1}\frac {(\frac {4} {3})_i\,(\frac {5} {6})_i\,
(\frac {5} {3})_i\,(\frac {7} {6})_i} 
{(\frac {3} {2})_{2i}\,(\frac {5} {2})_{2i}}\(\frac {27} {4}\)^{2i},
\end{equation}
\begin{equation} \label{eq:3nA}
\det_{0\le i,j\le n-1}\(\frac {1} {3i+3j+1}\binom {3i+3j+1}{i+j}\) 
=
\prod _{i=0} ^{n-1}\frac {(\frac {2} {3})_i\,(\frac {1} {6})_i\,
(\frac {4} {3})_i\,(\frac {5} {6})_i} 
{(\frac {1} {2})_{2i}\,(\frac {3} {2})_{2i}}\(\frac {27} {4}\)^{2i},
\end{equation}
\begin{equation}
\label{eq:3nB}
\det_{0\le i,j\le n-1}\(\frac {1} {3i+3j+4}\binom {3i+3j+4}{i+j+1}\) 
=
\prod _{i=0} ^{n-1}\frac {(\frac {4} {3})_i\,(\frac {5} {6})_i\,
(\frac {5} {3})_i\,(\frac {7} {6})_i} 
{(\frac {3} {2})_{2i}\,(\frac {5} {2})_{2i}}\(\frac {27} {4}\)^{2i},
\end{equation}
\begin{equation}
\label{eq:3nC}
\det_{0\le i,j\le n-1}\(\frac {1} {3i+3j+2}\binom {3i+3j+2}{i+j+1}\) 
=
\prod _{i=0} ^{n-1}\frac {(\frac {4} {3})_i\,(\frac {5} {6})_i\,
(\frac {5} {3})_i\,(\frac {7} {6})_i} 
{(\frac {3} {2})_{2i}\,(\frac {5} {2})_{2i}}\(\frac {27} {4}\)^{2i},
\end{equation}
\begin{equation}
\label{eq:3nD}
\det_{0\le i,j\le n-1}\(\frac {1} {3i+3j+5}\binom {3i+3j+5}{i+j+2}\) 
=
\prod _{i=0} ^{n}\frac {(\frac {2} {3})_i\,(\frac {1} {6})_i\,
(\frac {4} {3})_i\,(\frac {5} {6})_i} 
{(\frac {1} {2})_{2i}\,(\frac {3} {2})_{2i}}\(\frac {27} {4}\)^{2i},
\end{equation}
\begin{equation}
\label{eq:3nE}
\det_{0\le i,j\le n-1}\(\frac {2} {3i+3j+1}\binom {3i+3j+1}{i+j+1}\) 
=
\prod _{i=0} ^{n}\frac {(\frac {2} {3})_i\,(\frac {1} {6})_i\,
(\frac {4} {3})_i\,(\frac {5} {6})_i} 
{(\frac {1} {2})_{2i}\,(\frac {3} {2})_{2i}}\(\frac {27} {4}\)^{2i},
\end{equation}
\begin{equation}
\label{eq:3nF}
\det_{0\le i,j\le n-1}\(\frac {2} {3i+3j+4}\binom {3i+3j+4}{i+j+2}\) 
=
\prod _{i=0} ^{n}\frac {(\frac {4} {3})_i\,(\frac {5} {6})_i\,
(\frac {5} {3})_i\,(\frac {7} {6})_i} 
{(\frac {3} {2})_{2i}\,(\frac {5} {2})_{2i}}\(\frac {27} {4}\)^{2i},
\end{equation}
\begin{equation}
\label{eq:3nG}
\det_{0\le i,j\le n-1}\(\frac {2} {(3i+3j+1)(3i+3j+2)}
       \binom {3i+3j+2}{i+j+1}\) 
=
\prod _{i=0} ^{n-1}2\frac {(\frac {5} {3})_i\,(\frac {1} {6})_i\,
(\frac {7} {3})_i\,(\frac {5} {6})_i} 
{(\frac {3} {2})_{2i}\,(\frac {5} {2})_{2i}}\(\frac {27} {4}\)^{2i},
\end{equation}
\begin{multline}
\label{eq:3nH}
\det_{0\le i,j\le n-1}\(\frac {2} {(3i+3j+4)(3i+3j+5)}
      \binom {3i+3j+5}{i+j+2}\) \\
=(-1)^n
\prod _{i=1} ^{n}\frac {(\frac {5} {3})_i\,(\frac {1} {6})_i\,
(\frac {4} {3})_i\,(-\frac {1} {6})_i} 
{(\frac {1} {2})_{2i}\,(\frac {3} {2})_{2i}}\(\frac {27} {4}\)^{2i},
\end{multline}
\begin{equation}
\label{eq:3nI}
\det_{0\le i,j\le n-1}\(\frac {(9i+9j+5)} {(3i+3j+1)(3i+3j+2)}
       \binom {3i+3j+2}{i+j+1}\) 
=
\prod _{i=0} ^{n-1}5\frac {(\frac {2} {3})_i\,(\frac {7} {6})_i\,
(\frac {4} {3})_i\,(\frac {11} {6})_i} 
{(\frac {3} {2})_{2i}\,(\frac {5} {2})_{2i}}\(\frac {27} {4}\)^{2i},
\end{equation}
\begin{equation}
\label{eq:3nJ}
\det_{0\le i,j\le n-1}\(\frac {(9i+9j+14)} {(3i+3j+4)(3i+3j+5)}
      \binom {3i+3j+5}{i+j+2}\) =
\prod _{i=1} ^{n}2\frac {(\frac {2} {3})_i\,(\frac {7} {6})_i\,
(\frac {1} {3})_i\,(\frac {5} {6})_i} 
{(\frac {1} {2})_{2i}\,(\frac {3} {2})_{2i}}\(\frac {27} {4}\)^{2i}.
\end{equation}
Let $a_0=-2$ and $a_m=\frac {1} {3m+1}\binom {3m+1}m$ for $m\ge1$. Then
\begin{equation} \label{eq:3nK} 
\det_{0\le i,j\le n-1}(a_{i+j})=
\prod _{i=0} ^{n-1}(-2)\frac {(\frac {1} {3})_i\,(-\frac {1} {6})_i\,
(\frac {5} {3})_i\,(\frac {7} {6})_i} 
{(\frac {1} {2})_{2i}\,(\frac {3} {2})_{2i}}\(\frac {27} {4}\)^{2i}.
\end{equation}
Let $b_0=10$ and $b_m=\frac {2} {3m+2}\binom {3m+2}m$ for $m\ge1$. Then
\begin{equation} \label{eq:3nL} 
\det_{0\le i,j\le n-1}(b_{i+j})=
\prod _{i=0} ^{n-1}10\frac {(\frac {2} {3})_i\,(\frac {1} {6})_i\,
(\frac {7} {3})_i\,(\frac {11} {6})_i} 
{(\frac {3} {2})_{2i}\,(\frac {5} {2})_{2i}}\(\frac {27} {4}\)^{2i}.
\end{equation}
Furthermore,
\begin{equation} \label{eq:3nM}
\det_{0\le i,j\le n-1}\(\frac {2} {3i+3j+5}\binom {3i+3j+5}{i+j+1}\) 
=
\prod _{i=0} ^{n}\frac {(\frac {4} {3})_i\,(\frac {5} {6})_i\,
(\frac {5} {3})_i\,(\frac {7} {6})_i} 
{(\frac {3} {2})_{2i}\,(\frac {5} {2})_{2i}}\(\frac {27} {4}\)^{2i}.
\end{equation}
Let $c_0=\frac {7} {2}$ and 
$c_m=\frac {2} {3m+1}\binom {3m+1}{m+1}$ for $m\ge1$. Then
\begin{equation} \label{eq:3nN} 
\det_{0\le i,j\le n-1}(c_{i+j})=
\prod _{i=0} ^{n}\frac {(\frac {4} {3})_i\,(\frac {5} {6})_i\,
(\frac {5} {3})_i\,(\frac {7} {6})_i} 
{(\frac {3} {2})_{2i}\,(\frac {5} {2})_{2i}}\(\frac {27} {4}\)^{2i}.
\end{equation}
Let $d_0=-5$ and $d_m=\frac {8} {(3m+1)(3m+2)}\binom {3m+2}m$ for $m\ge1$. Then
\begin{equation} \label{eq:3nO} 
\det_{0\le i,j\le n-1}(d_{i+j})=
\prod _{i=0} ^{n-1}(-5)\frac {(\frac {4} {3})_i\,(-\frac {1} {6})_i\,
(\frac {8} {3})_i\,(\frac {7} {6})_i} 
{(\frac {3} {2})_{2i}\,(\frac {5} {2})_{2i}}\(\frac {27} {4}\)^{2i}.
\end{equation}
Furthermore,
\begin{equation} \label{eq:3nP}
\det_{0\le i,j\le n-1}\(\frac {8} {(3i+3j+4)(3i+3j+5)}
\binom {3i+3j+5}{i+j+1}\) 
=
\prod _{i=0} ^{n-1}2\frac {(\frac {7} {3})_i\,(\frac {5} {6})_i\,
(\frac {8} {3})_i\,(\frac {7} {6})_i} 
{(\frac {5} {2})_{2i}\,(\frac {7} {2})_{2i}}\(\frac {27} {4}\)^{2i}.
\end{equation}
Let $e_0=14$ and 
$e_m=\frac {2(9m+5)} {(3m+1)(3m+2)}\binom {3m+2}m$ for $m\ge1$. Then
\begin{equation} \label{eq:3nQ} 
\det_{0\le i,j\le n-1}(e_{i+j})=
\prod _{i=0} ^{n-1}14\frac {(\frac {1} {3})_i\,(\frac {5} {6})_i\,
(\frac {5} {3})_i\,(\frac {13} {6})_i} 
{(\frac {3} {2})_{2i}\,(\frac {5} {2})_{2i}}\(\frac {27} {4}\)^{2i}.
\end{equation}
Furthermore,
\begin{equation} \label{eq:3nR}
\det_{0\le i,j\le n-1}\(\frac {2(9i+9j+14)} {(3i+3j+4)(3i+3j+5)}
\binom {3i+3j+5}{i+j+1}\) 
=
\prod _{i=1} ^{n}2\frac {(\frac {2} {3})_i\,(\frac {7} {6})_i\,
(\frac {1} {3})_i\,(\frac {5} {6})_i} 
{(\frac {1} {2})_{2i}\,(\frac {3} {2})_{2i}}\(\frac {27} {4}\)^{2i}.
\end{equation}
\quad \quad \qed
\end{Theorem}

Some of the numbers appearing on the right-hand sides of the formulae
in this theorem have combinatorial significance, although no intrinsic
explanation is known why this is the case. More precisely, the numbers on the
right-hand sides of \eqref{eq:3nA}, 
\eqref{eq:3nD} and \eqref{eq:3nE} count {\it cyclically symmetric 
transpose-complementary plane partitions} (cf.\ 
\machSeite{MiRRAD}\cite{MiRRAD} and
\machSeite{BresAO}\cite{BresAO}), whereas those on the right-hand
sides of \eqref{eq:3nAAAA}, 
\eqref{eq:3nB}, \eqref{eq:3nC} and \eqref{eq:3nF} count {\em vertically
symmetric alternating sign matrices} (cf.\
\machSeite{KupeAH}\cite{KupeAH} and
\machSeite{BresAO}\cite{BresAO}).

In \machSeite{EgRiAA}\cite[Theorem~4]{EgRiAA},
E\u gecio\u glu, Redmond and Ryavec prove also the following
common generalisation of \eqref{eq:3nB} and \eqref{eq:3nC}.
(The first identity is the special case $x=0$, while the second is the
special case $x=1$ of the following theorem.)

\begin{Theorem} \label{thm:EgRR}
For $m\ge0$, let $s_m(x)=
\sum _{k=0} ^{m}\frac {k+1} {m+1}\binom {3m-k+1}{m-k}x^k$. Then, for
any positive integer $n$, there holds
\begin{equation} \label{eq:EgRR}
\det_{0\le i,j\le n-1}(s_{i+j}(x))= 
\prod _{i=0} ^{n-1}\frac {(\frac {4} {3})_i\,(\frac {5} {6})_i\,
(\frac {5} {3})_i\,(\frac {7} {6})_i} 
{(\frac {3} {2})_{2i}\,(\frac {5} {2})_{2i}}\(\frac {27} {4}\)^{2i}.
\end{equation}  
\quad \quad \qed
\end{Theorem}

As I mentioned above,
in \machSeite{KratBN}\cite{KratBN} I only stated a few special Hankel
determinant evaluations explicitly, because there are too many ways to
apply Theorems~\ref{cor:cfracHankel} and \ref{thm:cfrac2}. 
I realise, however, that I should have stated the
evaluation of the {\it Hankel determinant of Catalan numbers} there.
I make this up now by doing this in the theorem below. I did not do it
then because orthogonal polynomials are not needed for its evaluation
(the orthogonal polynomials which are tied to Catalan numbers as
moments are {\it Chebyshev polynomials}, but, via
Theorems~\ref{cor:cfracHankel} and \ref{thm:cfrac2}, one would only
cover the cases $m=0,1,2$ in the theorem below). In fact, the Catalan number
$C_n=\frac {1} {n+1}\binom {2n}n$ can be alternatively written as
$C_n=(-1)^n2^{2n+1}\binom {1/2}{n+1}$, and therefore the Hankel
determinant evaluation below follows from 
\machSeite{KratBN}\cite[Theorem~26, (3.12)]{KratBN}.
This latter observation shows that 
even a more general determinant, namely
$\det_{0\le i,j\le n-1}(C_{\lambda_i+j})$, can be evaluated in closed form.
For historical remarks on this ubiquitous determinant see
\machSeite{GhKrAA}\cite[paragraph before the Appendix]{GhKrAA}.

\begin{Theorem} \label{thm:Catalan}
\begin{equation} \label{eq:Catalan} 
\det_{0\le i,j\le n-1}(C_{m+i+j})=
\prod _{1\le i\le j\le m-1} ^{}\frac {2n+i+j} {i+j}.
\end{equation}
\quad \quad \qed
\end{Theorem}

As in 
\machSeite{KratBN}\cite{KratBN}, 
let me conclude the part on Hankel determinants 
by pointing the reader to further
papers containing interesting results on them, high-lighting
sometimes the point of view of orthogonal polynomials that I
explained above, sometimes a combinatorial point of view.
The first point of view is put forward in 
\machSeite{WimpAB}\cite{WimpAB} (see 
\machSeite{KratZZ}\cite{KratZZ} for the solution of the conjectures in
that paper) in order to present Hankel determinant evaluations of
matrices with {\it hypergeometric $_2F_1$-series} as entries. 
The orthogonal polynomials approach is also used in 
\machSeite{CvRIAA}\cite{CvRIAA} to show that a certain Hankel determinant
defined by {\it Catalan numbers} evaluates to {\it Fibonacci numbers}.
In \machSeite{AnWiAA}\cite{AnWiAA}, one finds Hankel determinant
evaluations involving generalisations of the {\it Bernoulli numbers}.
The combinatorial point of view dominates in
\machSeite{AignAA}%
\machSeite{CiglAM}%
\machSeite{CiglAO}%
\machSeite{CiglAV}%
\machSeite{EhreAB}\cite{AignAA,CiglAM,CiglAO,CiglAV,EhreAB}, where Hankel
determinants involving {\it $q$-Catalan numbers, $q$-Stirling numbers}, and
{\it $q$-Fibonacci numbers} are considered.

A very interesting new direction, which seems to have much potential,
is opened up by Luque and Thibon in
\machSeite{LuThAB}\cite{LuThAB}. They show that {\it Selberg-type
integrals} can be evaluated by means of {\it Hankel
hyperdeterminants}, and they prove many hyperdeterminant
generalisations of classical Hankel determinant evaluations.

At last, (but certainly not least!),
I want to draw the reader's attention to Lascoux's
``unorthodox"\footnote{It could easily be that it
is the ``modern" treatment of the theory which must be labelled with
the attribute ``unorthodox." As Lascoux documents in
\machSeite{LascAZ}\cite{LascAZ}, in his treatment {\it he} follows the
tradition of great masters such as Cauchy, Jacobi or Wro\'nski \dots}
approach to Hankel determinants and orthogonal polynomials through
symmetric functions which he presents in detail in
\machSeite{LascAZ}\cite[Ch.~4, 5, 8]{LascAZ}. 
In particular, Theorem~\ref{cor:cfracHankel}, Eq.~\eqref{eq:Hankel2}
are the contents of Theorem~8.3.1 in
\machSeite{LascAZ}\cite{LascAZ} (see also the end of Section~5.3 there), and 
Theorem~\ref{thm:cfrac2}, Eqs.~\eqref{eq:Hankel5} and
\eqref{eq:Hankel6} are the contents of Theorem~4.2.1 in
\machSeite{LascAZ}\cite{LascAZ}.
The usefulness of this symmetric function approach is, 
for example, demonstrated in
\machSeite{HoLMAA}%
\machSeite{HoLMAB}\cite{HoLMAA,HoLMAB} in order to evaluate Hankel
determinants of matrices the entries of which are {\it Rogers--Szeg\H
o}, respectively {\it Meixner polynomials}.

\subsection{More binomial determinants}
A vast part of Section~3 in 
\machSeite{KratBN}\cite{KratBN} is occupied by binomial
determinants. As I mentioned in Section~\ref{sec:comb} of the present article,
an extremely rich source for binomial determinants is rhombus tiling
enumeration. I want to present here some which did not already appear in
\machSeite{KratBN}\cite{KratBN}. 

To begin with, I want to remind the reader of an old problem posed by
Andrews in 
\machSeite{AndrAO}\cite[p.~105]{AndrAO}. The determinant in this problem is a
variation of a determinant which enumerates {\it cyclically
symmetric plane partitions} and {\it descending plane partitions}, 
which was evaluated by Andrews in
\machSeite{AndrAN}\cite{AndrAN} (see also 
\machSeite{KratBN}\cite[Theorem~32]{KratBN}; the latter determinant
arises from the one in \eqref{eq:desc-var} by replacing $j+1$ by $j$
in the bottom of the binomial coefficient).

\begin{Problem} \label{prob:7}
Evaluate the determinant
\begin{equation} \label{eq:desc-var} 
D_1(n):=\det_{0\le i,j\le n-1}\left(\delta_{ij}+\binom {\mu+i+j}{j+1}\right),
\end{equation}
where $\de_{ij}$ is the Kronecker delta.
In particular, show that
\begin{equation} \label{eq:fmu} 
\frac {D_1(2n)} {D_1(2n-1)}=(-1)^{\binom {n-1}2}
\frac {2^n\,\(\frac\mu2+n\)_{\cl{n/2}}\,\(\frac {\mu}
{2}+2n+\frac {1} {2}\)_{n-1}} {(n)_n\,\(-\frac {\mu} {2}-2n+\frac {3}
{2}\)_{\cl{(n-2)/2}}}.
\end{equation}
\end{Problem}

The determinants $D_1(n)$ are rather intriguing. Here are the first 
few values:
\begin{align*}
D_1(1)&=\mu+1,\\
D_1(2)&=(\mu+1)(\mu+2),\\
D_1(3)&= \frac {1} {12}{( \mu + 1  )  ( \mu + 2  )  ( \mu + 3  )  
     ( \mu + 14  ) },\\
D_1(4)&=\frac {1} {72}
{( \mu + 1  )  ( \mu + 2  )  ( \mu + 3  )  
     ( \mu + 4  )  ( \mu + 9  )  ( \mu + 14  ) },\\
D_1(5)&= \frac {1} {8640}
( \mu + 1  )  ( \mu + 2  )  ( \mu + 3  )  
     ( \mu + 4  )  ( \mu + 5  )  ( \mu + 9 
) \\
&\kern3cm
\times
     ( 3432 + 722 \mu + 45 \mu^2 + \mu^3 ) ,\\
D_1(6)&=\frac {1} {518400}
( \mu + 1  )  ( \mu + 2  )  ( \mu + 3  )  
     ( \mu + 4  )  ( \mu + 5  )  ( \mu + 6  )  
     ( \mu + 8  )  ( \mu + 13  ) \\
&\kern3cm
\times
( \mu + 15  )  
     ( 3432 + 722 \mu + 45 \mu^2 + \mu^3 ) ,\\
D_1(7)&=\frac {1} {870912000}
( \mu + 1  )  ( \mu + 2  )  ( \mu + 3  )  
     ( \mu + 4  )  ( \mu + 5  )  ( \mu + 6  )  
     ( \mu + 7  )  ( \mu + 8  ) ^2\\
&\kern3cm
\times
     ( \mu + 13  )  ( \mu + 15  ) ^2 
     ( \mu + 34  )  (  \mu^3 + 47 \mu^2 + 954 \mu + 5928) 
     ,\\
D_1(8)&=\frac {1} {731566080000}
( \mu + 1  )  ( \mu + 2  )  ( \mu + 3  )  
     ( \mu + 4  )  ( \mu + 5  )  ( \mu + 6  )  
     ( \mu + 7  )  ( \mu + 8  ) ^3 \\
&\kern3cm\times
     ( \mu + 10  )  ( \mu + 15  ) ^2 
     ( \mu + 17  )  ( \mu + 19  ) \\
&\kern3cm\times
 ( \mu + 21  )  
     ( \mu + 34  )  (  \mu^3 + 47 \mu^2 + 954 \mu + 5928 ) 
     .
\end{align*}
``So," these determinants factor almost completely, there is only a
relatively small (in degree) irreducible factor which is not
linear. (For example, this factor is of degree $6$ for $D_1(9)$ and
$D_1(10)$, and of degree $7$ for $D_1(11)$ and $D_1(12)$.)
Moreover, this ``bigger" factor is always the same for $D_1(2n-1)$ and
$D_1(2n)$. Not only that, the quotient which is predicted in
\eqref{eq:fmu} is at the same time a building block in the result 
of the evaluation
of the determinant which enumerates the cyclically symmetric and
descending plane partitions (see 
\machSeite{AndrAO}\cite{AndrAO}). 
All this begs for an explanation in terms
of a factorisation of the matrix of which the determinant is taken
from. In fact, for the plane partition matrix there is such a
factorisation, due to Mills, Robbins and Rumsey
\machSeite{MiRRAD}\cite[Theorem~5]{MiRRAD} (see also 
\machSeite{KratBN}\cite[Theorem~36]{KratBN}).
The question is whether there is a similar one for the matrix in
\eqref{eq:desc-var}. 

Inspired by this conjecture and by the variations in
\machSeite{CiEKAA}\cite[Theorems~11--13]{CiEKAA} (see
\machSeite{KratBN}\cite[Theorem~35]{KratBN}) on Andrews' original
determinant evaluation in
\machSeite{AndrAN}\cite{AndrAN},
Guoce Xin (private communication) 
observed that, if we change the sign in front of the
Kronecker delta in \eqref{eq:desc-var}, then the resulting determinant
factors completely into linear factors.

\begin{Conjecture} \label{conj:Xin1}
Let $\mu$ be an indeterminate and $n$ be a non-negative integer.
The determinant
\begin{equation} \label{eq:Xin-det1} 
\det_{0\le i,j\le n-1}\left(-\delta_{ij}+\binom {\mu+i+j}{j+1}\right)
\end{equation}
is equal to
\begin{multline} \label{eq:Xin-erg1a} 
            (-1)^{n/2}2^{n(n+2)/4}
\frac {\(\frac \mu2\)_{n/2}} {\(\frac n2\)!}        
  \(  \prod_{i=0}^{(n-2)/2}\frac {i!^2}{(2i)!^2}
\)\\
\times
\(  \prod_{i=0}^{\fl{(n-4)/4}}  \(\frac\mu2+3i+\frac52\)_{(n-4i-2)/2}^2
\(-\frac\mu2-\frac{3n}2+3i+3\)_{(n-4i-4)/2}^2\)
\end{multline}
if $n$ is even, and it is equal to
\begin{multline} \label{eq:Xin-erg1b} 
            (-1)^{(n-1)/2}2^{(n+3)(n+1)/4}\,\(\frac {\mu-1} {2}\)_{(n+1)/2}
        \( \prod_{i=0}^{(n-1)/2}\frac{i!\,(i+1)!} {(2i)!\,(2i+2)!}\)\\
\times
\(\prod_{i=0}^{\fl{(n-3)/4}}\(\frac\mu2+3i+\frac52\)_{(n-4i-3)/2}^2
\(-\frac\mu2-\frac{3n}2+\frac32+3i\)_{(n-4i-1)/2}^2\)
\end{multline}
if $n$ is odd.
\end{Conjecture}

In fact, it seems that also the ``next" determinant, the determinant
where one replaces $j+1$ at the bottom of the binomial coefficient in
\eqref{eq:Xin-det1} by $j+2$ factors completely when $n$ is odd. 
(It does not when $n$ is even, though.)

\begin{Conjecture} \label{conj:Xin2}
Let $\mu$ be an indeterminate. For any odd non-negative integer $n$
there holds
\begin{multline} \label{eq:Xin-det2} 
\det_{0\le i,j\le n-1}\left(-\delta_{ij}+\binom
{\mu+i+j}{j+2}\right)\\=
            (-1)^{(n-1)/2}2^{(n-1)(n+5)/4}(\mu+1)
\frac {\(\frac \mu2-1\)_{(n+1)/2}} {\(\frac{n+1}2\)!}            
\(\prod_{i=0}^{(n-1)/2}\frac {i!^2} {(2i)!^2}
       \(\frac\mu2+3i+\frac32\)_{(n-4i-1)/2}^2\)\\
\times
\( \prod_{i=0}^{\fl{(n-3)/4}}
\(-\frac\mu2-\frac{3n}2+3i+\frac52\)_{(n-4i-3)/2}^2\) .
\end{multline}
\end{Conjecture}

For the combinatorialist
I add that all the determinants in Problem~\ref{prob:7}
and Conjectures~\ref{conj:Xin1} and \ref{conj:Xin2} count 
certain rhombus tilings, as do the original determinants in
\machSeite{AndrAN}%
\machSeite{AndrAO}%
\machSeite{CiEKAA}\cite{AndrAN,AndrAO,CiEKAA}.

Alain Lascoux (private communication)
did not understand why we should stop here, and he
hinted at a parametric family of determinant evaluations into which the
case of odd $n$ of Conjecture~\ref{conj:Xin1} is embedded as a special
case.

\begin{Conjecture} \label{conj:Xin3}
Let $\mu$ be an indeterminate. For any odd non-negative integers $n$
and $r$ there holds
\begin{multline} \label{eq:Xin-det3} 
\det_{0\le i,j\le n-1}\left(-\delta_{i,j+r-1}+\binom
{\mu+i+j}{j+r}\right)\\\kern-1pt
=
            (-1)^{(n-r)/2}2^{(n^2+6n-2nr+r^2-4r+2)/4}
\(
\prod _{i=0} ^{r-2}i!\)
\(
\prod _{i=0} ^{(r-3)/2}\frac {(n-2i-2)!^2} {\(\frac {n-2i-3} {2}\)!^2\,
(n+2i)!\, (n+2i+2)!}\)\kern-1pt\\
\times
(\mu-r)\(\frac {\mu+1} {2}\)_{(n-r)/2}\(
\prod _{i=1} ^{r-1}(m-r+i)_{n+r-2i+1}\)
        \( \prod_{i=0}^{(n-1)/2}\frac{i!\,(i+1)!} {(2i)!\,(2i+2)!}\)\\
\times
\(\prod_{i=0}^{\fl{(n-r-2)/4}}\(\frac\mu2+3i+r+\frac32\)_{(n-4i-r-2)/2}^2
\(-\frac\mu2-\frac{3n}2+\frac {r} {2}+3i+1\)_{(n-4i-r)/2}^2\).
\end{multline}
\end{Conjecture}

\medskip
The next binomial determinant that I want to mention is, strictly
speaking, not a determinant but a Pfaffian (see \eqref{eq:Pfaff}
for the definition). 
While doing {\it $(-1)$-enumeration of sef-complementary plane partitions},
Eisenk\"olbl 
\machSeite{EisTAF}\cite{EisTAF} encountered an, I admit, complicated looking
Pfaffian,
\begin{equation} \label{eq:mn} 
\underset{1\le i,j\le a}\Pf\big(M(m_1,m_2,n_1,n_2,a,b)\big),
\end{equation}
where $a$ is even and $b$ is odd, and where
\begin{multline*}
M_{ij}(m_1, m_2, n_1, n_2, a, b) \\ 
=  \sum_{l=1} ^{\fl{(a + b - 1)/4)}} 
(-1)^{i + j} \(\binom {n_1}{ (b - 1)/2 + \fl{(i - 1)/2} - l + 1}
              \binom {m_1} { -a/2 + \fl{(j - 1)/2} + l} \right.\\
           - 
\left.     \binom {n_1} {(b - 1)/2 + \fl{(j - 1)/2} - l + 1} 
              \binom {m_1} {-a/2 + \fl{(i - 1)/2} + l}\) \\
+ 
 \sum _{l= 1} ^{\cl{(a + b - 1)/4}}
\( \binom {n_2} { (b - 1)/2 + \fl{i/2} - l +1}
            \binom {m_2} {-a/2 + \fl{j/2} + l - 1}
 \right. \\
- 
 \left.         \binom {n_2} {(b - 1)/2 + \fl{j/2} - l+1} 
            \binom {m_2} {-a/2 + \fl{i/2} + l-1}\).
\end{multline*}
Remarkably however, experimentally this Pfaffian, first of all,
factors completely into factors which are linear in the variables
$m_1,m_2,n_1,n_2$, but not only that, there seems to be complete
separation, that is, each linear factor contains only one of
$m_1,m_2,n_1,n_2$. One has the impression that this phenomenon should
have an explanation in a factorisation of the matrix in
\eqref{eq:mn}. However, the task of finding one does not seem to be an
easy one in view of the ``entangledness" of the parameters in the sums
of the matrix entries.

\begin{Problem} \label{prob:6}
\leavevmode
\kern-5pt\footnote{Theresia Eisenk\"olbl has recently
solved this problem in
\machSeite{EisTAG}\cite{EisTAG}.}
Find and prove the closed form evaluation of the Pfaffian in
\eqref{eq:mn}. 
\end{Problem}

Our next determinants can be considered as {\it shuffles} of two
binomial determinants. Let us first consider
\begin{equation} \label{det-ep}
\det_{1\le i,j\le a+m} \begin{pmatrix} \dbinom{b+c+m}{b-i+j}&
\text {\scriptsize $1\le i\le a$}\\
\dbinom{\frac {b+c} {2}}{\frac {b+a} {2}-i+j+\varepsilon}&
\text {\scriptsize $a+1\le i \le a+m$}
\end{pmatrix}.
\end{equation}
In fact, if $\varepsilon=0$, and if $a,b,c$ all have the same parity,
then this is exactly the determinant in \eqref{mat1}, the evaluation
of which proves Theorem~\ref{enum}, as we explained in
Section~\ref{sec:comb}. If $\varepsilon=1/2$ and $a$ has parity
different from that of $b$ and $c$, then the corresponding determinant
was also evaluated in 
\machSeite{CiEKAA}\cite{CiEKAA}, and this evaluation implied the
companion result to Theorem~\ref{enum} that we mentioned immediately
after the statement of the theorem. In the last section of
\machSeite{CiEKAA}\cite{CiEKAA}, 
it is reported that, apparently, there are also nice
closed forms for the determinant in \eqref{det-ep} for $\varepsilon=1$
and $\varepsilon=3/2$, both of which imply as well enumeration
theorems for {\it rhombus tilings of a hexagon with an equilateral 
triangle removed from
its interior} (see Conjectures~1 and 2 in 
\machSeite{CiEKAA}\cite{CiEKAA}). We reproduce the
conjecture for $\varepsilon=1$ here, the one for $\varepsilon=3/2$ is
very similar in form.

\begin{Conjecture} \label{conj:loch1}
Let $a,b,c,m$ be non-negative integers, $a,b,c$ having the same parity. 
Then for $\varepsilon=1$ the determinant in \eqref{det-ep} is equal to
\begin{multline} \label{eq:1-step}
\frac {1} {4}\frac {\h(a + m)\h(b + m)\h(c + m)\h(a + b + c + m)
} 
{\h(a + b + m)\h(a + c + m)\h(b + c + m)
}\\
\times
\frac {\h(m + \left \lceil {\frac{a + b + c}{2}} \right \rceil)
\h(m + \left \lfloor {\frac{a + b + c}{2}} \right \rfloor)
} {\h({\frac{a + b}{2}} + m+1)     \h({\frac{a + c}{2}} + m-1)\h({\frac{b + c}{2}} + m)
}
\\
\times\frac {\h(\left \lceil {\frac{a}{2}} \right \rceil)
\h(\left \lceil {\frac{b}{2}} \right \rceil)
     \h(\left \lceil {\frac{c}{2}} \right \rceil)
     \h(\left \lfloor {\frac{a}{2}} \right \rfloor)\,
     \h(\left \lfloor {\frac{b}{2}} \right \rfloor)\,
     \h(\left \lfloor {\frac{c}{2}} \right \rfloor)\,
} 
{\h({\frac{m}{2}} + \left \lceil {\frac{a}{2}} \right \rceil)\,
     \h({\frac{m}{2}} + \left \lceil {\frac{b}{2}} \right \rceil)\,
     \h({\frac{m}{2}} + \left \lceil {\frac{c}{2}} \right \rceil)\,
\h({\frac{m}{2}} + \left \lfloor {\frac{a}{2}} \right \rfloor)\,
     \h({\frac{m}{2}} + \left \lfloor {\frac{b}{2}} \right \rfloor)\,
     \h({\frac{m}{2}} + \left \lfloor {\frac{c}{2}} \right \rfloor)\,
     }\\
\times
\frac {\h(\frac{m}{2})^2 \h({\frac{a + b + m}{2}})^2 
\h({\frac{a + c + m}{2}})^2 \h({\frac{b + c +
m}{2}})^2
} 
{\h({\frac{m}{2}} + \left \lceil {\frac{a + b + c}{2}} \right \rceil)
\h({\frac{m}{2}} + \left \lfloor {\frac{a + b + c}{2}} \right \rfloor)
\h({\frac{a + b}{2}}-1)\h({\frac{a + c}{2}}+1)\h({\frac{b + c}{2}})
               }P_1(a,b,c,m),
\end{multline}
where $P_1(a,b,c,m)$ is the polynomial given by
$$P_1(a,b,c,m)=\begin{cases} (a+b)(a+c)+2am&\text {if $a$ is
even,}\\
(a+b)(a+c)+2(a+b+c+m)m&\text {if $a$ is odd,}\end{cases}$$
and where $\h(n)$ is the hyperfactorial defined in \eqref{eq:hyperfac}.
\end{Conjecture}

Two other examples of determinants in which the upper part is given
by one binomial matrix, while the lower part is given by a different
one, arose in 
\machSeite{CiKrAC}\cite[Conjectures~A.1 and A.2]{CiKrAC}. Again, both of
them seem to factor completely into linear factors, and both of them
imply enumeration results for {\it rhombus tilings of a certain V-shaped
region}. The right-hand sides of the (conjectured) 
results are the weirdest ``closed" forms in enumeration that I am
aware of.\footnote{No non-trivial simplifications seem to be possible.}
We state just the first of the two conjectures, the other is very similar.

\begin{Conjecture}
Let $x,y,m$ be non-negative integers. Then the determinant
\begin{equation} \label{eq:A.2}
\det_{1\le i,j\le m+y}\left(\left\{\begin{matrix} \binom
{x+i}{x-i+j}&i=1,\dots,m\hfill\\
\binom{x+2m-i+1}{m+y-2i+j+1}&i=m+1,\dots,m+y
\end{matrix}\right\}\right).
\end{equation}
is equal to
\begin{multline} \label{eq:A.1}
\prod _{i=1} ^{m}\frac {(x+i)!} {(x-i+m+y+1)!\,(2i-1)!}
\prod _{i=m+1} ^{m+y}\frac {(x+2m-i+1)!}
{(2m+2y-2i+1)!\,(m+x-y+i-1)!}\\
\times
{2^{\binom {m}2 + \binom y2 }}
      \prod_{i = 1}^{m-1}i!
      \prod_{i = 1}^{y-1}i!
      \prod_{i \ge 0}^{}
        ({ \textstyle x+i+{\frac{3}{2}} }) _{m-2i-1}
      \prod_{i \ge 0}^{}
        ({ \textstyle   x - y+{\frac{5}{2}} + 3 i}) _{ \left \lfloor
         {\frac{3 y}{2}} -{\frac{9 i}{2}}  \right \rfloor-2}\\
\times      \prod_{i \ge 0}^{}
        ({ \textstyle   x + {\frac{3 m}{2}} - y  + \left \lceil
         {\frac{3 i}{2}} \right \rceil+\frac {3} {2}})
      _{ 3 \left \lceil {\frac{y}{2}}
         \right \rceil - \left \lceil
         {\frac{9 i}{2}} \right \rceil -2}
           \prod_{i \ge 0}
        ^{}
        ({ \textstyle {  x+ {\frac{3 m}{2}} - y + \left \lfloor
         {\frac{3 i}{2}} \right \rfloor+2}}) _{ 3 \left \lfloor
         {\frac{y}{2}} \right \rfloor - \left \lfloor {\frac{9 i}{2}} \right
         \rfloor-1}\\
\times      \prod_{i \ge 0}^{}
        ({ \textstyle  x+m - \left \lfloor {\frac{y}{2}} \right
         \rfloor}+i+1) _{  2 \left \lfloor {\frac{y}{2}} \right
         \rfloor-m - 2 i }
     \prod_{i \ge 0}
        ^{}
        ({ \textstyle  x + \left \lfloor {\frac{y}{2}} \right \rfloor+i+2})
         _{m - 2 \left \lfloor {\frac{y}{2}} \right \rfloor-2i-2}
\\\times
{\frac{ \displaystyle
      \prod_{i = 0}^{y}
        ({ \textstyle x - y+3i+1}) _{m + 2 y-4i}
      \prod_{i = 0}^{ \left \lceil {\frac{y}{2}} \right \rceil-1}
        ({ \textstyle x+m - y+i+1}) _{3 y-m-4i}
   }
{\displaystyle
      \prod_{i \ge 0}^{}
        ({ \textstyle x+ {\frac{m}{2}}  - {\frac{y}{2}}+i+1}) _{y-2i}\,
  ({ \textstyle  x + {\frac{m}{2}}-
          {\frac{y}{2}}+i+{\frac{3}{2}}}) _{y-2i-1}   }}\\
\times\frac {\displaystyle
      \prod_{i = 0}^{y}
        ({ \textstyle x+i+2 }) _{2m - 2 i - 1}    }
 {  ({ \textstyle  x + y+2}) _{ m - y-1}  \,(m+x-y+1)_{m+y} }.
\end{multline}
Here, shifted factorials occur with positive as well as with negative
indices. The convention with respect to which these have
to be interpreted is
$$(\alpha)_k:=\begin{cases} \alpha(\alpha+1)\cdots(\alpha+k-1)&\text {if
}k>0,\\
1&\text {if }k=0,\\
1/(\alpha-1)(\alpha-2)\cdots(\alpha+k)&\text {if }k<0.
\end{cases}$$
All products $\prod
_{i\ge0} ^{}(f(i))_{g(i)}$ in \eqref{eq:A.1} 
have to interpreted as the products over
all $i\ge0$ for which $g(i)\ge0$.
\end{Conjecture}

For further conjectures of determinants of shuffles of two binomial
matrices I refer the reader to Conjectures~1--3 in Section~4 of
\machSeite{FuKrAC}\cite{FuKrAC}. 
All of them imply also enumeration results for rhombus
tilings of hexagons. This time, these would be results about the number of
{\it rhombus tilings of a symmetric hexagon with some fixed rhombi on the
symmetry axis}.

\subsection{Determinants of matrices with recursive entries}
Binomial coefficients $\binom {i+j}i$ 
satisfy the basic recurrence of the Pascal triangle,
\begin{equation} \label{eq:Pascal} 
p_{i,j}=p_{i,j-1}+p_{i-1,j}.
\end{equation}
We have seen many determinants of matrices with entries containing
binomial coefficients in the preceding subsection and in 
\machSeite{KratBN}\cite[Sec.~3]{KratBN}.
In \machSeite{BacRAA}\cite{BacRAA}, Bacher reports an experimental
study of determinants of matrices $(p_{i,j})_{0\le i,j\le n-1}$, where
the coefficients $p_{i,j}$ satisfy the recurrence \eqref{eq:Pascal}
(and sometimes more general recurrences), but where the initial
conditions for $p_{i,0}$ and $p_{0,i}$, $i\ge 0$, are different from
the ones for binomial coefficients. He makes many interesting
observations. The most intriguing one says that these determinants
satisfy also a linear recurrence (albeit a much longer one). 
It is intriguing because it
points towards the possibility of {\it automatising determinant
evaluations}\footnote{The reader should recall that the successful
automatisation 
\machSeite{PeWZAA}%
\machSeite{WegsAA}%
\machSeite{WiZeAC}%
\machSeite{ZeilAM}%
\machSeite{ZeilAV}%
\cite{PeWZAA,WegsAA,WiZeAC,ZeilAM,ZeilAV}
of the evaluation of binomial and hypergeometric sums
is fundamentally based on producing recurrences by the computer.}, 
something that several authors (cf.\ e.g.\
\machSeite{AmZeAB}%
\machSeite{KratBN}%
\machSeite{PeWiAA}%
\cite{AmZeAB,KratBN,PeWiAA}) have been aiming at (albeit, with only
limited success up to now). The
conjecture (and, in fact, a generalisation thereof) has been proved by
Petkov\v sek and Zakraj\v sek in 
\machSeite{ZaPeAA}\cite{ZaPeAA}. Still, there remains a large gap to
fill until computers will replace humans doing determinant evaluations.

The paper
\machSeite{BacRAA}\cite{BacRAA} contained as well several
pretty conjectures
on closed form evaluations of special cases of such
determinants. These were subsequently proved in
\machSeite{KratBU}\cite{KratBU}. We state three of them in the 
following three theorems. The first two are proved in
\machSeite{KratBU}\cite{KratBU} by working out the
LU-factorisation (see ``Method~1" in Section~\ref{sec:eval}) 
for the matrices of
which the determinant is computed. The third one is derived by simple
row and column operations.

\begin{Theorem} \label{thm:Pasc1}
Let $(a_{i,j})_{i,j\ge0}$ be the sequence given by the
recurrence
$$a_{i,j}=a_{i-1,j}+a_{i,j-1}+x\,a_{i-1,j-1},\quad \quad i,j\ge1,
$$
and the initial conditions $a_{i,0}=\rho^i$ and $a_{0,i}=\si^i$,
$i\ge0$. Then 
\begin{equation} \label{eq:Pasc1} 
\det_{0\le i,j\le n-1}(a_{i,j})=(1+x)^{\binom
{n-1}2}(x+\rho+\si-\rho\si)^{n-1}.
\end{equation}
\quad \quad \qed
\end{Theorem}

\begin{Theorem} \label{thm:Pasc2}
Let $(a_{i,j})_{i,j\ge0}$ be the sequence given by the
recurrence
$$a_{i,j}=a_{i-1,j}+a_{i,j-1}+x\,a_{i-1,j-1},\quad \quad i,j\ge1,
$$
and the initial conditions $a_{i,i}=0$, $i\ge0$, $a_{i,0}=\rho^{i-1}$ and
$a_{0,i}=-\rho^{i-1}$, $i\ge1$. Then 
\begin{equation} \label{eq:Pasc2} 
\det_{0\le i,j\le 2n-1}(a_{i,j})=(1+x)^{2(n-1)^2}
(x+\rho)^{2n-2}.
\end{equation}
\quad \quad \qed
\end{Theorem}

\begin{Theorem}\label{thm:Pasc3}
Let $(a_{i,j})_{i,j\ge0}$ be the sequence given by the recurrence
$$a_{i,j}=a_{i-1,j}+a_{i,j-1}+x\,a_{i-1,j-1},\quad \quad i,j\ge1,
$$
and the initial conditions $a_{i,0}=i$ and $a_{0,i}=-i$,
$i\ge0$. Then
\begin{equation} \label{eq:Pasc3} 
\det_{0\le i,j\le 2n-1}(a_{i,j})=(1+x)^{2n(n-1)}.
\end{equation}
\quad \quad \qed
\end{Theorem}

Certainly, the proofs in
\machSeite{KratBU}\cite{KratBU} are not very illuminating. Neuwirth
\machSeite{NeuwAE}\cite{NeuwAE} has looked more carefully into the
structure of recursive sequences of the type as those in
Theorems~\ref{thm:Pasc1}--\ref{thm:Pasc3}. Even more generally, he
looks at sequences $(f_{i,j})_{i,j\ge0}$
satisfying the recurrence relation
\begin{equation} \label{eq:rec} 
f_{i,j}=c_{j}f_{i-1,j}+d_jf_{i,j-1}+e_jf_{i-1,j-1},\quad i,j\ge1,
\end{equation}
for some given sequences $(c_j)_{j\ge1}$, $(d_j)_{j\ge1}$,
$(e_j)_{j\ge1}$. He approaches the problem by finding appropriate
{\it matrix decompositions} for the (infinite) matrix 
$(f_{i,j})_{i,j\ge0}$. In two special
cases, he is able to apply his
decomposition results to work out the LU-factorisation of the matrix
$(f_{i,j})_{i,j\ge0}$ explicitly, which then yields an elegant
determinant evaluation in both of these cases. Neuwirth's first 
result 
\machSeite{NeuwAE}\cite[Theorem~5]{NeuwAE}
addresses the case where the initial values $f_{0,j}$ satisfy a first
order recurrence determined by the coefficients $d_j$ from
\eqref{eq:rec}. It generalises Theorem~\ref{thm:Pasc1}. There is no
restriction on the initial values $f_{i,0}$ for $i\ge1$.

\begin{Theorem} \label{thm:Neuw1}
Let $(c_j)_{j\ge1}$, $(d_j)_{j\ge1}$ and $(e_j)_{j\ge1}$ be given
sequences, and let $(f_{i,j})_{i,j\ge0}$ be the doubly indexed 
sequence given by the recurrence
\eqref{eq:rec} and the initial conditions $f_{0,0}=1$ and
$f_{0,j}=d_jf_{0,j-1}$, $j\ge1$. Then
\begin{equation} \label{eq:Neuw1}
\det_{0\le i,j\le n-1}(f_{i,j})=
\prod _{0\le i<j\le n-1} ^{}(e_{i+1}+c_jd_{i+1}).
\end{equation}
\quad \quad \qed
\end{Theorem}

Neuwirth's second
result 
\machSeite{NeuwAE}\cite[Theorem~6]{NeuwAE}
also generalises Theorem~\ref{thm:Pasc1}, but in a different way. 
This time, the initial values $f_{0,j}$, $j\ge1$, are free, whereas
the initial values $f_{i,0}$ satisfy a first order recurrence
determined by
the coefficients $c_j$ from the recurrence \eqref{eq:rec}.
Below, we state its most attractive special case, in which all the
$c_j$'s are identical.

\begin{Theorem} \label{thm:Neuw2}
Let $(d_j)_{j\ge1}$ and $(e_j)_{j\ge1}$ be given
sequences, and let $(f_{i,j})_{i,j\ge0}$ be the doubly indexed 
sequence given by the recurrence
\eqref{eq:rec} and the initial conditions $f_{0,0}=1$ and
$f_{i,0}=cf_{i-1,0}$, $i\ge1$. Then
\begin{equation} \label{eq:Neuw2}
\det_{0\le i,j\le n-1}(f_{i,j})=
\prod _{i=1} ^{n-1}(cd_{i}+e_{i})^{n-i}.
\end{equation}
\quad \quad \qed
\end{Theorem}

\subsection{Determinants for signed permutations} \label{sec:signed}
The next class of determinants that we consider are determinants of
matrices in which rows and columns are indexed by 
{\it elements of reflection groups}
(the latter being groups generated by reflections of hyperplanes in real
$n$-dimensional space; see 
\machSeite{HumpAC}\cite{HumpAC} for more information on these
groups, and, more generally, on Coxeter groups).
The prototypical example of a reflection group is the {\it symmetric group}
$\mathfrak S_n$ of permutations of an $n$-element set. In 
\machSeite{KratBN}\cite{KratBN}, there
appeared two determinant evaluations associated to the symmetric group,
see Theorems~55 and 56 in 
\machSeite{KratBN}\cite{KratBN}. They concerned evaluations of
determinants of the type
\begin{equation} \label{eq:stat} 
\det_{\si,\pi\in \mathfrak S_n}\(q^{\stat(\si\pi^{-1})}\),
\end{equation}
due to Varchenko, Zagier, and Thibon, respectively, in which stat is
the statistic {\it ``number of inversions,"} respectively {\it``major index."}
We know that in many fields of mathematics there exist certain
diseases which are typical for that field. Algebraic combinatorics is
no exception. Here, I am {\it not\/} talking
of the earlier mentioned ``$q$-disease"
(see Footnote~\ref{foot:q}; 
although, due to the presence of $q$ we might also count it as a case
of $q$-disease), but of the disease which manifests itself by the
question ``And what about the other types?"\footnote{\label{foot:root}%
In order
to give a reader who is not acquainted with the language and theory of
reflection groups an idea what this question is referring to, 
I mention that all finite reflection groups have
been classified, each having been assigned a certain ``type." So,
usually one proves something for the symmetric group $\mathfrak S_n$,
which, according to this classification, has type $A_{n-1}$, and
then somebody (which could be oneself) will ask the question ``Can you
also do this for the other types?", meaning whether or not there exists
an analogous result for the other finite reflection groups.} 
So let us ask this question, that is, 
are there theorems similar to the two theorems
which we mentioned above for other
reflection groups? 

So, first of all we need analogues of the statistics ``number of
inversions" and ``major index" for other reflection groups. Indeed,
these are available in the literature. The analogue of ``number of
inversions" is the so-called {\it length} of an element in a Coxeter
group (see 
\machSeite{HumpAC}\cite{HumpAC} for the definition). 
As a matter of fact, a closed form evaluation
of the determinant \eqref{eq:stat}, where $\mathfrak S_n$ is replaced
by any finite or affine reflection group, and where stat is the
length, is known (and was already implicitly mentioned in
\machSeite{KratBN}\cite{KratBN}). This result is due to Varchenko
\machSeite{VarcAC}\cite[Theorem~(1.1), 
where $a(H)$ is specialised to $q$]{VarcAC}.
His result is actually much more general, as it is valid for real
hyperplane arrangements in which each hyperplane is assigned a different
weight. I will not state it here explicitly because I do
not want to go through the definitions and notations which would be
necessary for doing that.

So, what about analogues of the ``major index" for other reflection
groups? These are also available, and there are in fact several of
them. The first person to introduce a major index for reflection groups
other than the symmetric groups was Reiner in 
\machSeite{ReivAC}\cite{ReivAC}.
He proposed a major index for the {\it hyperoctahedral group} $B_n$, which
arose naturally in his study of {\it $P$-partitions for signed
posets}. The elements of $B_n$ are often called {\it signed
permutations}, and they are all elements of the form 
$\pi_1\pi_2\dots\pi_n$, where
$\pi_i\in\{\pm1,\pm2,\dots,\pm n\}$, $i=1,2,\dots,n$, and where 
$\vert\pi_1\vert\vert\pi_2\vert\dots\vert\pi_n\vert$ is a permutation
in $\mathfrak S_n$. To define their multiplicative structure, it is
most convenient to view $\pi=\pi_1\pi_2\dots\pi_n$ as a linear operator on
$\R^n$ acting by permutation and sign changes of the co-ordinates. To
be precise, the action is given by
$\pi(e_i)=(\sgn\pi_i)e_{\vert\pi_i\vert}$, where $e_i$ is the $i$-th
standard basis vector in $\R^n$, $i=1,2,\dots,n$. The multiplication
of two signed permutations is then simply the composition of the
corresponding linear operators. 

The major index $\maj_B\pi$ of an element $\pi\in B_n$
which Reiner defined is, as in the symmetric group case, the sum of
all positions of {\it descents} in $\pi$. (There is a natural notion
of ``descent" for any Coxeter group.) Concretely, it is
$$\maj_B\pi:=\chi(\pi_n<0)+
\sum _{i=1} ^{n-1}i\cdot\chi(\pi_i>_B\pi_{i+1}),$$
where we impose the order $1<_B2<_B\cdots<_Bn<_B-n<_B\cdots<_B-2<_B-1$
on our ground set $\{\pm1,\pm2,\dots,\pm n\}$,
and where $\chi(\mathcal A)=1$ if $\mathcal A$ is true and
$\chi(\mathcal A)=0$ otherwise. 

There is overwhelming computational
evidence\footnote{\label{foot:maj}The 
reader may wonder what this computational evidence
could actually be. After all, we are talking about a determinant of a
matrix of size $2^nn!$. More concretely, for $n=1,2,3,4,5$ these are
matrices of size $2$, $8$, $48$, $384$, $3840$, respectively.
While {\sl Maple} or {\sl Mathematica} have no problem to compute
these determinants for $n=1$ and $n=2$, it takes already considerable time to
do the computation for $n=3$, and it is, of course, completely hopeless to
let them compute the one for $n=4$, a determinant of a matrix of size
$384$ which has polynomial entries (cf.\
Footnote~\ref{foot:kompl}). However, the results for $n=1,2,3$
already ``show" that the determinant will factor completely into
factors of the form $1-q^i$, $i=1,2,\dots,2n$. One starts to expect 
the same to be true for higher $n$. To get a formula for $n=4$, one
would then apply the tricks explained in
Footnote~\ref{foot:tricks}. That is, one specialises $q$ to $4$, at
which value the first $8$ cyclotomic polynomials (in fact, even more)
are clearly distinguishable by their prime factorisations, and one computes
the determinant. The exponents of the various factors $1-q^i$ can then
be extracted from the exponents of the prime factors in the prime
factorisation of the determinant with $q=4$. Unfortunately, the data
collected for $n=1,2,3,4$ do not suffice to come up with a
guess, and, on the other hand, {\sl Maple} and {\sl Mathematica} will
certainly be incapable to compute a determinant of a matrix of size
$3840$ (which, just to store it on the disk, occupies already 10
megabytes \dots). So then, what did I mean when I said that the conjecture is
based on data including $n=5$? This turned out to become a ``test
case" for {\sl LinBox}, a C++ template library for exact
high-performance linear algebra 
\machSeite{DGGGAA}\cite{DGGGAA}, which is freely available under
{\tt http://linalg.org}. To be honest, 
I was helped by Dave Saunders and
Zhendong Wan (two of the developers) who applied {\sl LinBox}
to do rank and Smith normal form computations for the specialised matrix with
respect to various prime powers (each of which taking several hours). 
The specific computational approach that worked here is quite recent
(thus, it came just in time for our purpose), 
and is documented in 
\machSeite{SaWaAA}\cite{SaWaAA}.
The results of the computations made it possible to come
up with a ``sure" prediction for the exponents with which the various
prime factors occur in the prime factorisation of the specialised determinant.
As in the case $n=4$, the exponents of the various factors $1-q^i$,
$i=1,2,\dots,2n$ can then easily be extracted. (The guesses were 
subsequently also
tested with special values of $q$ other than $q=4$.)}
that the ``major-de\-ter\-min\-ant" for $B_n$, i.e., the
determinant \eqref{eq:stat} with $\text{stat}=\maj_B$ and with
$\mathfrak S_n$ replaced by $B_n$, factors completely into cyclotomic
polynomials. 

\begin{Conjecture} \label{prob:1}
For any positive integer $n$, we have
\begin{equation} \label{eq:maj-Bn} 
\det_{\si,\pi\in B_n}\(q^{\maj_B(\si\pi^{-1})}\)=
\prod_{i=1}^n (1-q^{2i})^{2^{n-1}n!/i}
\prod_{i=2}^n (1-q^i)^{2^n n!(i-1)/i}.
\end{equation}
\end{Conjecture}


A different major index for $B_n$ was proposed by Adin and Roichman in
\machSeite{AdRoAC}\cite{AdRoAC}. 
It arises there naturally in a combinatorial study of
{\it polynomial algebras which are diagonally invariant under $B_n$}.
(In fact, more generally, {\it wreath products} of the form $C_m\wr
\mathfrak S_n$, where $C_m$ is the cyclic group of order $m$, 
and their diagonal actions on polynomial algebras are
studied in 
\machSeite{AdRoAC}\cite{AdRoAC}. These groups are also sometimes called {\it
generalised reflection groups}.
In this context, $B_n$ is the special case
$C_2\wr\mathfrak S_n$.) 
If we write $\neg\pi$ for the number of $i$ for which $\pi_i$ is
negative, then the {\it flag-major index} fmaj of Adin and Roichman is
defined by 
\begin{equation} \label{eq:fmaj-def} 
\fmaj\pi:=2\maj_A\pi+\neg\pi,
\end{equation}
where $\maj_A$ is the ``ordinary" major index due to MacMahon,
$$\maj_A\pi:=\sum _{i=1} ^{n-1}i\cdot\chi(\pi_i>\pi_{i+1}).$$
If one now goes to the computer and calculates the determinant on the
left-hand side of
\eqref{eq:maj-Bn} with $\maj_B$ replaced by $\fmaj$ 
for $n=1,2,3,4,5$ (see
Footnote~\ref{foot:maj} for the precise meaning of ``calculating the
determinant for $n=1,2,3,4,5$"), then again the results factor
completely into cyclotomic factors. Even more generally, it seems that
one can treat the two parts on the right-hand side of
\eqref{eq:fmaj-def}, that is ``major index" and
``number of negative letters," separately.

\begin{Conjecture} \label{prob:2}
For any positive integer $n$, we have
\begin{equation} \label{eq:fmaj-Bn} 
\det_{\si,\pi\in B_n}\(q^{\maj_A(\si\pi^{-1})}p^{\neg(\si\pi^{-1})}\)=
\prod _{i=1} ^{n}(1-p^{2i})^{2^{n-1}n!/i}
\prod _{i=2} ^{n}(1-q^i)^{2^{n}n!(i-1)/i}.
\end{equation}
\end{Conjecture}


I should remark that Adin and Roichman have shown in 
\machSeite{AdRoAC}\cite{AdRoAC}
that the statistics fmaj is equidistributed with the statistics length
on $B_n$. However, even in the case where we just look at the
flag-major determinant (that is, the case where $q=p^2$ in
Conjecture~\ref{prob:2}), this does not seem to help. (Neither length nor
flag-major index satisfy a simple law with respect to multiplication
of signed permutations.) In fact, from the
data one sees that the flag-major determinants are different from the
length determinants (that is, the determinants \eqref{eq:stat}, where
$\mathfrak S_n$ is replaced by $B_n$ and stat is flag-major,
respectively length). 

Initially, I had my program wrong, and,
instead of taking the (ordinary) major index $\maj_A$ of the signed
permutation $\pi=\pi_1\pi_2\dots\pi_n$ in \eqref{eq:fmaj-def}, 
I computed taking the major index of the {\it absolute value}
of $\pi$. This absolute value is obtained
by forgetting all signs of the letters of $\pi$, that is,
writing $\vert\pi\vert$ for the absolute value of $\pi$,
$\vert\pi\vert=\vert\pi_1\vert\,\vert\pi_2\vert\,\dots\,\vert\pi_n\vert$.
Curiously, it seems that also this ``wrong"
determinant factors nicely. (Again, the evidence for this
conjecture is based on data which were obtained in the way described in
Footnote~\ref{foot:maj}.)

\begin{Conjecture} \label{prob:3}
For any positive integer $n$, we have
\begin{equation} \label{eq:fabsmaj-Bn} 
\det_{\si,\pi\in B_n}\(q^{\maj\vert\si\pi^{-1}\vert}p^{\neg(\si\pi^{-1})}\)=
(1-p^{2})^{2^{n-1}n\cdot n!}
\prod _{i=2} ^{n}(1-q^i)^{2^{n}n!(i-1)/i}.
\end{equation}
\end{Conjecture}


Since, as I indicated earlier, Adin and Roichman actually define a
flag-major index for wreath
products $C_m\wr\mathfrak S_n$, a question that suggests itself is
whether or not we can expect closed product formulae for the
corresponding determinants. Clearly, since we are now dealing with
determinants of the size $m^nn!$, computer computations will
exhaust our computer's resources even faster if $m>2$. The
calculations that I was able to do suggest strongly that there is indeed an
extension of the statement in Conjecture~\ref{prob:2} to the case
of arbitrary $m$, if one uses the definition of major index 
and the ``negative" statistics for $C_m\wr\mathfrak S_n$ as
in \machSeite{AdRoAC}\cite{AdRoAC}. (See 
\machSeite{AdRoAC}\cite[Section~3]{AdRoAC} 
for the definition of the major index.
The sum on the right-hand side of 
\machSeite{AdRoAC}\cite[(3.1)]{AdRoAC} must be taken as the extension of
the ``negative" statistics neg to $C_m\wr\mathfrak S_n$.)

\begin{Problem} \label{prob:4}
Find and prove the closed form evaluation of 
\begin{equation} \label{eq:maj-wreath} 
\det_{\si,\pi\in C_m\wr \mathfrak S_n}\(q^{\maj(\si\pi^{-1})}
p^{\neg(\si\pi^{-1})}\),
\end{equation}
where $\maj$ and $\neg$ are the extensions to $C_m\wr \mathfrak S_n$  
of the statistics $\maj_A$ and
$\neg$ in Conjecture~{\em\ref{prob:2}}, as described in
the paragaph above.
\end{Problem}

Together with Brenti, Adin and Roichman proposed another major
statistics for signed permutations in 
\machSeite{AdBRAA}\cite{AdBRAA}. They call it the
{\it negative major index}, denoted nmaj, and it is defined as the sum of the
ordinary major index and the sum of the absolute values of the 
negative letters, that is,
$$\nmaj\pi:=\maj_A\pi+\sneg\pi,$$
where $\sneg\pi:=-\sum _{i=1} ^{n}\chi(\pi_i<0)\pi_i$.
Also for this statistics, the corresponding determinant seems to
factor nicely. In fact, it seems that one can again treat the two
components of the definition of the statistics, that is, ``major index" and
``sum of negative letters," separately. (Once more, 
the evidence for this
conjecture is based on data which were obtained in the way described in
Footnote~\ref{foot:maj}.)

\begin{Conjecture} \label{prob:8}
For any positive integer $n$, we have
\begin{equation} \label{eq:nmaj-Bn} 
\det_{\si,\pi\in
B_n}\(q^{\maj_A(\si\pi^{-1})}p^{\sneg(\si\pi^{-1})}\)=
\prod _{i=1} ^{n}(1-p^{2i^2})^{2^{n-1}n!/i}
\prod _{i=2} ^{n}(1-q^i)^{2^{n}n!(i-1)/i}.
\end{equation}
\end{Conjecture}

If one compares the (conjectured) result with the (conjectured) one 
for the ``flag-major determinant" in Conjecture~\ref{prob:2}, 
then one notices the somewhat mind-boggling fact that
one obtains the right-hand side of \eqref{eq:nmaj-Bn} from the one
of \eqref{eq:fmaj-Bn} by simply replacing (in the factored form
of the latter) $1-p^{2i}$ by $1-p^{2i^2}$, everything else, the
exponents, the ``$q$-part", is identical. It is difficult to imagine
an intrinsic explanation why this should be the case.

Since Thibon's proof of the evaluation of the determinant
\eqref{eq:stat} with stat being the (ordinary) major index for
permutations (see 
\machSeite{KratBN}\cite[Appendix~C]{KratBN})
involved the {\em descent algebra} of the symmetric group, viewed
in terms of {\em non-commutative symmetric functions}, one might
speculate that to prove Conjectures~\ref{prob:1}--\ref{prob:3} and
\ref{prob:8} it may be necessary to work with $B_n$ versions of descent
algebras (which exist, see
\machSeite{SoloAA}\cite{SoloAA}) 
and non-commutative symmetric functions (which also exist, see
\machSeite{ChoCAA}\cite{ChoCAA}).

Adriano Garsia points out that all the determinants in \eqref{eq:stat},
\eqref{eq:maj-Bn}, \eqref{eq:fmaj-Bn}--\eqref{eq:nmaj-Bn} are special
instances of {\it group determinants}. (See the excellent survey article
\machSeite{LamTAA}\cite{LamTAA} for information on group
determinants.)
The main theorem on group determinants, due to Frobenius, says that
a general group determinant factorises into irreducible factors, each
of which corresponding to an irreducible representation of the group,
and the exact exponent to which the irreducible factor is raised is
the degree of the corresponding irreducible representation. In view of
this, Garsia poses the following problem, a solution of which
would refine the above
conjectures, Problem~\ref{prob:4}, 
and the earlier mentioned results of Varchenko, Zagier,
and Thibon.

\begin{Problem} \label{prob:Frob}
For each of the above special group determinants, determine
the closed formula for the value of the irreducible factor corresponding
to a fixed irreducible representation of the involved group
{\em(}$\mathfrak S_n$, $B_n$, $C_m\wr \mathfrak S_n$, respectively{\em)}.
\end{Problem}

It seems that a solution to this problem has not even been worked out
for the determinant which is the subject of the results of Varchenko
and Zagier, that is, for the determinant
\eqref{eq:stat} with stat being the number of inversions.

For further work on statistics for (generalised) reflection groups
(thus providing further prospective candidates for forming interesting
determinants), I refer the reader to
\machSeite{AdBRAA}%
\machSeite{AdBRAB}%
\machSeite{BagnAA}%
\machSeite{BernAA}%
\machSeite{BernAB}%
\machSeite{BiagAA}%
\machSeite{BiCaAA}%
\machSeite{FoHaAM}%
\machSeite{FoHaAN}%
\machSeite{FoHaAO}%
\machSeite{HaLRAA}%
\machSeite{ReRoAA}%
\cite{AdBRAA,AdBRAB,BagnAA,BernAA,BernAB,BiagAA,BiCaAA,FoHaAM,FoHaAN,FoHaAO,HaLRAA,ReRoAA}.
I must report that, somewhat disappointingly, it seems that the
various major indices proposed for the group $D_n$ of even signed
permutations (see 
\machSeite{BiagAA}%
\machSeite{BiCaAA}%
\machSeite{ReivAC}%
\cite{BiagAA,BiCaAA,ReivAC}) 
apparently do not give rise to determinants in the same
way as above that have nice product formulae. This remark seems to 
also apply to determinants formed in an anlogous way by using
the various statistics proposed for the alternating group in
\machSeite{ReRoAB}\cite{ReRoAB}.

\subsection{More poset and lattice determinants} \label{sec:poset}
Continuing the discussion of determinants which arise under the
influence of the above-mentioned ``reflection group disease," we turn
our attention to two miraculous
determinants which were among the last things Rodica Simion was able
to look at. Some of her considerations in this direction are reported in
\machSeite{SchFAA}\cite{SchFAA}.

The first of the two is a determinant of a matrix the rows and columns of
which are indexed by {\it type $B$ non-crossing partitions}. 
This determinant is inspired by the evaluation of 
an analogous one for {\it ordinary}
non-crossing partitions (that is, in ``reflection group language," 
{\it type $A$ non-crossing
partitions}), due to Dahab 
\machSeite{DahaAA}\cite{DahaAA} and Tutte 
\machSeite{TuttAC}\cite{TuttAC} (see
\machSeite{KratBN}\cite[Theorem~57, (3.69)]{KratBN}). Recall (see 
\machSeite{StanAP}\cite[Ch.~1 and
3]{StanAP} for more information) that a {\it {\rm(}set{\rm)}
partition}
of a set $S$ is a collection $\{B_1,B_2,\dots,B_k\}$ of
pairwise disjoint non-empty subsets of $S$ such that
their union is equal to $S$. The subsets $B_i$ are also called {\it
blocks} of the partition. One partially orders 
partitions by refinement. With respect to this partial order, the 
partitions form a lattice. We write $\pi\lor_{A}\ga$ (the $A$ stands
for the fact that, in ``reflection group language", we are looking at
``type $A$ partitions") for
the join of $\pi$ and $\ga$ in this lattice. Roughly speaking, the
join of $\pi$ and $\ga$ is formed by considering altogether
all the blocks of $\pi$ and $\ga$. Subsequently,
whenever we find two blocks which have a non-empty
intersection, we merge them into a bigger block, and we keep doing
this until all the (merged) blocks are pairwise disjoint.

If $S=\{1,2,\dots,n\}$, we call a
partition {\it non-crossing} if for any $i<j<k<l$ the elements $i$ and
$k$ are in the same block at the same time as the elements $j$ and $l$
are in the same block only if these two blocks are the same. (I refer
the reader to 
\machSeite{SimiAD}\cite{SimiAD} for a survey on non-crossing partitions.)

Reiner 
\machSeite{ReivAD}\cite{ReivAD} 
introduced non-crossing partitions
in type $B$. Partitions of type $B_n$ are (ordinary) partitions of 
$\{1,2,\dots,n,-1,-2,\dots,-n\}$ with the property that whenever $B$
is a block then so is $-B:=\{-b:b\in B\}$, and that there is at most
one block $B$ with $B=-B$. A block $B$ with $B=-B$, if present, 
is called the {\it zero-block} of the partition. We denote the set of
all type $B_n$ partitions by $\Pi_n^B$, the number of zero blocks of a
partition $\pi$by $\zbk \pi$, and we write $\nzbk\pi$ for half of the
number of the non-zero blocks. Type $B_n$
non-crossing partitions are a subset of type $B_n$ partitions.
Imposing the order $1<2<\dots<n<-1<-2<\dots<-n$ on our ground-set,
the definition of type $B_n$ non-crossing partitions is identical with
the one for type $A$ non-crossing partitions, that is, given this
order on the ground-set, a $B_n$ partition is called {\it non-crossing} 
if for any $i<j<k<l$ the elements $i$ and
$k$ are in the same block at the same time as the elements $j$ and $l$
are in the same block only if these two blocks are the same. 
We write $\NC_n^B$ for the set of all $B_n$ non-crossing partitions.

The determinant defined by type $B_n$ non-crossing partitions that
Simion tried to evaluate was the one in 
\eqref{eq:Simion1} below.\footnote{In fact, instead of $\lor_A$, the
``ordinary" join, she
used the join in the type $B_n$ partition lattice $\Pi_n^B$. However,
the number of {\it non-zero} blocks will be the same regardless of whether
we take the join of two type $B_n$ non-crossing partitions with
respect to ``ordinary" join or with respect to ``type $B_n$" join. This
is in contrast to the numbers of zero blocks, which can differ
largely. (To be more precise, one way to form the ``type $B_n$" join is
to first form the ``ordinary" join, and then merge all zero blocks
into one big block.) The reason that I insist on using $\lor_A$ is
that this is crucial for the more general
Conjecture~\ref{conj:Simion2}. To tell the truth, the discovery of the
latter conjecture is due to a programming error on my behalf (that is,
originally I aimed to program the ``type $B$" join, but it happened 
to be the ``ordinary" join \dots).} 
The use of the ``type $A$" join $\lor_A$ for two type $B$ non-crossing
partitions in \eqref{eq:Simion1} may seem strange. However, this is
certainly a well-defined operation. The result may neither be a type
$B$ partition nor a non-crossing one, it will just be an ordinary partition
of the ground-set $\{1,2,\dots,n,-1,-2,\dots,-n\}$. We extend the
notion of ``zero block" and ``non-zero block" to these objects in the
obvious way.
Simion observed that, as in the case of the
analogous type $A_n$ partition determinant due to Dahab and Tutte, it
factors apparently completely 
into factors which are Chebyshev polynomials. Based on some
additional numerical calculations,\footnote{\label{foot:Simion}%
Evidently, more than 
five years later, thanks to technical progress since then,
one can go much farther when doing computer
calculations. The evidence for Conjecture~\ref{conj:Simion1} which I
have is based on, similar to the conjectures and calculations
on determinants for signed permutations in Subsection~\ref{sec:signed} (see
Footnote~\ref{foot:maj}), the exact form of the determinants for $n=1,2,3,4$,
which were already computed by Simion, and, essentially, the exact
form of the determinants for $n=5$ and $6$. By ``essentially" 
I mean, as earlier, that I computed the determinant for many special values of
$q$, which then let me make a guess on the basis of comparison of the prime
factors in the factorised results with the prime factors of the
candidate factors, that is the irreducible factors of the Chebyshev
polynomials. 
Finally, for guessing the general form of the exponents, the available
data were not sufficient for {\tt Rate} (see
Footnote~\ref{foot:Rate}). So I consulted the
fabulous {\it On-Line Encyclopedia of Integer Sequences}
({\tt http://www.research.att.com/\~{}njas/sequences/Seis.html}),
originally created by Neil Sloane and Simon Plouffe 
\machSeite{SlPlAA}\cite{SlPlAA}, 
and since many years continuously further developed by Sloane and his
team 
\machSeite{SloaAA}\cite{SloaAA}. An appropriate
selection from the results turned up by the Encyclopedia then 
led to the exponents on the right-hand side of
\eqref{eq:Simion1}.} I propose the following conjecture.

\begin{Conjecture} \label{conj:Simion1}
For any positive integer $n$, we have
\begin{equation} \label{eq:Simion1} 
\det_{\pi,\ga\in\NC_n^B}\(q^{\nzbk (\pi\lor_{A}\ga)}\)=
\prod _{i=1} ^{n}\left(\dfrac {U_{3i-1}(\sqrt q/2)} 
{U_{i-1}(\sqrt q/2)}\)^{\binom {2n}{n-i}},
\end{equation}
where $U_{m}(x):=\sum _{j\ge0} ^{}(-1)^j\binom
{m-j}j(2x)^{m-2j}$ is the $m$-th
{\em Chebyshev polynomial of the second kind}.
\end{Conjecture}

If proved, this would solve Problem~1 in 
\machSeite{SchFAA}\cite{SchFAA}. It would also
solve Problem~2 from 
\machSeite{SchFAA}\cite{SchFAA},
because $U_{m-1}(\sqrt q/2)$ is, up to multiplication by a power of
$q$, equal to the product $\prod _{j\mid m} ^{}f_j(q)$, where the
polynomials $f_j(q)$ are the ones of 
\machSeite{SchFAA}\cite{SchFAA}. A simple
computation then shows that, when the right-hand side product of
\eqref{eq:Simion1} is expressed in terms of the $f_j(q)$'s, one obtains
\begin{equation} \label{eq:fprod} 
\prod _{k=1} ^{n}f_{3k}(q)^{e_{n,k}},
\end{equation}
where
$$e_{n,k}=\underset {\ell\not\equiv 0\,(\text{mod }3)}
{\sum _{\ell=1} ^{\fl{n/k}}}\binom
{2n}{n-\ell k}.$$
This agrees with the data in 
\machSeite{SchFAA}\cite{SchFAA} and with the further ones I
have computed (see Footnote~\ref{foot:Simion}).

Even more seems to be true. The following conjecture predicts the
evaluation of the more general determinant where we also keep track of
the zero blocks.

\begin{Conjecture} \label{conj:Simion2}
For any positive integer $n$, we have
\begin{equation} \label{eq:Simion2} 
\det_{\pi,\ga\in\NC_n^B}\(q^{\nzbk (\pi\lor_{A}\ga)}z^{\zbk
(\pi\lor_{A}\ga)}\)=z^{\frac {1} {2}\binom {2n}{n}}
\prod _{i=1} ^{n}\Big(2T_{2i}(\sqrt q/2)+2-z\Big)^{\binom {2n}{n-i}},
\end{equation}
where $T_{m}(x):=\frac {1} {2}\sum _{j\ge0} ^{}(-1)^j\frac {m} {m-j}\binom
{m-j}j(2x)^{m-2j}$ is the $m$-th
{\em Chebyshev polynomial of the first kind}.
\end{Conjecture}

Again, the conjecture is supported by extensive numerical
calculations. It is not too difficult to show, by using some
identities for Chebyshev polynomials, that
Conjecture~\ref{conj:Simion2} implies Conjecture~\ref{conj:Simion1}.

The other determinant which Simion looked at (cf.\ 
\machSeite{SchFAA}\cite[Problem~9ff]{SchFAA}), 
was the $B_n$ analogue
of a determinant of a matrix the rows and columns of which are indexed
by {\it non-crossing matchings}, due to Lickorish 
\machSeite{LickAA}\cite{LickAA}, and 
evaluated by Ko and Smolinsky
\machSeite{KoSmAA}\cite{KoSmAA} and independently by Di~Francesco 
\machSeite{DiFrAA}\cite{DiFrAA}
(see 
\machSeite{KratBN}\cite[Theorem~58]{KratBN}). As we may regard (ordinary)
non-crossing matchings as partitions all the blocks of which consist of two
elements, we {\it define} a {\it $B_n$ non-crossing matching} to be a $B_n$
non-crossing partition all the blocks of which consist of two elements.
We shall be concerned with $B_{2n}$ non-crossing matchings, which we
denote by $\match(2n)$. With this notation, the following seems to be
true. 

\begin{Conjecture} \label{conj:Simion3}
For any positive integer $n$, we have
\begin{equation} \label{eq:Simion3} 
\det_{\pi,\ga\in\match(2n)}\(q^{\nzbk (\pi\lor_{A}\ga)}z^{\zbk
(\pi\lor_{A}\ga)}\)=
\prod _{i=1} ^{n}\Big(2T_{2i}(q/2)+2-z^2\Big)^{\binom {2n}{n-i}}.
\end{equation}
\end{Conjecture}

The reader should notice the remarkable fact that,
in the case that Conjectures~\ref{conj:Simion2}
and \ref{conj:Simion3} are true, the right-hand side of \eqref{eq:Simion3}
is, up to a power of $z$, equal to the right-hand side of
\eqref{eq:Simion2} with $q$ replaced by $q^2$ and $z$ replaced by
$z^2$. An intrinsic explanation why this should be the
case is not known.
An analogous relation between the determinants of 
Tutte and of Lickorish, respectively, was observed, and proved, in
\machSeite{CoSSAA}\cite{CoSSAA}. Also here, no intrinsic
explanation is known.

The reader is referred to 
\machSeite{SchFAA}\cite{SchFAA} for further open problems
related to the determinants in
Conjectures~\ref{conj:Simion1}--\ref{conj:Simion3}. 
Finally, it may also be worthwhile to look at determinants defined
using $D_n$ non-crossing partitions and non-crossing matchings,
see 
\machSeite{AtReAA}\cite{AtReAA} and
\machSeite{ReivAD}\cite{ReivAD} for two possible definitions of those.

\medskip
The reader may have wondered why in
Conjectures~\ref{conj:Simion1} and \ref{conj:Simion2} we considered
determinants defined by type $B_n$ non-crossing partitions, which form
in fact a lattice, but used the extraneous type $A$ join\footnote{and
not even the one in the type $A$ {\it non-crossing} partition lattice!} in the
definition of the determinant, instead of the join which is
intrinsic to the lattice of type $B_n$ non-crossing partitions.
In particular, what would happen if we would make the latter choice?
As it turns out, for that situation there exists an elegant general
theorem due to Lindstr\"om 
\machSeite{LindAB}\cite{LindAB}, which I missed to state in 
\machSeite{KratBN}\cite{KratBN}. I refer to 
\machSeite{StanAP}\cite[Ch.~3]{StanAP} for the explanation of the
poset terminology used in the statement.

\begin{Theorem} \label{thm:Lind}
Let $L$ be a finite meet semilattice, $R$ be a commutative ring, and 
$f:L\times L\to R$ be an incidence function, that is, $f(x,y)=0$ unless
$x\land y=x$. Then
\begin{equation} \label{eq:Lind}
\det_{x,y\in L}\big(f(x\land y,x)\big)=
\prod _{y\in L} ^{}
\(\sum _{x\in L} ^{}\mu(x,y)f(x,y) \),
\end{equation}
where $\mu$ is the M\"obius function of $L$.\quad \quad \qed
\end{Theorem}

Clearly, this does indeed answer our question, we just have to
specialise $f(x,y)=h(x)$ for $x\land y=x$, where $h$ is
some function from $L$ to $R$. The fact that the above
theorem talks about meets instead of joins is of course
no problem because this is just a matter
of convention. 

Having an answer in such a great generality, one is tempted to pose
the problem of finding a general theorem that would encompass the
above-mentioned determinant evaluations due to Tutte, Dahab, Ko and
Smolinsky, Di~Francesco, as well as Conjectures~\ref{conj:Simion1} and
\ref{conj:Simion2}. This problem is essentially Problem~6 in 
\machSeite{SchFAA}\cite{SchFAA}.

\begin{Problem} \label{prob:9}
Let $L$ and $L'$ be two lattices {\em(}semilattices?{\em)} with
$L'\subseteq L$. Furthermore, let $R$ be a commutative ring, and 
let $f$ be a function from $L$ to $R$.
Under which conditions is there a compact formula for the determinant
\begin{equation} \label{eq:det9}
\det_{x,y\in L'}\big(f(x\land_L y)\big) ,
\end{equation}
where $\land_L$ is the meet operation in $L$?
\end{Problem}

By specialisation in Theorem~\ref{thm:Lind}, 
one can derive numerous corollaries. For example, a
very attractive one is the evaluation of the ``GCD determinant" due to Smith 
\machSeite{SmitAA}\cite{SmitAA}. (In fact, Smith's result is a more
general one for {\it factor closed subsets} of the positive integers.)

\begin{Theorem} \label{thm:Smith}
For any positive integer $n$, we have
$$\det_{1\le i,j\le n}\big(\gcd(i,j)\big)=
\prod _{i=1} ^{n}\phi(i),$$
where $\phi$ denotes the Euler totient function.\quad \quad \qed
\end{Theorem}

An interesting generalisation of Theorem~\ref{thm:Lind}
to posets was given by
Altini\c sik, Sagan and Tu\u glu 
\machSeite{AlSTAA}\cite{AlSTAA}. 
Again, all undefined terminology can be found in 
\machSeite{StanAP}\cite[Ch.~3]{StanAP}.

\begin{Theorem} \label{thm:AlSaTu}
Let $P$ be a finite poset, $R$ be a commutative ring, and $f,g:P\times
P\to R$ be two incidence functions, that is, $f(x,y)=0$ unless $x\le y$
in $P$, the same being true for $g$. Then
$$\det_{x,y\in P}\(
\sum _{z\in P} ^{}f(z,x)g(z,y)\)=
\prod _{x\in P} ^{}f(x,x)g(x,x).$$
\quad \quad \qed
\end{Theorem}

The reader is referred to Section~3 of
\machSeite{AlSTAA}\cite{AlSTAA} for the
explanation why this theorem implies Lindstr\"om's.

\subsection{Determinants for compositions}
Our next family of determinants consists of determinants of matrices
the rows and columns of which are indexed by {\it
compositions}. Recall that a composition of a non-negative integer $n$ 
is a vector $(\al_1,\al_2,\dots,\al_k)$ of non-negative integers such that
$\al_1+\al_2+\dots+\al_k=n$, for some $k$. For a fixed $k$, 
let $\mathcal C(n,k)$ denote the corresponding set of compositions of $n$.
While working on a problem in global optimisation, 
Brunat and Montes 
\machSeite{BrMoAA}\cite{BrMoAA}
discovered the following surprising determinant evaluation.
It allowed them to show how to explicitly express a multivariable
polynomial as a {\it difference of convex functions}.
In the statement, we use standard multi-index notation: 
if $\boldsymbol\al=(\al_1,\al_2,\dots,\al_k)$ and
$\boldsymbol \be=(\be_1,\be_2,\dots,\be_k)$ are two compositions, we
let
$$\boldsymbol\al^{\boldsymbol \be}:=\al_1^{\be_1}\al_2^{\be_2}\cdots
\al_k^{\be_k}.$$

\begin{Theorem} \label{thm:Brunat}
For any positive integers $n$ and $k$, we have
\begin{equation} \label{eq:Brunat} 
\det_{\boldsymbol\al,\boldsymbol\be\in\mathcal
C(n,k)}\left(\boldsymbol\al^{\boldsymbol \be}\right)=
n^{{\binom{ n + k-1} k}+k-1}\,
 \prod_{i = 1}^{ n-1}
      i^{( n -i+1) {\binom {n + k-i-1} { k-2}}}.
\end{equation}
\quad \quad \qed
\end{Theorem}

In recent joint work 
\machSeite{BrKMAA}\machSeite{BrMoAB}\cite{BrKMAA,BrMoAB}, 
Brunat, Montes and the author showed that there is in fact
a polynomial generalisation of this determinant evaluation.

\begin{Theorem} \label{conj:Brunat1}
Let $\mathbf x=(x_1,x_2,\dots,x_k)$ be a vector of indeterminates,
and let $\la$ be an indeterminate.
Then, for any non-negative integers $n$ and $k$, we have
\begin{equation} \label{eq:Brunat1} 
\det_{\boldsymbol\al,\boldsymbol\be\in\mathcal
C(n,k)}\left((\mathbf x+\la\boldsymbol\al)^
    {\boldsymbol \be}\right)=\la^{(k-1)\binom {n+k-1}k}
 \left({ \vert \mathbf x\vert+\la n }\right)
      ^{\binom{n+k-1} {k}}
\prod _{i=1} ^{n}i^{(k-1)\binom {n+k-i-1}{k-1}},
\end{equation}
where $\mathbf x+\la\boldsymbol\al$ is short for
$(x_1+\la\al_1,x_2+\la\al_2,\dots,x_k+\la\al_k)$, and
where $\vert \mathbf x\vert=x_1+x_2+\dots+x_k$. 
\end{Theorem}

As a matter of fact,  
there is actually a binomial variant which implies the
above theorem. Extending our
multi-index notation, let
$$\binom{\boldsymbol\al}{\boldsymbol \be}:=
\binom{\al_1}{\be_1}\binom{\al_2}{\be_2}\cdots
\binom{\al_k}{\be_k}.$$

\begin{Theorem} \label{conj:Brunat2}
Let $\mathbf x=(x_1,x_2,\dots,x_k)$ be a vector of indeterminates,
and let $\la$ be an indeterminate.
Then, using the notation from Theorem~{\em\ref{conj:Brunat1}},
for any non-negative integers $n$ and $k$, we have
\begin{equation} \label{eq:Brunat2} 
\det_{\boldsymbol\al,\boldsymbol\be\in\mathcal
C(n,k)}\left(\binom{\mathbf x+\la\boldsymbol\al}
    {\boldsymbol \be}\right)=\la^{(k-1)\binom {n+k-1}k}
\prod _{i=1} ^{n} \left(\frac{{ \vert \mathbf x\vert+(\la-1)n+i }}
 i\right)
      ^{\binom{n+k-i-1} {k-1}}.
\end{equation}
\end{Theorem}

The above theorem is proved in
\machSeite{BrKMAA}\cite{BrKMAA} by 
identification of factors (see ``Method~3" in Section~\ref{sec:eval}). 
Theorem~\ref{conj:Brunat2} follows by extracting the 
highest homogeneous component in \eqref{eq:Brunat1}. 
I report that, if one naively replaces ``compositions" by ``integer
partitions" in the above considerations, then the arising determinants
do not have nice product formulae.

\medskip
Another interesting determinant of a matrix with rows indexed by
compositions appears in the work of Bergeron, Reutenauer, Rosas and
Zabrocki 
\machSeite{BeRRAA}\cite[Theorem~4.8]{BeRRAA} on 
{\it Hopf algebras of non-commutative symmetric functions}. It
was used there to show that a certain set of generators of 
non-commutative symmetric functions were algebraically independent.
To state their determinant evaluation, we need to introduce some
notation. Given a composition 
$\boldsymbol\al=(\al_1,\al_2,\dots,\al_k)$ of $n$
with all summands $\al_i$ positive, 
we let
$D(\boldsymbol\al)=\{\al_1,\al_1+\al_2,\dots,\al_1+\al_2+\dots+\al_{k-1}\}$.
Furthermore, for two compositions $\boldsymbol\al$ and 
$\boldsymbol\be$ of $n$, we write
$\boldsymbol\al\cup\boldsymbol\be$ for the (unique) composition 
$\boldsymbol\ga$ of $n$ with
$D(\boldsymbol\ga)=D(\boldsymbol\al)\cup D(\boldsymbol\be)$. 
Finally, $\boldsymbol\al!$ is short for
$\al_1!\,\al_2!\cdots\al_k!$. 

\begin{Theorem} \label{thm:BRRZ}
Let $\Comp(n)$ denote the set of all compositions of $n$ all summands
of which are positive. Then
\begin{equation} \label{eq:BRRZ} 
\det_{\boldsymbol\al,\boldsymbol\be\in\Comp(n)}
\big((\boldsymbol\al\cup\boldsymbol\be)!\big)=
\prod _{\boldsymbol\ga\in\Comp(n)} ^{}
\prod _{i=1} ^{\ell(\boldsymbol\ga)}a_{\ga_{i}},
\end{equation}
where $\ell(\boldsymbol\ga)$ is the number of summands 
{\em(}components{\em)} of
the composition $\boldsymbol\ga$, and where $a_m$ denotes the number of
{\em indecomposable permutations} of $m$ {\em(}cf.\ 
\machSeite{StanBI}\cite[Ex.~5.13(b)]{StanBI}{\em)}. 
These numbers can be computed
recursively by $a_1=1$ and 
$$a_n=n!-\sum _{i=1} ^{n-1}a_i(n-1)!, \quad \quad n>1.$$
\quad \quad \qed
\end{Theorem}

As explained in
\machSeite{BeRRAA}\cite{BeRRAA}, one proves the theorem by factoring the matrix
in \eqref{eq:BRRZ} in
the form $CDC^t$, where $C$ is the ``incidence matrix" of ``refinement of
compositions," and where $D$ is a diagonal matrix. Thus, in
particular, the LU-factorisation of the matrix is determined.

\subsection{Two partition determinants}
On the surface,
{\it integer partitions} (see below for their definition) seem to be
very closely related to compositions, as they can be considered as
``compositions where the order of the summands is without importance."
However, experience shows that integer partitions are much more
complex combinatorial objects than compositions. This may be the
reason that the ``composition determinants" from the preceding
subsection do not seem to have analogues for integer partitions.
Leaving aside this disappointment, 
here {\it is} a determinant of a matrix 
in which rows and columns are indexed by integer
partitions. This determinant arose in work on {\it linear forms of values
of the Riemann zeta function evaluated at positive integers}, 
although the traces of it have now been completely erased in the
final version of the article 
\machSeite{KrRiAA}\cite{KrRiAA}. (The symmetric function
calculus in Section~12 of the earlier version 
\machSeite{KrRiAAA}\cite{KrRiAAA}
gives a vague idea where it may have come from.) 

Recall that the {\it power symmetric function} of degree $d$ in
$x_1,x_2,\dots,x_k$ is given by $x_1^d+x_2^d+\dots+x_k^d$, and is
denoted by $p_d(x_1,x_2,\dots,x_k)$. (See 
\machSeite{LascAZ}\cite[Ch.~1 and 2]{LascAZ},
\machSeite{MacdAC}\cite[Ch.~I]{MacdAC} and 
\machSeite{StanBI}\cite[Ch.~7]{StanBI} for 
in-depth expositions of the theory of 
{\it symmetric functions}.) Then, while working on 
\machSeite{KrRiAA}\cite{KrRiAA},
Rivoal and the author needed to evaluate the determinant
\begin{equation} \label{eq:Part} 
\det_{\la,\mu\in\Part(n,k)}\big(p_{\la}(\mu_1,\mu_2,\dots,\mu_k)\big),
\end{equation}
where $\Part(n,k)$ is the set of integer partitions of $n$ with at
most $k$ parts, that is, the set of all possibilities to write $n$ as
a sum of non-negative integers,
$n=\la_1+\la_2+\dots+\la_k$, with $\la_1\ge\la_2\ge\dots\ge\la_k\ge0$
(the non-zero $\la_i$'s being called the {\it parts of $\la$}),
and where 
$$p_\la(x_1,x_2,\dots,x_k)=p_{\la_1}(x_1,\dots,x_k)
p_{\la_1}(x_1,\dots,x_k)\cdots p_{\la_k}(x_1,\dots,x_k).
$$ 

Following the advice given in Section~\ref{sec:comb}, we went to the
computer and let it calculate the prime factorisations of the values
of this determinant for small values of $n$ and $k$. Indeed, the prime
factors turned out be always small so that we were sure that a
``nice" formula exists for the determinant. However, even more seemed
to be true. Recall that, in order to facilitate a proof of a (still
unknown) formula, it is (almost) always a good idea to try to
introduce more parameters (see 
\machSeite{KratBN}\cite[Sec.~2]{KratBN}). 
This is what we did. It led us consider the following determinant,
\begin{equation} \label{eq:Part-gen} 
\det_{\la,\mu\in\Part(n,k)}\big(p_{\la}(\mu_1+X_1,\mu_2+X_2,\dots,\mu_k+X_k)
\big),
\end{equation}
where $X_1,X_2,\dots,X_k$ are indeterminates.
Here are some values of the determinant \eqref{eq:Part-gen} for
special values of $n$ and $k$. For $n=k=3$ we obtain
$$ 6 \,( X_1 - X_2+2 )  ( X_1 - X_3+1 )  
   ( X_2 - X_3+1 )  (  X_1 + X_2 + X_3+3) ^4,
$$
for $n=4$ and $k=3$ we obtain
\begin{multline*}
 8 \,(  X_1 - X_2 +1)  (  X_1 - X_2+2 )  
   (  X_1 - X_2+3 )  (  X_1 - X_3+2 )  \\
\times
   (  X_2 - X_3+1 )  (  X_1 + X_2 + X_3+4) ^7, 
\end{multline*}
for $n=5$ and $k=3$ we get
\begin{multline*}
8 \,( X_1 - X_2+1) ( X_1 - X_2+2) ( X_1 - X_2+3)
    ( X_1 - X_2+4) ( X_1 - X_3+2)\\
\times ( X_1 - X_3+3) ( X_2 - X_3+1)
    ( X_2 - X_3+2) ( X_1 + X_2 + X_3+5)^{11},
\end{multline*}
for $n=6$ and $k=3$ we get
\begin{multline*}
 576 \,(  X_1 - X_2+1 )  (  X_1 - X_2+2 ) ^2 
   (  X_1 - X_2+3 ) ^2 (  X_1 - X_2+4 )  
   (  X_1 - X_2+5 )\\
\times  (  X_1 - X_3+1 )  
   (  X_1 - X_3+2 )  (  X_1 - X_3+3 )  
   (  X_1 - X_3+4 )  (  X_2 - X_3+1 ) ^2 \\
\times
   (  X_2 - X_3+2 ) ^2 
   (  X_1 + X_2 + X_3+6 ) ^{16},
\end{multline*}
for $n=4$ and $k=4$ we get
\begin{multline*}
 192\,(  X_1 - X_2 +1) (  X_1 - X_2+2 )
   (  X_1 - X_2+3 ) (  X_1 - X_3+2 )
   (  X_2 - X_3+1 ) \\
\times (  X_1 - X_4+1 )
   (  X_2 - X_4+1 ) (  X_3 - X_4+1 )
   (  X_1 + X_2 + X_3 + X_4+4 ) ^7
\end{multline*}
while for $n=5$ and $k=4$ we get
\begin{multline*}
 48 \,(  X_1 - X_2+1 )  (  X_1 - X_2+2 )  
   (  X_1 - X_2+3 )  (  X_1 - X_2+4 )  
   (  X_1 - X_3+2 ) \\
\times (  X_1 - X_3+3 )   (  X_2 - X_3+1 )  (  X_2 - X_3+2 )  
   (  X_1 - X_4+2 ) \\
\times (  X_2 - X_4+1)     (  X_3 - X_4+1 )  
   (  X_1 + X_2 + X_3 + X_4+5 ) ^{12}.
\end{multline*}
It is ``therefore" evident that there will be one factor which is a
power of $n+\sum _{i=1} ^{k}X_i$, whereas the other factors will be of
the form $X_i-X_j+c_{i,j}$, again raised to some power. 
In fact, the determinant is easy to compute for $k=1$ because, in that
case, it reduces to a special case of the
Vandermonde determinant evaluation. However, we were not able to come even
close to computing \eqref{eq:Part-gen} in general. 
(As I already indicated, finally we managed to avoid
the determinant evaluation in our work in 
\machSeite{KrRiAA}\cite{KrRiAA}.) 


To proceed further, we remark that evaluating \eqref{eq:Part-gen} is
equivalent to evaluating the same determinant, 
but with the power symmetric functions
$p_\la$ replaced by the {\it Schur functions} $s_\la$, because the
transition matrix between these two bases of symmetric functions is
the {\it character table of the symmetric group} of the corresponding
order (cf.\ 
\machSeite{MacdAC}\cite[Ch.~I, Sec.~7]{MacdAC}), the determinant of
which is known (see Theorem~\ref{thm:chi} below;
since in our determinants \eqref{eq:Part} and \eqref{eq:Part-gen} the
indices range over all partitions of $n$ {\it with at most $k$ parts}, it is
in fact the refinement given in Theorem~\ref{thm:chik} 
which we have to apply). The determinant with
Schur functions has some advantages over the one with power symmetric
functions since the former decomposes into a finer block structure. 
Alain Lascoux observed that, in fact, there is a
generalisation of the Schur function determinant to {\it 
Gra\3mannian Schubert polynomials}, which contains another set of variables,
$Y_1,Y_2,\dots,Y_{n+k-1}$. More precisely, given
$\la=(\la_1,\la_2,\dots,\la_k)$, let $\mathbb
Y_\la(X_1,X_2,\dots,X_k;Y_1,Y_2,\dots)$ denote the polynomial in the
variables $X_1,X_2,\dots,X_k$ and $Y_1,Y_2,\dots$, defined by (see 
\machSeite{LascAZ}\cite[Sections~1.4 and 9.7; the order of the
$B_{k_i}$ should be reversed in (9.7.2) and analogous places]{LascAZ})
$$\mathbb
Y_\la(X_1,X_2,\dots,X_k;Y_1,Y_2,\dots):=\det_{1\le i,j\le k}(S_{\la_i-i+j}
(X_1,\dots,X_k;Y_1,\dots,Y_{\la_{i}+k-i})),$$
where the entries of the determinant are given by 
\begin{equation} \label{eq:XY} 
\sum _{m=0} ^{\infty}S_m
(X_1,\dots,X_k;Y_1,\dots,Y_l)x^m=
\frac {\prod _{i=1} ^{l}(1-Y_ix)} {\prod _{i=1} ^{k}(1-X_ix)}.
\end{equation}
The Gra\3mannian Schubert polynomial $\mathbb Y_\la$ reduces to the Schur
function $s_\la$ when all the variables $Y_i$, $i=1,2,\dots$, are set equal
to $0$. Given these definitions, Alain Lascoux (private communication)
established the following result. 

\begin{Theorem} \label{thm:Lascoux}
Let $X_1,X_2,\dots,X_k,Y_1,Y_2,\dots,Y_{n+k-1}$ be indeterminates.
Then,
\begin{multline} \label{eq:Part-gena} 
\det_{\la,\mu\in\Part(n,k)}
\big(\mathbb Y_{\la}(\mu_1+X_1,\mu_2+X_2,\dots,\mu_k+X_k;
Y_1,Y_2,\dots)
\big)\\=\Pi(n,k)
\prod _{\si\in G(n,k)} ^{}\(n+\sum _{i=1} ^{k}X_i-\sum _{i=1}
^{k}Y_{\si(i)}\),
\end{multline}
where $G(n,k)$ denotes the set of
all {\it Gra\3mannian permutations},\footnote{A permutation $\si$ in
$\mathfrak S_\infty$ (the set of all permutations of the natural
numbers $\mathbb N$ which fix all but a finite number of elements of 
$\mathbb N$)
is called {\it Gra\3mannian} if $\si(i)<\si(i+1)$ for all $i$ except 
possibly for one, the latter being called the {\it descent} of $\si$
(see \machSeite{MacdAE}\cite[p.~13]{MacdAE}). An inversion of $\si$ is
a pair $(i,j)$ such that $i<j$ but $\si(i)>\si(j)$. We remark that the
number of Gra\3mannian permutations with descent (if existent) at $k$
and at most $n-1$ inversions is equal to the number of
partitions of at most $n-1$ 
with at most $k$ parts (including the empty partition). 
More concretely, if we denote this number by $g_{n,k}$, the
generating function of the numbers $g_{n,k}$ is given by
$$
\sum _{n=1} ^{\infty}g_{n,k}x^{n-1}=
\frac {1} {1-q}\prod _{i=1} ^{k}\frac {1} {1-q^i}.$$
} 
the descent of which {\em(}if existent{\em)} is at $k$, and
which contain at most $n-1$ inversions, and where $\Pi(n,k)$ is given
recursively by
\begin{equation} \label{eq:Pi} 
\Pi(n,k)=\Pi(n,k-1)\Pi(n-k,k)
\prod _{\mu\in\Part^+(n,k)} ^{}
\prod _{j=1} ^{k}(X_j+\mu_j-X_k),
\end{equation}
$\Part^+(n,k)$ denoting the set of partitions of $n$ into {\em exactly}
$k$ {\em(}positive{\em)} parts,
with initial conditions $\Pi(n,k)=1$ if $k=1$ or $n\le1$. Explicitly,
\begin{equation} \label{eq:Piexpl} 
\Pi(n,k)=\underset{\vert w\vert_B<k-1}
{\prod _{w} ^{}}
\prod _{\mu\in\Part(n-\vert w\vert_A k+\inv w,k-\vert w\vert_B)} ^{}
\prod _{j=1} ^{k-\vert w\vert_B}(X_j+\mu_j-X_{k-\vert w\vert_B}),
\end{equation}
where the product over $w$ runs over all finite-length words $w$ with letters
from $\{A,B\}$, including the empty word. 
The notation $\vert w\vert_B$ means the number of
occurrences of $B$ in $w$, with the analogous meaning for $\vert
w\vert_A$. The quantity $\inv w$ denotes the number of inversions of
$w=w_1w_2\dots$, 
which is the number of pairs of letters $(w_i,w_j)$, $i<j$, such
that $w_i=B$ and $w_j=A$.
\end{Theorem}

We obtain the evaluation of the determinant \eqref{eq:Part} if we set
$X_i=Y_i=0$ for all $i$ in the above theorem and multiply by \eqref{eq:chik}.
Likewise, we obtain the evaluation of the determinant
\eqref{eq:Part-gen} if we set $Y_i=0$ for all $i$ in the 
above theorem and multiply the result by \eqref{eq:chik}.

For example, here is the determinant \eqref{eq:Part-gena} for $n=3$ and $k=2$, 
\begin{multline*}
(X_1-X_2+2)(X_1+X_2-Y_1-Y_2+3)(X_1+X_2-Y_1-Y_3+3)\\
\times
(X_1+X_2-Y_2-Y_3+3)(X_1+X_2-Y_1-Y_4+3),
\end{multline*}
and the following is the one for $n=k=3$,
\begin{multline*}
(X_1-X_2+2)(X_1-X_3+1)(X_2-X_3+1)\\
\times
(X_1+X_2+X_3-Y_1-Y_2-Y_3+3)
(X_1+X_2+X_3-Y_1-Y_2-Y_4+3)\\
\times
(X_1+X_2+X_3-Y_1-Y_3-Y_4+3)
(X_1+X_2+X_3-Y_1-Y_2-Y_5+3).
\end{multline*}

In the sequel, I sketch Lascoux's proof of Theorem~\ref{thm:Lascoux}.
The first step consists in applying the 
{\it Monk formula for double Schubert polynomials}
in the case of Gra\3mannian Schubert polynomials (see
\machSeite{KoVeAA}\cite{KoVeAA}), 
 $$ \(\sum _{i=1} ^{k}X_i-\sum _{i=1} ^{k}Y_{\la_{k-i+1}+i}\)  \mathbb Y_\la 
= \mathbb Y_{\la+(1,0,0,\dots)} + 
\sum_\mu \mathbb Y_\mu, 
$$
where $\mathbb Y_\la$ is short for $\mathbb
Y_\la(X_1,X_2,\dots,X_k;Y_1,Y_2,\dots)$, 
and where the sum on the right-hand side 
is over all partitions $\mu$ of the same size as
$\la+(1,0,0,\dots)$ but lexicographically smaller. Clearly, by using
this identity, appropriate row operations
in the determinant \eqref{eq:Part-gena} show that 
\begin{equation} \label{eq:factor} 
n+\sum _{i=1} ^{k}X_i-\sum _{i=1} ^{k}Y_{\si(i)},
\end{equation}
is one of its factors, the relation between $\si$ and $\la$
being $\si(i)=\la_{k-i+1}+i$, $i=1,2,\dots,k$. 
In particular, $\si$ can be extended (in a unique way) 
to a Gra\3mannian permutation.
Moreover, doing these row operations, and taking out the
factors of the form \eqref{eq:factor}, we collect on the one hand the
product in \eqref{eq:Part-gena}, and
we may on the other hand reduce
the determinant \eqref{eq:Part-gena} to a determinant of the same
form, in which, however, the partitions
$\la=(\la_1,\la_2,\dots,\la_k)$ run over all partitions
of size {\it at most\/} $n$
with the {\it additional property that\/} $\la_1=\la_2$ (instead of over all
partitions from
$\Part(n,k)$). 

As it turns out, the determinant 
thus obtained is independent of the variables $Y_1,Y_2,\break
\dots,Y_{n+k-1}$.
Indeed, if we expand each Schubert polynomial
$\mathbb Y_\la(\mu_1+X_1,\mu_2+X_2,\dots,\mu_k+X_k;Y_1,Y_2,\dots)$ 
in the determinant as a
linear combination of Schur functions in $X_1,X_2,\break
\dots,X_k$ with
coefficients being polynomials in the $Y_1,Y_2,\dots,Y_{n+k-1}$, then, by also
using that 
$$S_1(\mu_1+X_1,\mu_2+X_2,\dots,\mu_k+X_k)=
n+\sum _{i=1} ^{k}X_i$$
is independent of $\mu$, it is not difficult to see that one can
eliminate all the $Y_i$'s by appropriate row operations.

To summarise the current state of the discussion: we have already
explained the occurrence of the product on the right-hand side of 
\eqref{eq:Part-gena} as a factor of the determinant. Moreover, the
remaining factor is given by a determinant of the same form as in
\eqref{eq:Part-gena}, 
\begin{equation} \label{eq:det-Y} 
\det_{\la,\mu}
\big(\mathbb Y_{\la}(\mu_1+X_1,\mu_2+X_2,\dots,\mu_k+X_k;
Y_1,Y_2,\dots)
\big),
\end{equation}
where, as before, $\mu$ runs over
all partitions in $\Part(n,k)$, but where $\la$ runs over partitions
$(\la_1,\la_2,\dots,\la_k)$
of size at most $n$, with the additional
restriction that $\la_1=\la_2$. The final observation was that this
latter determinant is in fact independent of the $Y_i$'s. This allows
us to specify $Y_1$ arbitrarily, say $Y_1=X_k$. Now another property
of double Schubert polynomials, namely that (this follows from the
definition of double Schubert polynomials by means of divided
differences, see \machSeite{LascAZ}\cite[(10.2.3)]{LascAZ}, 
and standard properties of divided differences)
\begin{equation} \label{eq:Y1} 
\mathbb Y_{\rho+(1,1,\dots,1)}(X_1,X_2,\dots,X_k;Y_1,Y_2,\dots)=
\mathbb Y_{\rho}(X_1,X_2,\dots,X_k;Y_2,\dots)\prod _{j=1} ^{k}(X_j-Y_1)
\end{equation}
comes in handy. (In the vector $(1,1,\dots,1)$ on the left-hand side
of \eqref{eq:Y1} there are $k$ occurrences of $1$.) 
Namely, if $Y_1=X_k$, the matrix in \eqref{eq:det-Y}
(of which the determinant is taken), 
$M(n,k;Y_1,Y_2,\dots)$ say,
decomposes in block form. If $\la$ is a partition with $\la_k>0$ and
$\mu$ is a partition with $\mu_k=0$, then, because of \eqref{eq:Y1},
the corresponding entry in \eqref{eq:det-Y} vanishes. Furthermore, in
the block where $\la$ and $\mu$ are partitions with $\la_k=\mu_k=0$,
because of the definition \eqref{eq:XY} of the quantities $S_m(.)$
the corresponding entry reduces to
\begin{multline*}
\mathbb Y_{\la}(\mu_1+X_1,\mu_2+X_2,\dots,\mu_{k-1}+X_{k-1},X_k;
X_k,Y_2,\dots)\\
=\mathbb Y_{\la}(\mu_1+X_1,\mu_2+X_2,\dots,\mu_{k-1}+X_{k-1};
Y_2,\dots).
\end{multline*}
In other words, this block is identical with $M(n,k-1;Y_2,\dots)$. 
Finally, in
the block indexed by partitions $\la$ and $\mu$ with $\la_k>0$ and
$\mu_k>0$, we may use \eqref{eq:Y1} to factor $
\prod _{j=1} ^{k}(X_j+\mu_j-X_k)$ out of the column indexed by
$\mu$, for all such $\mu$. 
What remains is identical with $M(n-k,k;Y_2,\dots)$. Taking determinants,
we obtain the recurrence \eqref{eq:Pi}. (Here we use again that the
determinants in \eqref{eq:det-Y}, that is, 
in particular, the determinants of
$M(n,k;Y_1,\dots)$, $M(n,k-1;Y_2,\dots)$, and of $M(n-k,k;Y_2,\dots)$, are
all independent of the $Y_i$'s.) The explicit form
\eqref{eq:Piexpl} for $\Pi(n,k)$ can be easily derived by induction on
$n$ and $k$.

\medskip
A few paragraphs above, we mentioned in passing another interesting determinant
of a matrix the rows and columns of which are indexed by (integer) partitions:
the {\it determinant of the character table of the symmetric group $\mathfrak
S_n$} (cf.\ 
\machSeite{JameAA}\cite[Cor.~6.5]{JameAA}). 
Since this is a classical and beautiful determinant evaluation which I missed
to state in \machSeite{KratBN}\cite{KratBN}, I present it now in the
theorem below. There, the notation $\la\vdash n$ stands for ``$\la$ is
a partition of $n$."
For all undefined notation, I refer the reader to standard
texts on the representation theory of symmetric groups, as for example
\machSeite{JameAA}%
\machSeite{JaKeAA}%
\machSeite{SagaAQ}%
\cite{JameAA,JaKeAA,SagaAQ}.

\begin{Theorem} \label{thm:chi}
For partitions $\la$ and $\rh$ of $n$, let $\chi^\la(\rho)$ denote the
value of the irreducible character $\chi^\la$ evaluated at a permutation of
cycle type $\rh$. Then
\begin{equation} \label{eq:chi}
\det_{\la,\rh\,\vdash n}\(\chi^\la(\rho)\)=
\prod _{\mu\,\vdash n} ^{}
\prod _{i\ge1} ^{}\mu_i. 
\end{equation}
In words: the determinant of the character table of the symmetric
group $\mathfrak S_n$ is equal to the products of all the parts of all
the partitions of $n$.\quad \quad \qed
\end{Theorem}

A refinement of this statement, where we restrict to partitions 
of $n$ with at most $k$ parts, is the following.

\begin{Theorem} \label{thm:chik}
With the notation of Theorem~{\em\ref{thm:chi}}, for all positive
integers $n$ and $k$, $n\ge k$, we have
\begin{equation} \label{eq:chik}
\det_{\la,\rh\in\Part(n,k)}\(\chi^\la(\rho)\)=
\prod _{\mu\in\Part(n,k)} ^{}
\prod _{i\ge1} ^{}m_i(\mu)!, 
\end{equation}
where $m_i(\mu)$ is the number of times $i$ occurs as a part
in the partition $\mu$.\quad \quad \qed
\end{Theorem}

This determinant evaluation follows from the decomposition of the {\it
full\/} character table of $\mathfrak S_n$ in the form
\begin{equation} \label{eq:LK} 
\(\chi^\la(\rho)\)_{\la,\rh\,\vdash n}=L\cdot K^{-1},
\end{equation} 
where $L=(L_{\la,\mu})_{\la,\mu\,\vdash n}$ 
is the transition matrix from power symmetric functions
to monomial symmetric functions, and where 
$K=(K_{\la,\mu})_{\la,\mu\,\vdash n}$ is the {\it Kostka matrix},
the transition matrix from Schur functions to monomial symmetric
functions (see 
\machSeite{MacdAC}\cite[Ch.~I, (6.12)]{MacdAC}). For, if we order the
partitions of $n$ so that the partitions in $\Part(n,k)$ come before
the partitions in $\Part(n,k+1)$, $k=1,2,\dots,n-1$, and within
$\Part(n,k)$ lexicographically, then, with respect to this order, the
matrix $L$ is lower triangular and the matrix $K$, and hence also
$K^{-1}$, is upper triangular. Furthermore,
the matrix $K$, and hence also $K^{-1}$, is even block diagonal,
the blocks along the diagonal being the ones which are formed by the
rows and columns indexed by the partitions in $\Part(n,k)$,
$k=1,2,\dots,n$. These facts together imply that the decomposition
\eqref{eq:LK} restricts to the submatrices indexed by partitions in
$\Part(n,k)$, 
\begin{equation*}
\(\chi^\la(\rho)\)_{\la,\rh\in \Part(n,k)}=
(L_{\la,\mu})_{\la,\mu\in\Part(n,k)}\cdot 
(K_{\la,\mu})_{\la,\mu\in\Part(n,k)}^{-1}.
\end{equation*} 
If we now take determinants on both sides, then, in view of
$K_{\mu,\mu}=1$ and of
$$L_{\mu,\mu}=\prod _{i\ge1} ^{}m_i(\mu)!$$
for all $\mu$, the theorem follows.

\medskip
Further examples of nice
determinant evaluations of tables of
{\it characters of representations
of symmetric groups and their double covers} can be found in
\machSeite{BeOSAA}%
\machSeite{OlssAB}%
\cite{BeOSAA,OlssAB}.
Determinants of tables of {\it characters of the alternating group}
can be found in
\machSeite{BeOlAE}\cite{BeOlAE}.

\subsection{Elliptic determinant evaluations}\label{sec:ell}
In special functions theory there is currently a disease rapidly
spreading, generalising the earlier mentioned $q$-disease (see
Footnote~\ref{foot:q}). It could be called the {\it ``elliptic
disease}." Recall that, during the $q$-disease, we replaced every
positive integer $n$ by $1+q+q^2+\dots+q^{n-1}=(1-q^n) /(1-q)$,
and, more generally, shifted factorials $a(a+1)\cdots(a+k-1)$ by
$q$-shifted factorials $(1-\al)(1-q\al)\cdots(1-\al q^{k-1})$. (Here,
$\al$ takes the role of $q^a$, and one drops the powers of $1-q$ in
order to ease notation.) Doing this with some ``ordinary" identity, we
arrived (hopefully) at its {\it $q$-analogue}. 
Now, once infected by the elliptic disease, we would replace every
occurrence of a term $1-x$ (and, looking at the definition of
$q$-shifted factorials, we can see that there will be many) by its
{\it elliptic analogue} $\theta(x;p)$:
$$\theta(x)=\theta(x;p)=\prod_{j=0}^\infty(1-p^jx)(1-p^{j+1}/x). $$
Here, $p$ is a complex number with $|p|<1$, which will be fixed
throughout. Up to a trivial factor,
$\theta(e^{2\pi ix};e^{2\pi i \tau})$ equals
the {\it Jacobi theta function} $\theta_1(x|\tau)$ (cf.\ 
\machSeite{WW}\cite{WW}). Clearly, $\theta(x)$ reduces to $1-x$ if $p=0$.

At first sight, one will be sceptical if this is a fruitful
thing to do. After all, for working with the functions $\theta(x)$,
the only identities which are available are
the (trivial) {\it inversion formula}
\begin{equation}\label{ti}\theta(1/x)=-\frac1 x\,\theta(x), \end{equation}
the (trivial) {\it quasi-periodicity}
\begin{equation}\label{tp}\theta(px)=-\frac1 x\,\theta(x),\end{equation}
and {\it Riemann's} (highly non-trivial) {\it addition formula} (cf.\ 
\machSeite{WW}\cite[p.~451, Example~5]{WW})
\begin{equation}\label{tadd}
\theta(xy)\,\theta(x/y)\,\theta(uv)\,\theta(u/v)-
\theta(xv)\,\theta(x/v)\,\theta(uy)\,\theta(u/y)
=\frac uy\,\theta(yv)\,\theta(y/v)\,\theta(xu)\,\theta(x/u).
\end{equation}
Nevertheless, it has turned out recently 
that a surprising number of identities from
the ``ordinary" and from the ``$q$-world" can be lifted to the
elliptic level. This is particularly true for series of hypergeometric
nature. We refer the reader to Chapter~11 of
\machSeite{GaRaAA}\cite{GaRaAA} for an account of the current state of
the art
in the theory of, as they are called now, {\it elliptic hypergeometric
series}. 

On the following pages,
I give elliptic determinant evaluations a rather extensive coverage
because, first of all, they were non-existent in
\machSeite{KratBN}\cite{KratBN} (with the exception of the mention of
the papers
\machSeite{MilnAO}%
\machSeite{MilnAP}\cite{MilnAP,MilnAO} by Milne), and, second, because
I believe that the ``elliptic research" is a research direction that
will further prosper in the next future and will have 
numerous applications in
many fields, also outside of just special functions theory and number
theory. I further believe that the determinant evaluations presented
in this subsection will turn out to be as fundamental as the
determinant evaluations in Sections~2.1 and 2.2 in
\machSeite{KratBN}\cite{KratBN}. For some of them this belief is
already a fact. For example,
determinant evaluations involving elliptic functions have
come into the picture in the theory of {\it
multiple elliptic hypergeometric series}, see
\machSeite{KN}%
\machSeite{RainAA}%
\machSeite{RoseAA}%
\machSeite{Ro}%
\machSeite{RoScAB}%
\machSeite{Sp}%
\machSeite{WarnAG}%
\cite{KN,RainAA,RoseAA,Ro,RoScAB,Sp,WarnAG}.
They have also an important role in the study of
\emph{Ruijsenaars operators} and related {\it integrable systems} 
\machSeite{H}%
\machSeite{Ru}%
\cite{H,Ru}. Furthermore, they have recently found applications
in number theory 
to the problem of counting the number of {\it representations of an integer
as a sum of triangular numbers} 
\machSeite{RoseAB}\cite{RoseAB}.

Probably the first elliptic determinant evaluation is due to Frobenius
\machSeite{Fr}\cite[(12)]{Fr}.
This identity  has found applications to {\it Ruijsenaars
operators}
\machSeite{Ru}\cite{Ru}, to {\it multidimensional elliptic hypergeometric
series} and {\it integrals} 
\machSeite{KN}\cite{KN}, 
\machSeite{RainAA}\cite{RainAA} and to {\it number theory} 
\machSeite{RoseAB}\cite{RoseAB}. For a generalisation
to {\it higher genus Riemann surfaces}, see 
\machSeite{F}\cite[Cor.~2.19]{F}.
Amdeberhan \machSeite{AmdeAC}\cite{AmdeAC} observed that it can be
easily proved using the condensation method (see ``Method~2" in
Section~\ref{sec:eval}). 

\begin{Theorem}
\label{froa}
Let $x_1,x_2,\dots,x_n$, $a_1,a_2,\dots,a_n$ and $t$ be indeterminates.
Then there holds
\begin{equation}\label{froaid}
\det_{1\leq i,j\leq n}\left(\frac{\theta(ta_jx_i)}{\theta(t)\,
\theta(a_jx_i)}\right)
=\frac{\theta(ta_1\dotsm a_nx_1\dotsm x_n)}{\theta(t)}
\frac{\displaystyle\prod_{1\leq i<j\leq n}
a_jx_j\,\theta(a_i/a_j)\,\theta(x_i/x_j)}
{\displaystyle\prod_{i,j=1}^n\theta(a_jx_i)}.
\end{equation}
\quad \quad \qed
\end{Theorem}

For $p=0$ and $t\to\infty$, this determinant identity reduces to 
Cauchy's evaluation \eqref{eq:Cauchy} of the double alternant, and,
thus, may be regarded as an ``elliptic analogue" of the latter.

Okada 
\machSeite{OkadAK}\cite[Theorem~1.1]{OkadAK} has recently found an
elliptic extension of Schur's evaluation
\eqref{eq:Schur} of a Cauchy-type Pfaffian.
His proof works by the Pfaffian version of the condensation method.

\begin{Theorem} \label{thm:Okada}
Let $x_1,x_2,\dots,x_n$, $t$ and $w$ be indeterminates.
Then there holds
\begin{multline}\label{eq:Okada}
\underset{1\leq i,j\leq 2n}\Pf\left(\frac{\theta(x_j/x_i)}{\theta(x_ix_j)}
\frac{\theta(tx_ix_j)}{\theta(t)}
\frac{\theta(wx_ix_j)}{\theta(w)}
\right)
\\=
\frac{\theta(tx_1\dotsm x_{2n})}{\theta(t)}
\frac{\theta(wx_1\dotsm x_{2n})}{\theta(w)}
\prod _{1\le i<j\le 2n} ^{}
\frac{\theta(x_j/x_i)}
{x_j\,\theta(x_ix_j)}.
\end{multline}
\quad \quad \qed
\end{Theorem}

The next group of determinant evaluations is from 
\machSeite{RoScAC}\cite[Sec.~3]{RoScAC}. As the Vandermonde determinant
evaluation, or the other Weyl denominator formulae (cf.\
\machSeite{KratBN}\cite[Lemma~2]{KratBN}), are fundamental {\it
polynomial} determinant evaluations, the evaluations in
Lemma~\ref{wp} below are equally
fundamental in the elliptic domain as they can be considered as the
elliptic analogues of the former. Indeed, Rosengren and Schlosser
show that they {\it imply} the {\it Macdonald identities associated to
affine root systems} 
\machSeite{MacdAA}\cite{MacdAA}, which are the affine analogues of the Weyl
denominator formulae. In particular, in this way they obtain new
proofs of the Macdonald identities.

In order to conveniently formulate Rosengren and
Schlosser's determinant evaluations, we shall adopt the following
terminology from
\machSeite{RoScAC}\cite{RoScAC}.
For $0<|p|<1$ and $t\neq 0$, 
an {\it $A_{n-1}$ theta function $f$ of norm $t$} is
a holomorphic function for $x\neq 0$ such that
\begin{equation}\label{ade}f(px)=\frac{(-1)^n}{tx^n}\,f(x).\end{equation}
Moreover, if $R$ denotes either of the root systems $B_n$,
$B^\vee_n$, $C_n$, $C^\vee_n$, $BC_n$ or $D_n$ (see
Footnote~\ref{foot:root} and \machSeite{HumpAC}\cite{HumpAC} for
information on root systems), 
  we call $f$ an {\it$R$ theta function} if
{\allowdisplaybreaks
\begin{align*}
f(px)&=-\frac{1}{p^{n-1}x^{2n-1}}\,f(x),& f(1/x)&=-\frac
1x\,f(x),& R&=B_n,\\
f(px)&=-\frac{1}{p^{n}x^{2n}}\,f(x),& f(1/x)&=-f(x),&
R&=B_n^\vee,\\
f(px)&=\frac{1}{p^{n+1}x^{2n+2}}\,f(x),& f(1/x)&=-f(x),&
R&=C_n,\\
f(px)&=\frac{1}{p^{n-\frac12}x^{2n}}\,f(x),& f(1/x)&=-\frac
1x\,f(x),& R&=C_n^\vee,\\
f(px)&=\frac{1}{p^{n}x^{2n+1}}\,f(x),& f(1/x)&=-\frac
1x\,f(x),& R&=BC_n,\\
f(px)&=\frac{1}{p^{n-1}x^{2n-2}}\, f(x),& f(1/x)&=f(x),&
R&=D_n.
\end{align*}
}

Given this definition,
Rosengren and Schlosser \machSeite{RoScAC}\cite[Lemma~3.2]{RoScAC}
show that a function $f$ is an $A_{n-1}$ theta function of norm $t$ 
if and only if there exist constants $C$, $b_1,\dots,b_{n}$ such that
$b_1\dotsm b_n=t$ and
$$f(x)=C\,\theta(b_1x)\cdots\theta(b_nx), $$
and for the other six cases, they show that $f$ is an 
$R$ theta function if and only if there exist
constants $C$, $b_1,\dots,b_{n-1}$ such that
\begin{align*}f(x)&=C\,\theta(x)\,\theta(b_1x)\,\theta(b_1/x)
\cdots\theta(b_{n-1}x)\,\theta(b_{n-1}/x),&
  R&=B_n,\\
f(x)&=C\,x^{-1}\theta(x^2;p^2)\,\theta(b_1x)\,\theta(b_1/x)
\cdots\theta(b_{n-1}x)\,\theta(b_{n-1}/x),&
R&=B_n^\vee,\\
f(x)&=C\,x^{-1}\theta(x^2)\,\theta(b_1x)\,\theta(b_1/x)
\cdots\theta(b_{n-1}x)\,\theta(b_{n-1}/x),& R&=C_n,\\
f(x)&=C\,\theta(x;p^{\frac12})\,\theta(b_1x)\,\theta(b_1/x)
\cdots\theta(b_{n-1}x)\,\theta(b_{n-1}/x),&
R&=C_n^\vee,\\
f(x)&=C\,\theta(x)\,\theta(px^2;p^2)\,\theta(b_1x)\,\theta(b_1/x)
\cdots\theta(b_{n-1}x)\,\theta(b_{n-1}/x),&
R&=BC_n,\\
f(x)&=C\,\theta(b_1x)\,\theta(b_1/x)
\cdots\theta(b_{n-1}x)\,\theta(b_{n-1}/x),& R&=D_n,
\end{align*}
where $\theta(x)=\theta(x;p)$.

If one puts $p=0$, then an $A_{n-1}$ theta function
of norm $t$ becomes a polynomial of degree $n$ such that the 
reciprocal of the product of its roots is equal to $t$. Similarly,
if one puts $p=0$, then a $D_n$ theta function 
becomes a polynomial in $(x+1/x)$ of degree $n$.
This is the
specialisation of some of the following results which is relevant for
obtaining the earlier Lemmas~\ref{lem:RS1}--\ref{cdetr1cor}.

The elliptic extension of the Weyl denominator formulae
is the following formula.
(See 
\machSeite{RoScAC}\cite[Prop.~3.4]{RoScAC}.)

\begin{Lemma}\label{wp}
Let $f_1,\dots,f_n$ be  $A_{n-1}$ theta functions of norm
$t$. Then, 
\begin{equation}\label{awpi}\det_{1\leq i,j\leq n}\left(f_j(x_i)\right)=
C\,W_{A_{n-1}}(x),
\end{equation}
for some constant $C$, where
$$W_{A_{n-1}}(x)=\theta(tx_1\dotsm x_n)
\,\prod_{1\leq i<j\leq n}x_j\theta(x_i/x_j).$$ 
Moreover, if $R$ denotes either  $B_n$,
$B^\vee_n$, $C_n$, $C^\vee_n$, $BC_n$ or $D_n$ and   $f_1,\dots,f_n$
are $R$ theta functions, we have
\begin{equation}\label{wpi}
\det_{1\leq i,j\leq n}\left(f_j(x_i)\right)=C\,W_R(x) ,
\end{equation}
for some constant $C$, where
{\allowdisplaybreaks
\begin{align*}
W_{B_n}(x)&=\prod_{i=1}^n\theta(x_i)
\prod_{1\leq i<j\leq n}x_i^{-1}\theta(x_ix_j)\,\theta(x_i/x_j), \\
W_{B_n^\vee}(x)&=\prod_{i=1}^n x_i^{-1}\theta(x_i^2;p^2)
\prod_{1\leq i<j\leq n}x_i^{-1}\theta(x_ix_j)\,\theta(x_i/x_j), \\
W_{C_n}(x)&=\prod_{i=1}^n x_i^{-1}\theta(x_i^2)
\prod_{1\leq i<j\leq n}x_i^{-1}\theta(x_ix_j)\,\theta(x_i/x_j),\\
W_{C_n^\vee}(x)&=\prod_{i=1}^n\theta(x_i;p^{\frac 12})
\prod_{1\leq i<j\leq n}x_i^{-1}\theta(x_ix_j)\,\theta(x_i/x_j), \\
W_{BC_n}(x)&=\prod_{i=1}^n\theta(x_i)\,\theta(px_i^2;p^2)
\prod_{1\leq i<j\leq n}x_i^{-1}\theta(x_ix_j)\,\theta(x_i/x_j), \\
W_{D_n}(x)&=\prod_{1\leq i<j\leq n}x_i^{-1}\theta(x_ix_j)\,\theta(x_i/x_j).
\end{align*}}%
\quad \quad \qed
\end{Lemma}

Rosengren and Schlosser show in
\machSeite{RoScAC}\cite[Prop.~6.1]{RoScAC} 
that the famous Macdonald identities for
affine root systems 
\machSeite{MacdAA}\cite{MacdAA} are equivalent to special cases of 
this lemma. We state the corresponding results below.

\begin{Theorem}\label{mdp}
The following determinant evaluations hold:
$$
\det_{1\leq i,j\leq n}\left(x_i^{j-1}
\theta((-1)^{n-1}p^{j-1}tx_i^n;p^n)\right)=\frac{(p;p)_\infty^{n}}{(p^n;p^n)_\infty^n}
\,W_{A_{n-1}}(x),
$$
\begin{multline*}
\det_{1\leq i,j\leq n}\left(x_i^{j-n}
\theta(p^{j-1}x_i^{2n-1};p^{2n-1})
-x_i^{n+1-j}
\theta(p^{j-1}x_i^{1-2n};p^{2n-1})
\right)\\
=\frac{2(p;p)_\infty^{n}}
{(p^{2n-1};p^{2n-1})_\infty^n}\,W_{B_{n}}(x),
\end{multline*}
\begin{multline*}
\det_{1\leq i,j\leq n}\left(x_i^{j-n-1}
\theta(p^{j-1}x_i^{2n};p^{2n})
-x_i^{n+1-j}
\theta(p^{j-1}x_i^{-2n};p^{2n})
\right)\\
=\frac{2(p^2;p^2)_\infty(p;p)_\infty^{n-1}}
{(p^{2n};p^{2n})_\infty^n}\,W_{B_{n}^\vee}(x),
\end{multline*}
\begin{multline*}
\det_{1\leq i,j\leq n}\left(x_i^{j-n-1}
\theta(-p^{j}x_i^{2n+2};p^{2n+2})
-x_i^{n+1-j}
\theta(-p^{j}x_i^{-2n-2};p^{2n+2})
\right)\\
=\frac{(p;p)_\infty^{n}}
{(p^{2n+2};p^{2n+2})_\infty^n}\,
W_{C_{n}}(x),
\end{multline*}
\begin{multline*}
\det_{1\leq i,j\leq n}\left(x_i^{j-n}
\theta(-p^{j-\frac12}x_i^{2n};p^{2n})
-x_i^{n+1-j}
\theta(-p^{j-\frac12}x_i^{-2n};p^{2n})
\right)\\
=\frac{(p^{\frac12};p^{\frac12})_\infty(p;p)_\infty^{n-1}}
{(p^{2n};p^{2n})_\infty^n}\,  
W_{C_{n}^\vee}(x),
\end{multline*}
\begin{multline*}
\det_{1\leq i,j\leq n}\left(x_i^{j-n}
\theta(-p^{j}x_i^{2n+1};p^{2n+1})
-x_i^{n+1-j}
\theta(-p^{j}x_i^{-2n-1};p^{2n+1})
\right)\\
=\frac{(p;p)_\infty^n}{(p^{2n+1};p^{2n+1})_\infty^n}\,
W_{BC_{n}}(x),
\end{multline*}
\begin{multline*}
\det_{1\leq i,j\leq n}\left(x_i^{j-n}
\theta(-p^{j-1}x_i^{2n-2};p^{2n-2})
+x_i^{n-j}
\theta(-p^{j-1}x_i^{2-2n};p^{2n-2})
\right)\\
=\frac{4(p;p)_\infty^n}{(p^{2n-2};p^{2n-2})_\infty^n}
\,W_{D_{n}}(x),
\qquad n\geq 2.
\end{multline*}
\quad \quad \qed
\end{Theorem}

Historically, aside from Frobenius' elliptic Cauchy identity \eqref{froaid}, 
the subject of elliptic determinant evaluations begins with Warnaar's
remarkable paper \machSeite{WarnAG}\cite{WarnAG}. While the main subject of
this paper is {\it elliptic hypergeometric series}, some elliptic
determinant evaluations turn out to be crucial for the proofs of
the results. Lemma~5.3 from
\machSeite{WarnAG}\cite{WarnAG} extends one of the basic
determinant lemmas listed in
\machSeite{KratBN}\cite{KratBN}, namely  
\machSeite{KratBN}\cite[Lemma~5]{KratBN}, to the elliptic world, to
which it reduces in the case $p=0$. We present this important
elliptic determinant evaluation in the theorem below.
 
\begin{Theorem}
\label{bcdet}
Let $x_1,x_2,\dots,x_n$, $a_1,a_2,\dots,a_n$ be indeterminates.
For each $j=1,\dots,n$, let $P_j(x)$ be a $D_j$ theta function. 
Then there holds
\begin{multline}\label{bcdetid}
\det_{1\le i,j\le n}\left(P_{j}(x_i)
\prod_{k=j+1}^n\theta(a_kx_i)\,\theta(a_k/x_i)\right)\\
=\prod_{i=1}^nP_{i}(a_i)
\prod_{1\le i<j\le n}a_jx_j^{-1}\,\theta(x_jx_i)\,\theta(x_j/x_i).
\end{multline}
\quad \quad \qed
\end{Theorem}

Warnaar used this identity to obtain a summation formula for a 
{\it multidimensional elliptic hyper\-geometric
series}. Further related applications may be found in 
\machSeite{RoseAA}%
\machSeite{Ro}%
\machSeite{RoScAB}%
\machSeite{Sp}%
\cite{RoseAA,Ro,RoScAB,Sp}.
The relevant special case of the above theorem is the following
(see \machSeite{WarnAG}\cite[Cor.~5.4]{WarnAG}). 
It is the elliptic generalisation of 
\machSeite{KratBN}\cite[Theorem~28]{KratBN}.
In the
statement, we use the notation
\begin{equation} \label{eq:pqell} 
(a;q,p)_m=\theta(a;p)\,\theta(aq;p)\cdots \theta(aq^{m-1};p),
\end{equation}
which extends the notation for $q$-shifted factorials to the elliptic
world. 

\begin{Theorem} \label{thm:Warn1}
Let $X_1,X_2,\dots,X_n$, $A$, $B$ and $C$ be indeterminates. Then, for
any non-negative integer $n$, there holds
\begin{multline} \label{eq:Warn1} 
\det_{1\le i,j\le n}\(\frac {(AX_i;q,p)_{n-j}\,(AC/X_i;q,p)_{n-j}} 
{(BX_i;q,p)_{n-j}\,(BC/X_i;q,p)_{n-j}}\)\\=(Aq)^{\binom n2}
\prod _{1\le i<j\le n} ^{}X_j\,\theta(X_i/X_j)\,\theta(C/X_iX_j)
\prod _{i=1} ^{n}\frac {(B/A;q,p)_{i-1}\,(ABCq^{2n-2i};q,p)_{i-1}} 
{(BX_i;q,p)_{n-1}\,(BC/X_i;q,p)_{n-1}}.
\end{multline}
\quad \quad \qed
\end{Theorem}

Theorem~29 from \machSeite{KratBN}\cite{KratBN}, which is slightly
more general than \machSeite{KratBN}\cite[Theorem~28]{KratBN}, can
also be extended to an elliptic theorem by suitably specialising
the variables in Theorem~\ref{bcdet}.

\begin{Theorem} \label{thm:Warn1a}
Let $X_1,X_2,\dots,X_n$, $Y_1,Y_2,\dots,Y_n$, $A$ and $B$ be 
indeterminates. Then, for
any non-negative integer $n$, there holds

\vbox{\noindent
\begin{multline} \label{eq:Warn1a} 
\det_{1\le i,j\le n}\(\frac {(X_iY_j;q,p)_{j}\,(AY_j/X_i;q,p)_{j}} 
{(BX_i;q,p)_{j}\,(AB/X_i;q,p)_{j}}\)\\=
q^{2\binom n3}(AB)^{\binom n2}
\prod _{1\le i<j\le n} ^{}\theta(X_jX_i/A)\,\theta(X_j/X_i)\\
\times
\prod _{i=1} ^{n}\frac {(ABY_iq^{i-2};q,p)_{i-1}\,(Y_i/Bq^{i-1};q,p)_{i-1}} 
{X_i^{i-1}\,(BX_i;q,p)_{n-1}\,(AB/X_i;q,p)_{n-1}}.
\end{multline}
\quad \quad \qed}
\end{Theorem}

Another, very elegant, special case of Theorem~\ref{bcdet} 
is the following elliptic Cauchy-type
determinant evaluation.
It was used by Rains
\machSeite{RainAA}\cite[Sec.~3]{RainAA} in the course
of deriving a {\it $BC_n\leftrightarrow BC_m$ integral transformation}.

\begin{Lemma}
\label{frobc}
Let $x_1,x_2,\dots,x_n$ and $a_1,a_2,\dots,a_n$ be indeterminates.
Then there holds
\begin{equation}\label{frobcid}
\det_{1\leq i,j\leq n}\left(\frac 1{\theta(a_jx_i)\,\theta(a_j/x_i)}\right)
=\frac{\displaystyle\prod_{1\leq i<j\leq n}
a_jx_j^{-1}\,\theta(x_jx_i)\,\theta(x_j/x_i)\,\theta(a_ia_j)\,\theta(a_i/a_j)}
{\displaystyle\prod_{i,j=1}^n\theta(a_jx_i)\,\theta(a_j/x_i)}.
\end{equation}
\quad \quad \qed
\end{Lemma}

The remaining determinant evaluations in the current subsection,
with the exception of the very last one, are all due to Rosengren and
Schlosser 
\machSeite{RoScAC}\cite{RoScAC}. The first one is
a further (however non-obvious) consequence of Theorem~\ref{bcdet} 
(see \machSeite{RoScAC}\cite[Cor.~4.3]{RoScAC}).
Two related determinant evaluations,
restricted to the polynomial case, were applied in 
\machSeite{SchlAB}\cite{SchlAB} and 
\machSeite{SchlAF}\cite{SchlAF}
to obtain {\it multidimensional matrix inversions} that played a major
role in the derivation of new {\it summation formulae for multidimensional
basic hypergeometric series}.

\begin{Theorem} \label{thm:cdet}
Let $x_1,x_2,\dots,x_n$, $a_1,a_2,\dots,a_{n+1}$, and
$b$ be indeterminates.
For each $j=1,\dots,n+1$, let $P_j(x)$ be a $D_j$ theta function. 
Then there holds
\begin{multline}
P_{n+1}(b)\det_{1\le i,j\le n}\!\Bigg(P_{j}(x_i)
\prod_{k=j+1}^{n+1}\big(\theta(a_kx_i)\,\theta(a_k/x_i)\big)\\
-\frac{P_{n+1}(x_i)}{P_{n+1}(b)}P_{j}(b)
\prod_{k=j+1}^{n+1}\big(\theta(a_kb)\,\theta(a_k/b)\big)\Bigg)\\
=\prod_{i=1}^{n+1}P_{i}(a_i)
\prod_{1\le i<j\le n+1}a_jx_j^{-1}\,\theta(x_jx_i)\,\theta(x_j/x_i),
\end{multline}
where $x_{n+1}=b$.
\quad \quad \qed
\end{Theorem}

The next determinant evaluation is Theorem~4.4 from
\machSeite{RoScAC}\cite{RoScAC}. It generalises another basic
determinant lemma listed in
\machSeite{KratBN}\cite{KratBN}, namely Lemma~6 from 
\machSeite{KratBN}\cite{KratBN}, to the elliptic case.
It looks as if it is a limit case of
Warnaar's in Theorem~\ref{bcdet}. However, limits are very problematic
in the elliptic world, and therefore it does not seem that
Theorem~\ref{bcdet} implies the theorem below.
For a generalisation in a different direction than Theorem~\ref{bcdet}
see \machSeite{TV}\cite[Appendix B]{TV} (cf.\ also
\machSeite{RoScAC}\cite[Remark~4.6]{RoScAC}).

\begin{Theorem}\label{adet}
Let $x_1,x_2,\dots,x_n$, $a_1,a_2,\dots,a_n$, and $t$ be indeterminates.
For each $j=1,\dots,n$, let $P_j(x)$ be an $A_{j-1}$ theta function of
norm $ta_1\dotsm a_j$. Then there holds
\begin{multline}\label{adetid}
\det_{1\le i,j\le n}\left(P_j(x_i)
\prod_{k=j+1}^n\theta(a_kx_i)\right)\\
=\frac{\theta(ta_1\dotsm a_nx_1\dotsm x_n)}{\theta(t)}
\prod_{i=1}^nP_i(1/a_i)
\prod_{1\le i<j\le n}a_jx_j\,\theta(x_i/x_j).
\end{multline}
\quad \quad \qed
\end{Theorem}

As is shown in \machSeite{RoScAC}\cite[Cor.~4.8]{RoScAC}, this
identity implies the following determinant evaluation.

\begin{Theorem} \label{adetcor}
Let $x_1,x_2,\dots,x_n$, $a_1,a_2,\dots,a_{n+1}$ and $b$ be indeterminates.
For each $j=1,\dots,n+1$, let $P_j(x)$ be an $A_{j-1}$ theta function of
norm $ta_1\dotsm a_j$. Then there holds
\begin{multline}
P_{n+1}(b)\;\det_{1\le i,j\le n}\left(P_j(x_i)
\prod_{k=j+1}^{n+1}\theta(a_kx_i)-\frac{P_{n+1}(x_i)}{P_{n+1}(b)}
P_j(b)\prod_{k=j+1}^{n+1}\theta(a_kb)\right)\\
=\frac{\theta(tba_1\dotsm a_{n+1}x_1\dotsm x_n)}{\theta(t)}
\prod_{i=1}^{n+1}P_i(1/a_i)
\prod_{1\le i<j\le n+1}a_jx_j\,\theta(x_i/x_j),
\end{multline}
where $x_{n+1}=b$.
\quad \quad \qed
\end{Theorem}

By combining Lemma~\ref{frobc} and Theorem~\ref{adet},
a determinant evaluation similar to the one in Theorem~\ref{thm:cdet},
but different, is obtained in
\machSeite{RoScAC}\cite[Theorem~4.9]{RoScAC}.

\begin{Theorem}\label{cdet}
Let $x_1,x_2,\dots,x_n$, $a_1,a_2,\dots,a_n$, and $c_1,\dots, c_{n+2}$ be
indeterminates. For each $j=1,\dots,n$, let $P_j(x)$ be an $A_{j-1}$
theta function of norm
$(c_1\dotsm c_{n+2}a_{j+1}\dotsm a_n)^{-1}$. Then there holds
\begin{multline}\label{dv}
\det_{1\leq i,j\leq n}\left(x_i^{-n-1}
P_j(x_i)\prod_{k=1}^{n+2}\theta(c_kx_i)\,
\prod_{k=j+1}^n\theta(a_kx_i)\right.\\
\left.-x_i^{n+1}P_j(x_i^{-1})\prod_{k=1}^{n+2}\theta(c_kx_i^{-1})\,
\prod_{k=j+1}^n\theta(a_kx_i^{-1})\right)\\
\kern-4cm=\frac{a_1\dotsm a_n}
{x_1\dotsm x_n\,\theta(c_1\dotsm c_{n+2}a_1\dotsm a_n)}
\prod_{i=1}^nP_i(1/a_i)\\\times
\prod_{1\leq i<j\leq n+2}\theta(c_ic_j)\prod_{i=1}^n\theta(x_i^2)
\prod_{1\le i<j\le n}a_jx_i^{-1}\,\theta(x_ix_j)\,\theta(x_i/x_j).
\end{multline}
\quad \quad \qed
\end{Theorem}

The last elliptic determinant evaluation which I present here is a
surprising elliptic extension of a determinant evaluation due to
Andrews and Stanton 
\machSeite{AnStAA}\cite[Theorem~8]{AnStAA}
(see \machSeite{KratBN}\cite[Theorem~42]{KratBN})
due to Warnaar \machSeite{WarnAG}\cite[Theorem~4.17]{WarnAG}.
It is surprising because in the former there appear $q$-shifted
factorials {\it and\/} $q^2$-shifted factorials at the same time, 
but nevertheless there exists
an elliptic analogue, and to obtain it one only has to add the $p$
everywhere in
the shifted factorials to convert them to elliptic ones.

\begin{Theorem} \label{thm:Warn2}
Let $x$ and $y$ be indeterminates. Then, for any non-negative integer
$n$, there holds
\begin{multline} \label{eq:Warn2}
\det_{0\le i,j\le n-1}\(\frac {(y/xq^{i};q^2,p)_{i-j}\,
(q/yxq^i;q^2,p)_{i-j}\, (1/x^2q^{2+4i};q^2,p)_{i-j}} {(q;q,p)_{2i+1-j}\,
(1/yxq^{2i};q,p)_{i-j}\, (y/xq^{1+2i};q,p)_{i-j}}\)\\
=\prod _{i=0} ^{n-1}\frac {(x^2q^{2i+1};q,p)_i\, (xq^{3+i}/y;q^2,p)_i\,
(yxq^{2+i};q^2,p)_i} {(x^2q^{2i+2};q^2,p)_i\, (q;q^2,p)_{i+1}\,
(yxq^{1+i};q,p)_i\, (xq^{2+i}/y;q,p)_i}.
\end{multline}
\quad \quad \qed
\end{Theorem}

In closing this final subsection, I remind the reader that,
as was already said before, many Hankel determinant
evaluations involving elliptic functions can be found in
\machSeite{MilnAO}\cite{MilnAP} and
\machSeite{MilnAP}\cite{MilnAO}.



\section*{Acknowledgments}
I would like to thank Anders Bj\"orner and Richard
Stanley, and the Institut Mittag--Leffler, for giving me the opportunity
to work in a relaxed and inspiring atmosphere during the
``Algebraic Combinatorics" programme in Spring 2005 at the Institut, 
without which this article would never have reached its present form. 
Moreover, I am extremely grateful to Dave Saunders and Zhendong Wan
who performed the {\sl LinBox} computations of determinants of size $3840$,
without which it would have been impossible for me to formulate
Conjectures~\ref{prob:1}--\ref{prob:3} and \ref{prob:8}.
Furthermore I wish to thank Josep Brunat, Adriano Garsia, Greg Kuperberg,
Antonio Montes,
Yuval Roichman, Hjalmar Rosengren, Michael Schlosser, 
Guoce Xin, and especially Alain
Lascoux, for the many useful comments and discussions which helped to
improve the contents of this paper considerably.

\end{document}

